\documentclass[12pt]{article}
\usepackage[usenames]{color}
\usepackage[german,english]{babel}
\usepackage{amsmath}
\usepackage{amssymb,latexsym}
\usepackage{moreverb}
\usepackage{url}
\usepackage{graphicx}
\usepackage{longtable}
\usepackage{rotating}
\usepackage{lscape}
\usepackage{booktabs}
\usepackage{tabularx}
\usepackage{psfrag}
\usepackage{dashrule}
\usepackage{microtype}

\usepackage{multicol}
\usepackage{gastex}
\newtheorem{deff}{Definition :}
\newtheorem{lem}[deff]{\it Lemma :}

\newtheorem{thm}[deff]{Theorem :}
\newtheorem{rem}[deff]{Remark :}

\newtheorem{example}[deff]{Example :}
\newcommand{\proof}{{\bf Proof :}}

\makeatletter

\title{\sc Semi - Equivelar Maps on the Torus and the Klein Bottle with few vertices}
\author{ Anand Kumar Tiwari and Ashish K. Upadhyay \\[1mm]
Department of Mathematics\\Indian Institute of Technology Patna\\
Patliputra Colony, Patna -- 800\,013, India\\
{\small \{upadhyay, anand\}@iitp.ac.in}}

\date{\today}

\begin{document}

\selectlanguage{english}

\maketitle

\hrule

\begin{abstract}

Semi-Equivelar maps are generalizations of maps on the surfaces of Archimedean solids to surfaces other than the $2$-sphere. The well known 11 types of normal tilings of the plane suggest the possible types of semi-equivelar maps on the torus and the Klein bottle. In this article we classify (up to isomorphism) semi-equivelar maps on the torus and the Klein bottle with few vertices.





\end{abstract}

\hrule

\bigskip

{\bf AMS Subject Classification\,:} 52B70, 57M20

{\bf Keywords\,:} Archimedean solids, planar tilings, torus, Klein bottle, semi-equivelar maps.


\section{Introduction and results}

As is well known, equivelar triangulations, also known as degree-regular triangulations (see \cite{dattaupadhyay1}, \cite{dattaupadhyay2}, \cite{LSTU}) of surfaces and in more generality equivelar maps (see \cite{Babai}, \cite{datta}, \cite{grunsheph}, \cite{karabasNedela}) on surfaces are generalizations of the maps on surfaces of Platonic solids to the surfaces other than the 2-sphere. In this article we attempt to study generalizations of maps on the surfaces of Archimedean solids to the torus and the Klein bottle. A similar study for the surfaces of Euler characteristics $-1$, $-2$, $-4$, and $-6$ was carried out in \cite{karabasNedela} and  \cite{upatiwamaity}. We call such objects semi-equivelar Maps or briefly SEM(s).

Most of the definitions  and the terminology used in this section are given in \cite{datta}, \cite{McSchWills} and \cite{upatiwamaity}. We reproduce them
here for the sake of completeness and ready reference. A {\em p-cycle}, denoted $C_p$,  is defined as a finite connected 2-regular graph with $p$ vertices. A 2-dimensional {\em polyhedral complex} $K$ is a collection of $p_{i}$-cycles, where $\{p_{i}\colon 1 \leq i \leq n\}$ is a set of positive integers $\geq 3$, together with vertices and edges in the cycles such that the intersection of any two cycles is empty, a vertex or is an edge. The cycles are called faces of $K$ and the symbols $V(K)$ and $EG(K)$ denote the set of vertices and edges of $K$ respectively. A polyhedral complex $K$ is called a {\em polyhedral 2-manifold} if for each vertex $v$ the faces containing $v$ are of the form $C_{p_1}, \ldots, C_{p_m}$ where $C_{p_1}\cap C_{p_2}, \ldots C_{p_{m-1}}\cap C_{p_m}$,  and  $C_{p_m}\cap C_{p_1}$ are edges for some $m \geq 3$. A connected finite polyhedral 2-manifold is called a {\em polyhedral map}. We will use the term {\em map} for a polyhedral map.  A geometric object $|K|$ is associated to a polyhedral complex $K$ as follows\,: corresponding to each $p$-cycle $C_p$ in $K$, consider a $p$-gon $D_p$ whose boundary cycle is $C_p$. Then $|K|$ is the union of all such $p$-gons and is called the {\em geometric carrier} of $K$. The complex $K$ is said to be connected (resp. orientable) if the topological space $|K|$ is connected (resp. orientable). Between any two polyhedral complexes $K_1$ and $K_2$ we define an isomorphism to be a map $f\colon K_1 \longrightarrow K_2$ such that $f|_{V(K_1)}\colon V(K_1) \longrightarrow V(K_2)$ is a bijection and $f(\sigma)$ is a cell in $K_2$ if and only if $\sigma$ is a cell in $K_1$. If $K_1 = K_2$ then $f$ is called an automorphism of $K_1$. The set of all automorphisms of a polyhedral complex $K$ form a group under the operation composition of maps. This group is called the group of automorphisms of $K$. If this group acts transitively on the set $V(K)$ then the complex is called a {\em vertex-transitive} complex. Some vertex-transitive maps of Euler characteristic 0 have been studied in \cite{Babai}.

The {\em face-sequence} of a vertex $v$ in a map $K$ is a finite sequence $(a^p, b^q, \ldots, m^r)$ of powers of positive integers  $a, b, \ldots, m\geq 3$ and $p, q, \ldots, r \geq 1$, such that through the vertex $v$, $p$ numbers of $C_a$, $q$ numbers of $C_b$, $\ldots$, $r$ numbers of $C_m$ are incidents. $K$ is said to be {\em semi-equivelar} if the face-sequence of each vertex of $K$ is the same up to a cyclic permutation. A SEM with the face-sequence $(a^p, b^q, \ldots, m^r)$, is also called SEM of type $(a^p, b^q, \ldots, m^r)$.

The present work is motivated by an attempt to search for the existence of SEMs on the torus and the Klein bottle. We construct and classify SEMs of all possible types on these surfaces. Due to computational constraints we classify these maps for $\leq 15$ vertices and obtain maps of the types $(3^3, 4^2)$ and $(3^2, 4, 3, 4)$ only. It is evident that in a SEM, as the size of faces grows the number of vertices increases. In this way we have been able to completely classify the SEMs of types $\{3^4, 6\}$, $\{3, 6, 3, 6\}$, $\{3, 4, 6, 4\}$, $\{4, 8^2\}$, $\{3, 12^2\}$ and $\{4, 6, 12\}$ on at most $20$ vertices. We have shown that the SEMs of types $\{3, 6, 3, 6\}$,  $\{3, 12^2\}$ and $\{4, 6, 12\}$ do not exist on $\leq 20$ vertices. We show that\,:




\begin{lem}\label{l1} If $M$ is a semi equivelar map of type $(3^3, 4^2)$ or $(3^2, 4, 3, 4)$ on the torus or Klein bottle with $\leq 15$ vertices then $M$ is isomorphic to one of $T_{1, 10}(3^3, 4^2)$, $T_{1, 12}(3^3, 4^2)$, $T_{2, 12}(3^3, 4^2)$, $T_{3, 12}(3^3, 4^2)$, $T_{1, 14}(3^3, 4^2)$, $T_{2, 14}(3^3, 4^2)$, $K_{1, 10}(3^3, 4^2)$, $K_{1, 12}(3^3, 4^2)$, $K_{2, 12}(3^3, 4^2)$, $K_{1, 14}(3^3, 4^2)$ or $K_{1, 12}(3^2, 4, 3, 4)$, given in Example \ref{ex1}.
\end{lem}




As a corollary we get\,:

\begin{thm}\label{t1} There are exactly 11 SEMs with $\leq 15$ vertices on the surfaces of Euler characteristic 0. Six of these maps are orientable and five are non-orientable.
\end{thm}

By a simple argument using Euler characteristic equations it is easy to see that the SEMs of the remaining types do not exist on $\leq 15$ vertices. Hence we now take the number of vertices to be at most 20. Then, we have\,:


\begin{lem}\label{l2} If $M$ is a SEM of type $(3, 4, 6, 4)$, $(4, 8^2)$, $(3^4, 6)$ or $(3, 6, 3, 6)$ on the torus or Klein bottle with $\leq 20$ vertices then $M$ is isomorphic to one of $T_{1, 18}(3, 4, 6, 4)$ or $K_{1, 18}(3, 4, 6, 4)$, $T_{1, 20}(4, 8^2)$,  $T_{1, 18}(3^4, 6)$, given in Example \ref{ex1}.
\end{lem}






As a corollary to the above Lemma \ref{l2} we have\,:

\begin{thm}\label{t2} There are at least 15 SEMs with $\leq 20$ vertices on the surfaces of Euler characteristic 0. Nine of these are orientable and six are non-orientable.
\end{thm}
\proof\,: From Example \ref{ex1} it is easy to see that the maps denoted by $T$ are all orientable and those by $K$ are non-orientable. The proof now follows from Lemma \ref{l1}, \ref{l2} and Theorem \ref{t1}. \hfill $\Box$

This article is organized in the following manner. In the first section we present examples of SEMs in the tabulated
form and also pictorially. The subsequent section gives proof of the results mentioned above. The proof is a case
by case exhaustive search for possible examples which might exist. Since the method involved in the proof for search of objects corresponding to different types of SEMs and different number of vertices is same and repetitive, we have skipped the unnecessary details and invite the interested reader to construct it for him(her)self or the reader may ask the authors to send the proof.

\section{Examples}

\noindent {\bf Table 1: Tabulated list of SEMs on torus and Klein bottle with at most 20 vertices}
\footnotesize
 \begin{center}
\renewcommand{\arraystretch}{1.23}
\begin{tabular}{|p{1cm}|p{2cm}|p{1cm}|p{1cm}|p{3cm}|p{3.5cm}|}
\hline
$S.No.$&$SEM$-$Type$&$|V|$&$|Maps|$&$Orientable$ $Maps$&$Non$-$orientable$ $Maps$

\\
\hline
~&~&8&0&-&-
\\
\cline{3-6}
~&~&10&2&$T_{1,{10}}(3^3,4^2)$&$K_{1,{10}}(3^3,4^2)$\\
\cline{3-6}
1.&$(3^3,4^2)$&12&5&$T_{1,{12}}(3^3,4^2)$&$K_{1,{12}}(3^3,4^2)$\\
~&~&~&~&$T_{2,{12}}(3^3,4^2)$&$K_{2,{12}}(3^3,4^2)$\\
~&~&~&~&$T_{3,{12}}(3^3,4^2)$&~\\

\cline{3-6}
~&~&14&3&$T_{1,{14}}(3^3,4^2)$&$K_{1,{14}}(3^3,4^2)$\\
~&~&~&~&$T_{2,{14}}(3^3,4^2)$&~\\
\hline

~&~&8&0&-&-\\

\cline{3-6}
2.&~&10&0&-&-\\
\cline{3-6}
~&$(3^3,4,3,4)$&12&1&-&$K_{1,{12}}(3^2,4,3,4)$\\
~&~&~&~&~&~\\

\hline
3.&~&12&0&-&-\\
\cline{3-6}
~&(3, 4, 6, 4)&18&2&$T_{1,{18}}(3, 4, 6, 4)$&$K_{1,{18}}(3,4,6, 4)$\\

\hline
4.&$(4,8^2)$&20&1&$T_{1, 20}(4,8^2)$&~\\
\hline

5.&$(3^4,6)$&12&0&-&-\\
\cline{3-6}
~&~&18&1&$T_{1, 18}(3^4,6)$&-\\
\hline


\end{tabular}
\end{center}
\bigskip

\vspace{0.19in}

\hrule
\bigskip
\noindent \begin{example}\label{ex1} Some semi-equivelar maps on torus and Klein bottle\,: The symbol $T_{i, n}(w, x, y, z)$(resp. $K_{i, n}(w, x, y, z)$ ) denotes that the figure is of the $i^{th}$ map of type $(w, x, y, z)$ on torus (resp. Klein bottle) with $n$ vertices
\begin{center}
\begin{picture}(60,-50)(3,132)
\setlength{\unitlength}{6.5mm}
\drawpolygon(5,0)(4,1)(3,1)(2,1)(1,1.5)(1,3)(0,2)(-1,1.5)(-3,2)(-3,0)(-4,0)(-5,-1)(-4,-2)(-3,-2)(-2,-2)(-1,-2.5)(-1,-4)(0,-3)(1,-2.5)(3,-3)(3,-1)(4,-1)

\drawpolygon(0,1)(1,0)(1,-1)(0,-2)(-1,-1)(-1,0)

\drawline[AHnb=0](0,1)(1,0)

\drawline[AHnb=0](0,1)(1,1.5)

\drawline[AHnb=0](0,1)(0,2)

\drawline[AHnb=0](0,1)(-1,1.5)

\drawline[AHnb=0](0,1)(-1,0)

\drawline[AHnb=0](0,2)(1,1.5)
\drawline[AHnb=0](-1,0)(-1,1.5)

\drawline[AHnb=0](-1,0)(-3,0)
\drawline[AHnb=0](-1,0)(-2,-.5)

\drawline[AHnb=0](-1,0)(-1,-1)
\drawline[AHnb=0](1,0)(1,1.5)

\drawline[AHnb=0](-1,-1)(-2,-2)
\drawline[AHnb=0](-1,-1)(-2,-.5)
\drawline[AHnb=0](-2,-.5)(-3,-2)
\drawline[AHnb=0](-2,-.5)(-2,-2)
\drawline[AHnb=0](-2,-.5)(-3,0)
\drawline[AHnb=0](-3,0)(-1,1.5)
\drawline[AHnb=0](-1,-1)(-1,-2.5)

\drawline[AHnb=0](-1,-2.5)(0,-2)
\drawline[AHnb=0](-1,-2.5)(0,-3)
\drawline[AHnb=0](0,-2)(0,-3)

\drawline[AHnb=0](0,-2)(1,-2.5)
\drawline[AHnb=0](0,-2)(1,-1)
\drawline[AHnb=0](1,-2.5)(0,-3)
\drawline[AHnb=0](1,-2.5)(3,-3)
\drawline[AHnb=0](1,-2.5)(3,-1)

\drawline[AHnb=0](1,-1)(3,-1)
\drawline[AHnb=0](1,-1)(1,-2.5)
\drawline[AHnb=0](2,-.5)(2,1)

\drawline[AHnb=0](2,-.5)(3,1)
\drawline[AHnb=0](2,-.5)(1,-1)
\drawline[AHnb=0](2,-.5)(3,-1)
\drawline[AHnb=0](1,0)(2,1)
\drawline[AHnb=0](1,0)(1,1.5)

\drawline[AHnb=0](1,0)(2,-.5)


\put(5.1,0){\scriptsize $v_7$}
\put(3.9,1.15){\scriptsize $v_{12}$}
\put(2.8,1.2){\scriptsize $v_{11}$}
\put(1.87,1.2){\scriptsize $v_{17}$}
\put(1.1,1.5){\scriptsize $v_6$}
\put(1,3){\scriptsize $v_{13}$}
\put(-.33,2.15){\scriptsize $v_7$}

\put(-1.3,1.68){\scriptsize $v_8$}
\put(-3.2,2.1){\scriptsize $v_{15}$}
\put(-3.5,.15){\scriptsize $v_9$}
\put(-4.3,.15){\scriptsize $v_{14}$}
\put(-5.66,-1){\scriptsize $v_{13}$}

\put(-4.2,-2.3){\scriptsize $v_6$}

\put(-3.23,-2.3){\scriptsize $v_{17}$}
\put(-2.35,-2.3){\scriptsize $v_{11}$}
\put(-1.65,-2.7){\scriptsize $v_{12}$}

\put(-1.3,-4.3){\scriptsize $v_7$}
\put(-.1,-3.3){\scriptsize $v_{13}$}
\put(3,-3.3){\scriptsize $v_9$}
\put(3.1,-1.35){\scriptsize $v_{15}$}
\put(4,-1.3){\scriptsize $v_8$}

\put(-.2,.55){\scriptsize $v_0$}
\put(.5,-.1){\scriptsize $v_5$}
\put(-.86,-.1){\scriptsize $v_1$}
\put(.55,-.9){\scriptsize $v_4$}
\put(-.93,-.9){\scriptsize $v_2$}
\put(-.18,-1.7){\scriptsize $v_3$}

\put(2.2,-.5){\scriptsize $v_{16}$}
\put(-2.8,-.6){\scriptsize $v_{10}$}
\put(.86,-2.85){\scriptsize $v_{14}$}

\put(-4,-3.5){\scriptsize $T_{1, 18}(3^4, 6)$}
\end{picture}

\begin{picture}(0,0)(-14,55)

\setlength{\unitlength}{7mm}

\drawpolygon(-2,1)(-2,-1)(3,-1)(3,1)

\drawline[AHnb=0](-2,0)(3,0)
\drawline[AHnb=0](-1,-1)(-1,1)
\drawline[AHnb=0](0,-1)(0,1)
\drawline[AHnb=0](1,-1)(1,1)
\drawline[AHnb=0](2,-1)(2,1)

\drawline[AHnb=0](-1,0)(-2,1)
\drawline[AHnb=0](0,0)(-1,1)
\drawline[AHnb=0](1,0)(0,1)
\drawline[AHnb=0](2,0)(1,1)
\drawline[AHnb=0](3,0)(2,1)

\put(-.1,1.1){\scriptsize $v_7$}
\put(.9,1.1){\scriptsize $v_1$}
\put(-1.1,1.1){\scriptsize $v_4$}
\put(-.95,-.25){\scriptsize $v_0$}
\put(.05,-.25){\scriptsize $v_3$}
\put(1.05,-.25){\scriptsize $v_8$}
\put(2.05,-.25){\scriptsize $v_9$}
\put(3.1,0){\scriptsize $v_6$}
\put(-2.45,0){\scriptsize $v_6$}
\put(-2.1,-1.3){\scriptsize $v_7$}
\put(-1.1,-1.3){\scriptsize $v_1$}
\put(-.1,-1.3){\scriptsize $v_2$}
\put(.9,-1.3){\scriptsize $v_5$}
\put(1.9,-1.3){\scriptsize $v_4$}
\put(2.9,-1.3){\scriptsize $v_7$}
\put(-2.2,1.1){\scriptsize $v_5$}
\put(1.9,1.2){\scriptsize $v_2$}
\put(2.9,1.1){\scriptsize $v_5$}

\put(-.5,-2){\scriptsize $T_{1, 10}(3^3, 4^2)$}

\end{picture}
\begin{picture}(0,0)(-66,3)
\setlength{\unitlength}{7mm}

\drawpolygon(1,0)(1,1)(-1,1)(-1,0)(-1.5,-1)(-1.5,-2)(-1,-3)(1,-3)(1.5,-2)(1.5,-1)
\drawline[AHnb=0](-1,0)(-1,0)
\drawline[AHnb=0](-1.5,-1)(1.5,-1)
\drawline[AHnb=0](-1.5,-2)(-1.5,-2)

\drawline[AHnb=0](0,0)(0,1)
\drawline[AHnb=0](0,0)(.5,-1)
\drawline[AHnb=0](0,0)(-.5,-1)
\drawline[AHnb=0](-1,0)(1,0)
\drawline[AHnb=0](-.5,-1)(-1,0)
\drawline[AHnb=0](.5,-1)(1,0)

\drawline[AHnb=0](-.5,-1)(-.5,-2)
\drawline[AHnb=0](.5,-1)(.5,-2)

\drawline[AHnb=0](-1.5,-2)(1.5,-2)

\drawline[AHnb=0](-.5,-2)(-1,-3)
\drawline[AHnb=0](-.5,-2)(0,-3)
\drawline[AHnb=0](.5,-2)(0,-3)
\drawline[AHnb=0](.5,-2)(1,-3)

\put(0.1,.1){\scriptsize $v_4$}
\put(1,.06){\scriptsize $v_6$}
\put(0,1.1){\scriptsize $v_9$}
\put(.75,1.1){\scriptsize $v_7$}
\put(-1.15,1.1){\scriptsize $v_8$}
\put(-1.46,0){\scriptsize $v_5$}
\put(-1.98,-.9){\scriptsize $v_6$}
\put(-1.98,-2){\scriptsize $v_7$}
\put(-1.3,-3.2){\scriptsize $v_8$}
\put(0,-3.3){\scriptsize $v_9$}
\put(1.1,-3.1){\scriptsize $v_7$}
\put(1.5,-2){\scriptsize $v_8$}
\put(1.5,-1){\scriptsize $v_5$}

\put(-.9,-1.2){\scriptsize $v_0$}
\put(.6,-1.3){\scriptsize $v_3$}

\put(-.9,-1.8){\scriptsize $v_1$}
\put(.6,-1.9){\scriptsize $v_2$}

\put(-1.1,-4){\scriptsize $K_{1, 10}(3^3, 4^2)$}

\end{picture}

\begin{picture}(0,0)(-54,65)
\setlength{\unitlength}{7mm}
\drawpolygon(0,0)(3,0)(3,4)(0,4)

\drawline[AHnb=0](0,1)(3,1)
\drawline[AHnb=0](0,2)(3,2)
\drawline[AHnb=0](0,3)(3,3)
\drawline[AHnb=0](0,0)(1,1)
\drawline[AHnb=0](1,0)(2,1)
\drawline[AHnb=0](2,0)(3,1)
\drawline[AHnb=0](1,2)(0,3)
\drawline[AHnb=0](2,2)(1,3)
\drawline[AHnb=0](3,2)(2,3)

\drawline[AHnb=0](1,0)(1,4)
\drawline[AHnb=0](2,0)(2,4)

\put(-.7,0){\scriptsize $v_{11}$}
\put(-.6,1){\scriptsize $v_{7}$}
\put(-.7,2){\scriptsize $v_{6}$}
\put(-.7,3){\scriptsize $v_{10}$}
\put(-.7,4){\scriptsize $v_{11}$}

\put(3.15,0){\scriptsize $v_{10}$}
\put(3.15,1){\scriptsize $v_{6}$}
\put(3.15,2){\scriptsize $v_{7}$}
\put(3.15,3){\scriptsize $v_{11}$}
\put(3.15,4){\scriptsize $v_{10}$}

\put(.85,4.2){\scriptsize $v_4$}
\put(1.85,4.2){\scriptsize $v_{5}$}

\put(.9,-.35){\scriptsize $v_4$}
\put(1.9,-.35){\scriptsize $v_5$}

\put(.5,1.2){\scriptsize $v_3$}
\put(2.15,1.2){\scriptsize $v_0$}

\put(.5,1.7){\scriptsize $v_2$}
\put(2.15,1.7){\scriptsize $v_1$}

\put(.5,3.2){\scriptsize $v_9$}
\put(2.15,3.2){\scriptsize $v_8$}

\put(.35,-1){\scriptsize $K_{1, 12}(3^3, 4^2)$}
\end{picture}

\begin{picture}(60,-40)(30,1)
\setlength{\unitlength}{7mm}
\drawpolygon(-1,2)(2,2)(2,-2)(-1,-2)

%
\drawline[AHnb=0](-1,1)(2,1)
\drawline[AHnb=0](-1,0)(2,0)
\drawline[AHnb=0](-1,-1)(2,-1)
\drawline[AHnb=0](0,2)(0,-2)
\drawline[AHnb=0](1,2)(1,-2)

\drawline[AHnb=0](-1,-1)(0,0)
\drawline[AHnb=0](0,-1)(1,0)
\drawline[AHnb=0](1,-1)(2,0)
\drawline[AHnb=0](-1,1)(0,2)
\drawline[AHnb=0](0,1)(1,2)
\drawline[AHnb=0](1,1)(2,2)

\put(2.1,2){\scriptsize $v_{11}$}
\put(-1.6,2){\scriptsize $v_{11}$}
\put(2.1,1){\scriptsize $v_6$}
\put(2.1,0){\scriptsize $v_{7}$}
\put(2.1,-1){\scriptsize $v_8$}
\put(2.1,-2){\scriptsize $v_4$}
\put(-1.5,-2){\scriptsize $v_4$}
\put(-1.5,-1){\scriptsize $v_8$}
\put(-1.5,0){\scriptsize $v_{7}$}
\put(-1.5,1){\scriptsize $v_{6}$}
\put(-.15,2.1){\scriptsize $v_5$}
\put(.85,2.1){\scriptsize $v_4$}
\put(-.15,-2.25){\scriptsize $v_5$}
\put(.85,-2.25){\scriptsize $v_{11}$}
\put(0.1,.8){\scriptsize $v_0$}
\put(1.1,.8){\scriptsize $v_3$}
\put(0.1,.2){\scriptsize $v_1$}
\put(1.1,.2){\scriptsize $v_2$}
\put(0.1,-1.2){\scriptsize $v_{9}$}
\put(1.1,-1.2){\scriptsize $v_{10}$}

\put(-.8,-2.9){\scriptsize $K_{2, 12}(3^3, 4^2)$}
\end{picture}

\begin{picture}(60,-40)(-12,-3)
\setlength{\unitlength}{7mm}
\drawpolygon(-1,2)(2,2)(2,-2)(-1,-2)

%
\drawline[AHnb=0](-1,1)(2,1)
\drawline[AHnb=0](-1,0)(2,0)
\drawline[AHnb=0](-1,-1)(2,-1)
\drawline[AHnb=0](0,2)(0,-2)
\drawline[AHnb=0](1,2)(1,-2)
\drawline[AHnb=0](-1,1)(0,2)
\drawline[AHnb=0](0,1)(1,2)
\drawline[AHnb=0](1,1)(2,2)
\drawline[AHnb=0](-1,-1)(0,0)
\drawline[AHnb=0](0,-1)(1,0)
\drawline[AHnb=0](1,-1)(2,0)

\put(.9,2.1){\scriptsize $v_{9}$}
\put(1.9,2.1){\scriptsize $v_8$}
\put(-.1,2.1){\scriptsize $v_{10}$}
\put(-1.2,2.1){\scriptsize $v_8$}
\put(-1.1,-2.3){\scriptsize $v_9$}
\put(-.1,-2.3){\scriptsize $v_8$}
\put(.9,-2.3){\scriptsize $v_{10}$}
\put(1.9,-2.3){\scriptsize $v_9$}
\put(2.1,-1){\scriptsize $v_{11}$}
\put(2.1,0){\scriptsize $v_6$}
\put(2.1,1){\scriptsize $v_7$}
\put(.05,.8){\scriptsize $v_2$}
\put(1.05,.8){\scriptsize $v_1$}
\put(.05,.1){\scriptsize $v_{3}$}
\put(1.05,.1){\scriptsize $v_0$}
\put(-1.5,1){\scriptsize $v_7$}
\put(-1.5,0){\scriptsize $v_6$}
\put(-1.6,-1){\scriptsize $v_{11}$}
\put(.05,-1.2){\scriptsize $v_4$}
\put(1.05,-1.2){\scriptsize $v_5$}

\put(-.7,-3){\scriptsize $T_{2, 12}(3^3, 4^2)$}
\end{picture}

\begin{picture}(0,0)(-5,90)
\setlength{\unitlength}{7mm}
\drawpolygon(0,0)(3,0)(3,4)(0,4)

\drawline[AHnb=0](1,0)(1,4)
\drawline[AHnb=0](2,0)(2,4)
\drawline[AHnb=0](0,1)(3,1)
\drawline[AHnb=0](0,2)(3,2)
\drawline[AHnb=0](0,3)(3,3)
\drawline[AHnb=0](0,2)(1,1)
\drawline[AHnb=0](1,2)(2,1)
\drawline[AHnb=0](2,2)(3,1)
\drawline[AHnb=0](0,4)(1,3)
\drawline[AHnb=0](1,4)(2,3)
\drawline[AHnb=0](2,4)(3,3)
%
%
\put(-.1,-.3){\scriptsize $v_5$}
\put(.9,-.3){\scriptsize $v_4$}
\put(1.9,-.3){\scriptsize $v_{11}$}
\put(2.8,-.3){\scriptsize $v_7$}

\put(-.2,4.15){\scriptsize $v_5$}
\put(.8,4.15){\scriptsize $v_{4}$}
\put(1.8,4.15){\scriptsize $v_{11}$}
\put(2.8,4.15){\scriptsize $v_7$}
\put(-.6,.9){\scriptsize $v_{10}$}
\put(-.5,1.9){\scriptsize $v_7$}
\put(-.5,2.9){\scriptsize $v_6$}

\put(3.1,.9){\scriptsize $v_{6}$}
\put(3.1,1.9){\scriptsize $v_5$}
\put(3.1,2.9){\scriptsize $v_{10}$}

\put(.56,.75){\scriptsize $v_8$}
\put(1.56,.75){\scriptsize $v_9$}

\put(.56,2.75){\scriptsize $v_0$}
\put(1.56,2.75){\scriptsize $v_{3}$}

\put(.56,1.75){\scriptsize $v_1$}
\put(1.56,1.75){\scriptsize $v_2$}
%
%
\put(.4,-1){\scriptsize $T_{1, 12}(3^3, 4^2)$}
\end{picture}

\begin{picture}(60,-40)(-53,-11)
\setlength{\unitlength}{7mm}
\drawpolygon(2,2)(-1,2)(-1,-2)(2,-2)

\drawline[AHnb=0](-1,1)(2,1)
\drawline[AHnb=0](-1,0)(2,0)
\drawline[AHnb=0](-1,-1)(2,-1)
\drawline[AHnb=0](0,-2)(0,2)
\drawline[AHnb=0](1,-2)(1,2)
\drawline[AHnb=0](-1,1)(0,2)
\drawline[AHnb=0](0,1)(1,2)
\drawline[AHnb=0](1,1)(2,2)
\drawline[AHnb=0](-1,-1)(0,0)
\drawline[AHnb=0](0,-1)(1,0)
\drawline[AHnb=0](1,-1)(2,0)

\put(2.1,0){\scriptsize $v_{6}$}
\put(2.1,1){\scriptsize $v_7$}
\put(1.9,2.1){\scriptsize $v_8$}
\put(.9,2.1){\scriptsize $v_9$}
\put(-.1,2.1){\scriptsize $v_{10}$}
\put(-1.2,2.1){\scriptsize $v_8$}
\put(-1.45,1){\scriptsize $v_{7}$}
\put(-1.45,0){\scriptsize $v_6$}
\put(-1.6,-1){\scriptsize $v_{11}$}
\put(-1.2,-2.3){\scriptsize $v_{10}$}
\put(-.1,-2.3){\scriptsize $v_{9}$}
\put(.9,-2.3){\scriptsize $v_8$}
\put(1.9,-2.3){\scriptsize $v_{10}$}
\put(2.1,-1){\scriptsize $v_{11}$}

\put(.05,.8){\scriptsize $v_2$}
\put(1.05,.8){\scriptsize $v_1$}
\put(.05,-.2){\scriptsize $v_3$}
\put(1.05,-.2){\scriptsize $v_0$}

\put(.05,-1.2){\scriptsize $v_4$}
\put(1.05,-1.2){\scriptsize $v_5$}

\put(-.7,-3){\scriptsize $T_{3, 12}(3^3, 4^2)$}

\end{picture}

\begin{picture}(60,-40)(17,30)
\setlength{\unitlength}{7mm}
\drawpolygon(4,1)(-3,1)(-3,-1)(4,-1)

\drawline[AHnb=0](-3,0)(4,0)
\drawline[AHnb=0](-2,-1)(-2,1)
\drawline[AHnb=0](-1,1)(-1,-1)
\drawline[AHnb=0](0,1)(0,-1)
\drawline[AHnb=0](1,1)(1,-1)
\drawline[AHnb=0](2,1)(2,-1)
\drawline[AHnb=0](3,1)(3,-1)
\drawline[AHnb=0](-2,0)(-3,1)
\drawline[AHnb=0](-1,0)(-2,1)
\drawline[AHnb=0](0,0)(-1,1)
\drawline[AHnb=0](1,0)(0,1)
\drawline[AHnb=0](2,0)(1,1)
\drawline[AHnb=0](3,0)(2,1)
\drawline[AHnb=0](4,0)(3,1)
\put(4.1,0){\scriptsize $v_6$}
\put(3.9,1.1){\scriptsize $v_5$}
\put(2.9,1.1){\scriptsize $v_{12}$}
\put(1.9,1.1){\scriptsize $v_{11}$}
\put(.9,1.1){\scriptsize $v_2$}
\put(-.1,1.1){\scriptsize $v_1$}
\put(-1.1,1.1){\scriptsize $v_7$}
\put(-2.1,1.1){\scriptsize $v_4$}
\put(-3.2,1.1){\scriptsize $v_5$}
\put(-3.5,0){\scriptsize $v_6$}
\put(-3.2,-1.3){\scriptsize $v_7$}
\put(-2.1,-1.3){\scriptsize $v_{1}$}
\put(-1.1,-1.3){\scriptsize $v_2$}
\put(-.1,-1.3){\scriptsize $v_{11}$}
\put(.9,-1.3){\scriptsize $v_{12}$}
\put(1.9,-1.3){\scriptsize $v_{5}$}
\put(2.9,-1.3){\scriptsize $v_4$}
\put(3.9,-1.3){\scriptsize $v_7$}
\put(.05,-.2){\scriptsize $v_8$}
\put(1.05,-.2){\scriptsize $v_{9}$}
\put(2.05,-.2){\scriptsize $v_{10}$}
\put(3.05,-.2){\scriptsize $v_{13}$}
\put(-.95,-.2){\scriptsize $v_3$}
\put(-1.95,-.2){\scriptsize $v_0$}

\put(-.3,-2){\scriptsize $T_{1, 14}(3^3, 4^2)$}
\end{picture}

\begin{picture}(60,-40)(-93,70)
\setlength{\unitlength}{7mm}
\drawpolygon(2,0)(2,1)(1.5,2)(.5,2)(-.5,2)(-1.5,2)(-2,1)(-2,0)(-1.5,-1)(-1.5,-2)(-.5,-2)(.5,-2)(1.5,-2)(1.5,-1)


\drawline[AHnb=0](-1.5,-1)(1.5,-1)

\drawline[AHnb=0](-2,0)(2,0)
\drawline[AHnb=0](-2,1)(2,1)
\drawline[AHnb=0](-.5,-2)(-.5,-1)

\drawline[AHnb=0](.5,-2)(.5,-1)
\drawline[AHnb=0](0,0)(0,1)
\drawline[AHnb=0](0,0)(.5,-1)

\drawline[AHnb=0](0,0)(-.5,-1)

\drawline[AHnb=0](1,0)(1,1)
\drawline[AHnb=0](1,0)(1.5,-1)
\drawline[AHnb=0](1,0)(.5,-1)

\drawline[AHnb=0](-1,0)(-1,1)
\drawline[AHnb=0](-1,0)(-.5,-1)
\drawline[AHnb=0](-1,0)(-1.5,-1)

\drawline[AHnb=0](0,1)(.5,2)
\drawline[AHnb=0](0,1)(-.5,2)

\drawline[AHnb=0](1,1)(1.5,2)

\drawline[AHnb=0](1,1)(.5,2)

\drawline[AHnb=0](-1,1)(-.5,2)
\drawline[AHnb=0](-1,1)(-1.5,2)

\put(2.1,0){\scriptsize $v_8$}
\put(2.1,1){\scriptsize $v_5$}

\put(1.6,2){\scriptsize $v_6$}
\put(.5,2.2){\scriptsize $v_{13}$}
\put(-.5,2.2){\scriptsize $v_4$}

\put(-2,2){\scriptsize $v_5$}
\put(-2.5,1){\scriptsize $v_6$}
\put(-2.5,0){\scriptsize $v_7$}
\put(-2,-1){\scriptsize $v_8$}
\put(-2,-2.2){\scriptsize $v_5$}
\put(-.5,-2.25){\scriptsize $v_4$}

\put(.5,-2.25){\scriptsize $v_{13}$}
\put(1.5,-2.2){\scriptsize $v_6$}
\put(1.7,-1){\scriptsize $v_7$}

\put(-.4,-1.3){\scriptsize $v_9$}
\put(.6,-1.3){\scriptsize $v_{10}$}

\put(.1,.2){\scriptsize $v_2$}
\put(1.1,.2){\scriptsize $v_{11}$}

\put(-.9,.2){\scriptsize $v_1$}

\put(.1,.7){\scriptsize $v_3$}
\put(1.1,.7){\scriptsize $v_{12}$}

\put(-.9,.7){\scriptsize $v_0$}

\put(-1.3,-2.8){\scriptsize $K_{1, 14}(3^3, 4^2)$}
\end{picture}

\begin{picture}(60,-40)(17,55)
\setlength{\unitlength}{7mm}
\drawpolygon(4,1)(-3,1)(-3,-1)(4,-1)
\drawline[AHnb=0](-3,0)(4,0)
\drawline[AHnb=0](-2,1)(-2,-1)
\drawline[AHnb=0](-1,1)(-1,-1)
\drawline[AHnb=0](0,1)(0,-1)
\drawline[AHnb=0](1,1)(1,-1)
\drawline[AHnb=0](2,1)(2,-1)
\drawline[AHnb=0](3,1)(3,-1)

\drawline[AHnb=0](-2,0)(-3,1)
\drawline[AHnb=0](-1,0)(-2,1)
\drawline[AHnb=0](0,0)(-1,1)
\drawline[AHnb=0](1,0)(0,1)
\drawline[AHnb=0](2,0)(1,1)
\drawline[AHnb=0](3,0)(2,1)
\drawline[AHnb=0](4,0)(3,1)

\put(4.1,0){\scriptsize $v_6$}
\put(3.9,1.1){\scriptsize $v_5$}
\put(2.9,1.1){\scriptsize $v_{11}$}
\put(1.9,1.1){\scriptsize $v_2$}
\put(.9,1.1){\scriptsize $v_1$}
\put(-.1,1.1){\scriptsize $v_7$}
\put(-1.1,1.1){\scriptsize $v_{13}$}
\put(-2.1,1.1){\scriptsize $v_4$}
\put(-3.2,1.1){\scriptsize $v_5$}
\put(-3.45,0){\scriptsize $v_6$}
\put(-3.2,-1.3){\scriptsize $v_7$}
\put(-2.1,-1.3){\scriptsize $v_1$}
\put(-1.1,-1.3){\scriptsize $v_2$}
\put(-.1,-1.3){\scriptsize $v_{11}$}
\put(.9,-1.3){\scriptsize $v_5$}
\put(1.9,-1.3){\scriptsize $v_{4}$}
\put(2.9,-1.3){\scriptsize $v_{13}$}
\put(3.9,-1.3){\scriptsize $v_7$}
\put(-1.95,-.3){\scriptsize $v_0$}
\put(-.95,-.3){\scriptsize $v_3$}
\put(.05,-.3){\scriptsize $v_{12}$}
\put(1.05,-.3){\scriptsize $v_{8}$}
\put(2.05,-.3){\scriptsize $v_9$}
\put(3.1,-.3){\scriptsize $v_{10}$}

\put(-.7,-2.2){\scriptsize $T_{2, 14}(3^3, 4^2)$}
\end{picture}


\begin{picture}(60,-40)(-32,132)
\setlength{\unitlength}{7mm}
\drawpolygon(2,1)(2,-2)(-2,-2)(-2,1)

\drawline[AHnb=0](-2,0)(2,0)
\drawline[AHnb=0](-2,-1)(2,-1)
\drawline[AHnb=0](-1,1)(-1,-2)
\drawline[AHnb=0](0,1)(0,-2)
\drawline[AHnb=0](1,1)(1,-2)
\drawline[AHnb=0](0,0)(-1,1)
\drawline[AHnb=0](2,0)(1,1)
\drawline[AHnb=0](-2,-1)(-1,0)
\drawline[AHnb=0](0,-1)(1,0)
\drawline[AHnb=0](0,-2)(-1,-1)
\drawline[AHnb=0](2,-2)(1,-1)

\put(1.9,1.2){\scriptsize $v_{6}$}
\put(.9,1.2){\scriptsize $v_9$}
\put(-.1,1.2){\scriptsize $v_8$}
\put(-1.1,1.2){\scriptsize $v_5$}
\put(-2.2,1.2){\scriptsize $v_6$}
\put(-2.5,0){\scriptsize $v_7$}
\put(-2.5,-1){\scriptsize $v_{1}$}
\put(2.1,0){\scriptsize $v_{7}$}
\put(2.1,-1){\scriptsize $v_1$}
\put(-1.1,-2.3){\scriptsize $v_9$}
\put(-.1,-2.3){\scriptsize $v_8$}
\put(.9,-2.3){\scriptsize $v_5$}
\put(2,-2.3){\scriptsize $v_6$}
\put(-2.2,-2.25){\scriptsize $v_6$}
\put(1.05,-.9){\scriptsize $v_{10}$}
\put(1,.1){\scriptsize $v_{11}$}
\put(.05,.1){\scriptsize $v_4$}
\put(.05,-1.22){\scriptsize $v_3$}
\put(-.98,.1){\scriptsize $v_0$}
\put(-1.5,-1.3){\scriptsize $v_2$}

\put(-1.5,-3){\scriptsize $K_{1, 12}(3^2, 4, 3, 4)$}
\end{picture}


\begin{picture}(60,-40)(-67,95)
\setlength{\unitlength}{7mm}
\drawpolygon(2,-.5)(2.5,0)(3.5,0)(4,.5)(4,1.5)(3.5,2)(2.5,2)(2,1.5)(1.5,2)(.5,2)(0,2.5)(-.5,2)(-1.5,2)(-2,1.5)(-2,.5)(-2.5,0)(-2,-.5)
(-2,-1.5)(-1.5,-2)(-.5,-2)(0,-2.5)(.5,-2)(1.5,-2)(2,-1.5)

\drawpolygon(.5,0)(0,.5)(-.5,0)(0,-.5)

\drawline[AHnb=0](.5,0)(1.5,0)
\drawline[AHnb=0](0,.5)(0,1.5)
\drawline[AHnb=0](-.5,0)(-1.5,0)
\drawline[AHnb=0](0,-.5)(0,-1.5)

\drawline[AHnb=0](0,1.5)(.5,2)
\drawline[AHnb=0](0,1.5)(-.5,2)

\drawline[AHnb=0](-1.5,0)(-2,.5)
\drawline[AHnb=0](-1.5,0)(-2,-.5)

\drawline[AHnb=0](0,-1.5)(-.5,-2)
\drawline[AHnb=0](0,-1.5)(.5,-2)

\drawline[AHnb=0](1.5,0)(2,-.5)
\drawline[AHnb=0](1.5,0)(2,.5)

\drawline[AHnb=0](2,.5)(2,1.5)
\drawline[AHnb=0](2,.5)(2.5,0)

\put(-.5,.5){\scriptsize $v_2$}
\put(.35,.2){\scriptsize $v_{15}$}
\put(1.05,.2){\scriptsize $v_{16}$}
\put(2.1,1.3){\scriptsize $v_{18}$}
\put(2.1,.5){\scriptsize $v_{17}$}
\put(1.2,2.1){\scriptsize $v_{12}$}
\put(.5,2.1){\scriptsize $v_{11}$}
\put(-.5,1.4){\scriptsize $v_3$}

\put(-.85,.1){\scriptsize $v_1$}
\put(-1.5,.1){\scriptsize $v_0$}
\put(-2.5,1.4){\scriptsize $v_6$}
\put(-2.5,.6){\scriptsize $v_7$}
\put(-1.6,2.1){\scriptsize $v_5$}
\put(-.9,2.1){\scriptsize $v_4$}

\put(-.6,-.6){\scriptsize $v_{14}$}
\put(-.65,-1.5){\scriptsize $v_{13}$}
\put(-2.65,-1.5){\scriptsize $v_{10}$}
\put(-2.5,-.7){\scriptsize $v_9$}
\put(-1.7,-2.3){\scriptsize $v_{11}$}
\put(-1,-2.3){\scriptsize $v_{12}$}

\put(2.1,-1.5){\scriptsize $v_7$}
\put(2.1,-.7){\scriptsize $v_6$}
\put(1.3,-2.3){\scriptsize $v_8$}
\put(.5,-2.3){\scriptsize $v_{19}$}

\put(2.5,-.3){\scriptsize $v_5$}
\put(3.2,-.3){\scriptsize $v_4$}
\put(4.1,1.4){\scriptsize $v_9$}
\put(4.1,.35){\scriptsize $v_{10}$}
\put(3.4,2.1){\scriptsize $v_8$}
\put(2.2,2.1){\scriptsize $v_{19}$}

\put(-.2,2.55){\scriptsize $v_{10}$}
\put(-3,0){\scriptsize $v_8$}
\put(-.2,-2.7){\scriptsize $v_{18}$}

\put(2.8,-1.5){\scriptsize $T_{1, 20}(4, 8^2)$}
\end{picture}

\begin{picture}(60,-40)(17,128)
\setlength{\unitlength}{9mm}
\drawpolygon(2,0)(1.7,.5)(1.7,1)(1,1.5)(.5,1)(-.5,1)(-1,1.5)(-1.7,1)(-1.7,.5)(-2,0)(-1.7,-.5)(-1.7,-1)(-1.3,-1.5)(-.5,-2)(.5,-2)
(1.3,-1.5)(1.7,-1)(1.7,-.5)

\drawpolygon(1,0)(.5,.5)(-.5,.5)(-1,0)(-1,-.5)(-1,-1)(0,-1.4)(1,-1)(1,-.5)

\drawline[AHnb=0](0,0)(1,-.5)
\drawline[AHnb=0](0,0)(.5,.5)
\drawline[AHnb=0](0,0)(-.5,.5)
\drawline[AHnb=0](0,0)(-1,-.5)

\drawline[AHnb=0](.5,.5)(.5,1)
\drawline[AHnb=0](-.5,.5)(-.5,1)

\drawline[AHnb=0](1,0)(1.7,.5)
\drawline[AHnb=0](-1,0)(-1.7,.5)

\drawline[AHnb=0](1,0)(1.7,-.5)
\drawline[AHnb=0](-1,0)(-1.7,-.5)

\drawline[AHnb=0](1,-.5)(1.7,-.5)
\drawline[AHnb=0](-1,-.5)(-1.7,-.5)

\drawline[AHnb=0](1,-1)(1.7,-1)
\drawline[AHnb=0](-1,-1)(-1.7,-1)

\drawline[AHnb=0](1,-1)(1.3,-1.5)
\drawline[AHnb=0](-1,-1)(-1.3,-1.5)

\drawline[AHnb=0](1.3,-1.5)(1.7,-1)
\drawline[AHnb=0](-1.3,-1.5)(-1.7,-1)

\drawline[AHnb=0](0,-1.4)(.5,-2)
\drawline[AHnb=0](0,-1.4)(-.5,-2)

\put(-.1,-.3){\scriptsize $v_1$}
\put(2.05,0){\scriptsize $v_{14}$}
\put(1.75,.5){\scriptsize $v_{13}$}
\put(1.75,1){\scriptsize $v_{12}$}
\put(1,1.6){\scriptsize $v_{16}$}
\put(.15,1.1){\scriptsize $v_{17}$}
\put(-.5,1.1){\scriptsize $v_{6}$}

\put(-2.5,0){\scriptsize $v_{16}$}
\put(-2.2,.5){\scriptsize $v_{15}$}
\put(-2.2,1){\scriptsize $v_{14}$}
\put(-1.1,1.6){\scriptsize $v_{7}$}

\put(-2.2,-.6){\scriptsize $v_{12}$}
\put(-2.2,-1.1){\scriptsize $v_{13}$}

\put(-1.7,-1.7){\scriptsize $v_{14}$}
\put(-.7,-2.2){\scriptsize $v_{15}$}
\put(.5,-2.2){\scriptsize $v_{16}$}
\put(1.2,-1.8){\scriptsize $v_{17}$}

\put(1.78,-1.1){\scriptsize $v_{6}$}
\put(1.78,-.6){\scriptsize $v_7$}

\put(.6,.55){\scriptsize $v_9$}
\put(-1,.6){\scriptsize $v_{10}$}
\put(-1.25,.2){\scriptsize $v_{11}$}

\put(-.9,-1){\scriptsize $v_3$}
\put(-.1,-1.25){\scriptsize $v_4$}
\put(.6,-.95){\scriptsize $v_5$}
\put(.6,-.6){\scriptsize $v_0$}

\put(.9,.2){\scriptsize $v_8$}

\put(-.9,-.6){\scriptsize $v_2$}

\put(-1.2,-2.7){\scriptsize $K_{1, 18} (3,4,6,4)$}
\end{picture}

\begin{picture}(-10,0)(-48,123)
\setlength{\unitlength}{9mm}

\drawpolygon(2,0)(1.7,.5)(1.7,1)(1,1.5)(.5,1)(-.5,1)(-1,1.5)(-1.7,1)(-1.7,.5)(-2,0)(-1.7,-.5)(-1.7,-1)(-1.3,-1.5)(-.5,-2)(.5,-2)
(1.3,-1.5)(1.7,-1)(1.7,-.5)

\drawpolygon(1,0)(.5,.5)(-.5,.5)(-1,0)(-1,-.5)(-1,-1)(0,-1.4)(1,-1)(1,-.5)

\drawline[AHnb=0](0,0)(1,-.5)
\drawline[AHnb=0](0,0)(.5,.5)
\drawline[AHnb=0](0,0)(-.5,.5)
\drawline[AHnb=0](0,0)(-1,-.5)

\drawline[AHnb=0](.5,.5)(.5,1)
\drawline[AHnb=0](-.5,.5)(-.5,1)

\drawline[AHnb=0](1,0)(1.7,.5)
\drawline[AHnb=0](-1,0)(-1.7,.5)

\drawline[AHnb=0](1,0)(1.7,-.5)
\drawline[AHnb=0](-1,0)(-1.7,-.5)

\drawline[AHnb=0](1,-.5)(1.7,-.5)
\drawline[AHnb=0](-1,-.5)(-1.7,-.5)

\drawline[AHnb=0](1,-1)(1.7,-1)
\drawline[AHnb=0](-1,-1)(-1.7,-1)

\drawline[AHnb=0](1,-1)(1.3,-1.5)
\drawline[AHnb=0](-1,-1)(-1.3,-1.5)

\drawline[AHnb=0](1.3,-1.5)(1.7,-1)
\drawline[AHnb=0](-1.3,-1.5)(-1.7,-1)

\drawline[AHnb=0](0,-1.4)(.5,-2)
\drawline[AHnb=0](0,-1.4)(-.5,-2)

\put(-.1,-.3){\scriptsize $v_0$}
\put(2.05,0){\scriptsize $v_{12}$}
\put(1.75,.5){\scriptsize $v_{11}$}
\put(1.8,1){\scriptsize $v_{10}$}
\put(1,1.6){\scriptsize $v_{15}$}
\put(.1,1.1){\scriptsize $v_{14}$}
\put(-.6,1.1){\scriptsize $v_{13}$}

\put(-2.5,0){\scriptsize $v_{15}$}
\put(-2.2,.5){\scriptsize $v_{16}$}
\put(-2.2,1){\scriptsize $v_{17}$}
\put(-1.1,1.6){\scriptsize $v_{12}$}

\put(-2.2,-.6){\scriptsize $v_{10}$}
\put(-2.2,-1.1){\scriptsize $v_{11}$}

\put(-1.4,-1.8){\scriptsize $v_{12}$}
\put(-.7,-2.2){\scriptsize $v_{13}$}
\put(.5,-2.2){\scriptsize $v_{14}$}
\put(1.2,-1.8){\scriptsize $v_{15}$}

\put(1.78,-1.1){\scriptsize $v_{16}$}
\put(1.78,-.6){\scriptsize $v_{17}$}

\put(.6,.55){\scriptsize $v_{7}$}
\put(-.9,.58){\scriptsize $v_{8}$}
\put(-1.2,.15){\scriptsize $v_9$}

\put(-.9,-1){\scriptsize $v_2$}
\put(-.1,-1.2){\scriptsize $v_3$}
\put(.5,-1){\scriptsize $v_4$}
\put(.5,-.6){\scriptsize $v_5$}

\put(.9,.2){\scriptsize $v_6$}

\put(-.9,-.6){\scriptsize $v_1$}

\put(-1,-2.7){\scriptsize $T_{1, 18} (3, 4, 6, 4)$}
\end{picture}

\end{center}
\end{example}

\newpage



\vspace{2cm}

\section{Proofs}\label{proof}

In the preceding section all possible face-sequences of semi-equivelar maps on the plane were given. Our aim is to find the SEMs on at most twenty vertices. Consider the case for the SEMs of face sequence $\{3, 12^2\}$. We
observe that link of a vertex in this map contains more than twenty vertices. Hence we discard this case. Similarly considering links of $0$ and $13$ one can see that link of vertex $13$ in the SEM of type $\{4, 6, 12\}$ requires more than 20 vertices. So this case is not possible in the present consideration. Hence we discard this case as well.



\noindent \proof\{Proof of Lemma \ref{l1}\}:\,\, Let $M$ be a SEM of type $(3^3, 4^2)$. The notation ${\rm lk}(i) = C_7(i_1, {\boldsymbol i_2},  i_3,  i_4,  i_5,  i_6,  {\boldsymbol i_7})$ for the link of a vertex $i$ will mean that  $[i,  i_3,  i_4]$,  $[i,  i_4,  i_5]$,  $[i,  i_5,  i_6]$ are triangular faces and $[i,  i_{1},  i_{2},  i_{3}]$ and  $[i,  i_{6},  i_{7},  i_{1}]$ are quadrangular faces. If $|V|$ denotes the number of vertices in $V(M)$, $Q$ denotes the number of quadrangular faces and $T$ denotes number of triangular faces in map $M$, respectively, then it is easy to see that $Q = \frac {2|V|}{4}$ and $T = \frac {3|V|}{3}$. From the expressions of $Q$, $T$ and ${\rm lk}(i)$, we see that a map exists if $|V|$ is an even positive integer and $|V| \geq 8$. Thus for $|V| \leq 15$, the possible values of $|V|$ are $8$,  $10$,  $12$ and  $14$. Without loss of generality,  we may assume ${\rm lk}(0) = C_7(1, {\bf 2},  3,  4,  5,  6,  {\bf 7})$. This implies ${\rm lk}(1) = C_7(0,  {\bf 3},  2,  y,  x,  7,  {\bf 6})$ for some $x,  y \in V$.

When $|V| = 8$,  $i.e.$,  $V = \{0, 1, \ldots, 7\}$ then $(x,  y)\in\{$(4,  5), (5,  4)$\}$. But,
for both the cases,  either 4 or 5 appears in more than three triangles. This is not possible. When $|V| = 10$,  $i.e.$,  $V = \{0,  1,  \ldots,  9\}$ then it is easy to see that $(x,  y)\in\{$(4, 5), (4, 8), (5, 4), (5, 8), (8, 4), (8, 5), (8,  9)$\}$. If $(x,  y)\in\{$(4, 5), (4, 8), (5, 4), (5, 8), (8, 4), (8, 5)$\}$ then as in the previous case, either 4 or 5 appear in more than three triangles, which is not possible. Hence $(x,  y) = (8,  9)$. This implies ${\rm lk}(1) = C_7(0,  {\bf 3},  2,   9,  8,  7,  {\bf 6})$. Then ${\rm lk}(2) = C_7(3,  {\bf 0},  1,  9,  a,  b,  {\bf c})$ for some $a, b, c \in V$. It is easy to see that $(a, b, c)\in\{$(4, 5, 8), (4, 6, 7), (4, 7, 6), (4, 8, 5), (5, 4, 8), (5, 6, 7), (5, 8, 4), (6, 4, 5), (6, 5, 4), (6, 5, 8), (6, 4, 8), (6, 8, 4),
(6, 8, 5), (7, 4, 5), (7, 4, 8), (7, 5, 4), (7, 5, 8), (7, 8, 4), (7, 8, 5)$\}$. But, when $(a, b, c)\in\{$(4, 5, 8), (4, 6, 7), (4, 7, 6), (4, 8, 5)$\}$ then ${\rm lk}(4)$ has more than seven vertices, which is a contradiction. When $(a,  b,  c)\in \{$(5, 4, 8), (5, 6, 7), (5, 8, 4)$\}$, then ${\rm lk}(5)$ has more than seven vertices, which is a contradiction. When $(a,  b,  c) \in \{$(6, 4, 5), (6, 4, 8), (6, 8, 4), (6, 8, 5)$\}$,  then ${\rm lk}(6)$ has more than seven vertices, which is a contradiction. When $(a, b, c) \in \{$(7, 4, 5), (7, 4, 8), (7, 5, 4), (7, 5, 8)$\}$,  then ${\rm lk}(7)$ has more than seven
vertices, which is a contradiction. If $(a, b, c) = (6, 5, 4)$ then we have ${\rm lk}(2) = C_7(3, {\bf 0}, 1, 9, 6, 5, {\bf 4})$, this implies $C_4(0, 3, 2, 5)\subseteq{\rm lk}(4)$, which is a contradiction. When $(a, b, c) = (7, 8, 4)$ then we have ${\rm lk}(2) = C_7(3, {\bf 0}, 1, 9, 7, 8, {\bf 4})$, this implies $C_5(0, 1, 2, 8, 4) \subseteq {\rm lk}(3)$, which is a contradiction. Hence we have $(a, b, c) \in \{$(6, 5, 8), (7, 8, 5)$\}$.

If $(a, b, c) = (6, 5, 8)$ then ${\rm lk}(2) = C_7(3, {\bf 0}, 1, 9, 6, 5, {\bf 8})$.
This implies ${\rm lk}(3) = C_7(2, {\bf 1}, 0, 4, 7, 8, {\bf 5})$, ${\rm lk}(8) = C_7(5, {\bf 2}, 3, 7, 1, 9, {\bf 4})$, completing successively, we get ${\rm lk}(9) = C_7(4, {\bf 5}, 8, 1, 2, 6, {\bf 7})$, ${\rm lk}(7) = C_7(6, {\bf 9}, 4 , 3, 8, 1, {\bf 0})$, ${\rm lk}(6) = C_7(7, {\bf 4}, 9, 2, 5, 0, {\bf 1})$,  ${\rm lk}(4) = C_7(9, {\bf 8}, 5, 0, 3, 7, {\bf 6})$ and ${\rm lk}(5) = C_7(8, {\bf 3}, 2, 6, 0, 4, {\bf 9})$. Then we get $M \cong T_{1, 10}(3^3, 4^2)$  by the map $i \mapsto v_i$, $0\leq i \leq 9$. If $(a, b, c) = (7, 8, 5)$,  then ${\rm lk}(2) = C_7(3, {\bf 0}, 1, 9, 7, 8, {\bf 5})$. Completing successively, we get ${\rm lk}(3) = C_7(2, {\bf 1}, 0, 4, 6, 5, {\bf 8})$, ${\rm lk}(6) = C_7(7, {\bf 1}, 0, 5, 3, 4, {\bf 9})$, ${\rm lk}(7) = C_7(6, {\bf 0}, 1, 8, 2, 9, {\bf 4})$, ${\rm lk}(9) = C_7(4, {\bf 6}, 7, 2, 1, 8, {\bf 5})$, ${\rm lk}(4) = C_7(9, {\bf 7}, 6, 3 , 0, 5, {\bf 8})$, ${\rm lk}(5) = C_7(8, {\bf 9}, 4, 0, 6, 3, {\bf 2})$ and ${\rm lk}(8) = C_7(5, {\bf 4}, 9, 1, 7, 2, {\bf 3})$. Then we get $M \cong K_{1, 10}(3^3, 4^2)$ by the map $i \mapsto v_i$, $0\leq i \leq 9$.

When $|V| = 12$, $i.e.$, $V = \{0,  1,  \ldots,  11\}$ then $(x,  y) = (8,  9)$. So ${\rm lk}(1) = C_7(0, {\bf 3}, 2, 9, 8, 7, {\bf 6})$. This implies ${\rm lk}(2) = C_7(3, {\bf 0}, 1, 9, a, b, {\bf c})$ for some $a, b, c \in V$.
It is easy to see that $(a, b, c) \in \{$(6, 5, 8), (6, 5, 10), (7, 8, 5), (10, 6, 7), (10, 7, 6), (10, 11, 5), (10, 11, 8)$\}$. Now, proceeding as in the case of $|V| = 10$ we get - for $(a, b, c) = (6, 5, 10)$, $M \cong T_{1, 12}(3^3, 4^2)$ by the map $i \mapsto v_i$, $0 \leq i \leq 11$. For $(a, b, c) = (10, 6, 7)$ we get $M \cong K_{1, 12}(3^3, 4^2)$ by the map $i \mapsto v_i$, $0 \leq i \leq 11$. For $(a, b, c) = (10, 7, 6)$ we get $M\cong K_{2, 12}(3^3, 4^2)$ by the map $i \mapsto v_i$,  $0 \leq i \leq 11$, $M \cong T_{2, 12}(3^3, 4^2)$ by the map $i \mapsto v_i$, $0 \leq i \leq 11$, $M \cong T_{3, 12}(3^3, 4^2)$ by the map $i \mapsto v_i$, $0 \leq i \leq 11$, $M \cong K_{2, 12}(3^3, 4^2)$ by the map $0 \mapsto v_8$,  $1 \mapsto v_4$,  $2 \mapsto v_{5}$,  $3\mapsto v_{9}$,  $4 \mapsto v_1$,  $5 \mapsto v_7$,  $6 \mapsto v_{10}$,  $7 \mapsto v_{11}$, $8\mapsto v_3$,  $9\mapsto v_0$,  $10\mapsto v_{6}$,  $11\mapsto v_{2}$, $M \cong K_{2, 12}(3^3, 4^2)$ by the map $0 \mapsto v_8$,  $1 \mapsto v_4$,  $2 \mapsto v_{11}$,  $3\mapsto v_{10}$,  $4 \mapsto v_7$,  $5 \mapsto v_1$,  $6 \mapsto v_9$,  $7 \mapsto v_5$, $8\mapsto v_0$,  $9\mapsto v_3$,  $10\mapsto v_{6}$,  $11\mapsto v_{2}$ and $M \cong T_{2, 12}(3^3, 4^2)$ by the map $0 \mapsto v_2$,  $1 \mapsto v_3$,  $2 \mapsto v_6$,  $3\mapsto v_7$,  $4 \mapsto v_{10}$,  $5 \mapsto v_9$,  $6 \mapsto v_1$,  $7 \mapsto v_0$, $8\mapsto v_4$,  $9\mapsto v_{11}$,  $10\mapsto v_{5}$,  $11\mapsto v_{8}$. For $(a, b, c) = (10, 11, 8)$ we get $M \cong T_{1, 12}(3^3, 4^2)$ by the map $0 \mapsto v_0$,  $1 \mapsto v_1$,  $2 \mapsto v_7$,  $3\mapsto v_6$,  $4 \mapsto v_5$,  $5 \mapsto v_4$,  $6 \mapsto v_3$,  $7 \mapsto v_2$, $8\mapsto v_9$,  $9\mapsto v_8$,  $10\mapsto v_{10}$,  $11\mapsto v_{11}$. For other values of $(a, b, c)$ we do not get any object.

When $|V| = 14$, $i.e.$, $V = \{0, 1, \ldots, 13\}$ then, for $(x, y) = (8, 9)$, we have
${\rm lk}(1) = C_7(0, {\bf 3}, 2, 9, 8, 7, {\bf 6})$. This implies ${\rm lk}(2) = C_7(3, {\bf 0}, 1, 9, a, b, {\bf c})$, for some $a, b, c \in V$. It is easy to see that $(a, b, c) \in \{$(6, 5, 8), (6, 5, 10), (7, 8, 10), (10, 6, 7), (10, 7, 6), (10, 11, 5), (10, 11, 8), (10, 11, 12)$\}$. For $(a, b, c) = (6, 5, 10)$ we get $M \cong T_{1, 14}(3^3, 4^2)$ by the map $ 0 \mapsto v_9,  1 \mapsto v_{12},  2 \mapsto v_{11},  3\mapsto v_8,  4 \mapsto v_1,  5 \mapsto v_2,  6 \mapsto v_{10},  7 \mapsto v_5, 8\mapsto v_6,  9\mapsto v_{13},  10 \mapsto v_{3},  11\mapsto v_{7},  12\mapsto v_{0},  13\mapsto v_{4}$. For $(a, b, c) = (10, 11, 8)$ we get $M \cong T_{1, 14}(3^3, 4^2)$ by the map $i \mapsto v_i, 0\leq i \leq 13$. For $(a, b, c) = (10, 11, 12)$ we get $M\cong T_{2, 14}(3^3, 4^2)$ by the map $i \mapsto v_i$ and $M \cong K_{1, 14}(3^3, 4^2)$ by the map $i \mapsto v_i$. For other values of $(a, b, c)$ we do not get any object.

Let $M$ be a SEM of type $(3^2, 4, 3, 4)$ on a closed surface of Euler characteristic 0. Link of a vertex $i$ of the map is denoted as ${\rm lk}(i) = C_7(i_1, i_2, {\boldsymbol i_3}, i_4, i_5, {\boldsymbol i_6}, i_7)$. The notation for ${\rm lk}(i)$, will mean that $[i, i_1, i_2]$, $[i, i_4, i_5]$, $[i, i_1, i_7]$ are triangular faces and $[i, i_{2}, i_{3}, i_{4}]$ and $[i, i_{5}, i_{6}, i_{7}]$ are quadrangular faces. Let $|V|$ denote the number of vertices in $V(M)$, $Q$ denote the number of quadrangular faces and $T$ denote number of triangular faces in map $M$, respectively. Then $Q = \frac {2|V|}{4}$ and $T = \frac {3|V|}{3}$. From these expressions of $Q$, $T$ we see that if a map exists in this case then $|V|$ is even and positive integer and $|V| \geq 8$.

Further, we may assume without loss of generality ${\rm lk}(0) = C_7(1, 2, {\bf 3}, 4, 5, {\bf 6}, 7)$. This implies ${\rm lk}(2)= C_7(c, 3, {\bf 4}, 0, 1, {\textbf {\textit a}}, b)$ or ${\rm lk}(2) = C_7(1, 0, {\bf 4}, 3, a, {\textbf {\textit b}}, c)$ for some $a, b, c \in V$. If ${\rm lk}(2) = C_7(1, 0, {\bf 4}, 3, a, {\textbf {\textit b}}, c)$ then we get three consecutive triangles incident with 1, which is not allowed. Hence ${\rm lk}(2)= C_7(c, 3, {\bf 4}, 0, 1, {\textbf {\textit a}}, b)$. Then, we make the following claim - For ${\rm lk}(2)= C_7(c, 3, {\bf 4}, 0, 1, {\bf \textit{a}}, b)$,  we have $(i)$ $b\neq7$ and $c\neq5$ $(ii)$ $c = 6$ implies $b = 5$ and $(iii)$ $c = 7$ implies $b = 6$. The reason for this can be seen as - If $b = 7$ then 17 is both an edge and a non-edge, which is not possible. If $c = 5$ then it is easy to see that  $b = 6$. For otherwise $\deg(5) > 5$.  Then considering ${\rm lk}(5)$ we see that 34 will be a non-edge. This is not possible. This proves $(i)$. If $c = 6$ then it is easy to see that $b \in \{5, 7\}$. For otherwise $\deg(6) > 5$. If $b = 7$ then 17 is both an edge and a non-edge. Hence $b = 5$,  this proves $(ii)$. If $c = 7$ then by the fact that degree of each vertex is 5, we have $b = 6$. This proves $(iii)$. Thus the claim. With this observation we proceed as follows\,:


When $|V| = 8$ then, considering the claim in previous paragraph, it is easy to see that $(c, b, a) \in \{(6, 5, 7), (7, 6, 5)\}$. But, then two quadrangles share more than one vertex. This is not possible. So $|V| \neq 8$.
When $|V| = 10$ then, considering the claim above, we have $(c, b, a) \in \{$(6, 5, 7), (6, 5, 8), (7, 6, 5), (7, 6, 8), (8, 5, 6), (8, 5, 7), (8, 5, 9), (8, 6, 5), (8, 6, 7), (8, 6, 9), (8, 9, 5), (8, 9, 6), (8, 9, 7)$\}$. We see that, for $(c, b, a)\in\{$(6, 5, 7), (7, 6, 5), (8, 5, 6), (8, 5, 7), (8, 6,  5), (8, 6, 7)$\}$, quadrangles [0, 5, 6, 7] and [1, 2, $b$, $a$] share more than one vertex, which is not possible. Thus $(c, b, a) \in \{$(6, 5, 8)$, $(7, 6 , 8)$, $(8, 5, 9)$, $(8, 6, 9)$, $(8, 9, 5)$, $(8, 9, 6)$, $(8, 9, 7)$\}$. If $(c, b, a) = (6, 5, 8)$ then ${\rm lk}(2)= C_7$(6, 3, {\bf 4}, 0, 1, {\bf 8}, 5 $)$. This implies ${\rm lk}(5)= C_7(4, 0, {\bf 7}, 6, 2, {\bf 1}, 8)$ and ${\rm lk}(4)= C_7(5, 0, {\bf 2}, 3, d, {\textbf {\textit e}}, 8)$ for some $d, e \in V$. It is easy to see that $d \in \{6, 7, 9\}$. If $d = 6$ then $e = 7$, for otherwise $\deg(6) > 5$. Then, there exist three quadrangles incident on 6,  which is not allowed. If $d = 7$ then $e \in \{1, 6\}$. If $e = 1$ then ${\rm lk}(7)= C_7(3, 6, {\bf 5}, 0, 1, {\bf 8}, 4)$. This implies $C_6(0, 2, 5, 8, 4, 7)\subseteq{\rm lk}(1)$. If $e = 6$ then there exist three quadrangles incident on 6. If $d = 9$ then we see that $e \in \{1, 6, 7\}$. In case $e \in \{1, 6\}$, then there exist three quadrangles incident on $e$. If $e = 7$ then there exist three consecutive triangles incident on 1, which is not allowed. Hence $(c, b, a) \neq (6, 5, 8)$. Proceeding in the same way for other values of $(c, b, a)$ we see that none of these values are admissible. Hence we conclude that $|V| \neq 10$.

When $|V| = 12$ then we have $(c, b, a) \in \{(6, 5, 8)$, $(7, 6, 8)$, $(8, 5, 9)$, $(8, 6, 9)$, $(8, 9, 5)$, $(8, 9, 6)$, $(8, 9, 6)$, $(8, 9, 7)$, $(8, 9, 10)\}$. When $(c, b, a) \in \{(8, 6, 9), (8, 9, 7)\}$,  then, as in case for $|V| = 10$ the map could not be constructed. So, we have $(c, b, a)\in\{(6, 5, 8)$, $(7, 6, 8)$, $(8, 5, 9)$, $(8, 9 , 5)$, $(8, 9, 6)$, $(8, 9, 10)\}$. For $(c, b, a) = (6, 5, 8)$ we get $M \cong K_{1, 12}(3^2, 4, 3, 4)$ by the map $0 \mapsto v_8$,  $1 \mapsto v_2$,  $2 \mapsto v_3$,  $3\mapsto v_{10}$,  $4 \mapsto v_5$,  $5 \mapsto v_4$,  $6 \mapsto v_{11}$,  $7 \mapsto v_9$, $8\mapsto v_0$, $9\mapsto v_6$,  $10\mapsto v_{7}$, $11\mapsto v_{1}$. For $(c, b, a) = (8, 9, 6)$. This is $M \cong K_{1, 12}(3^2, 4, 3, 4)$ by the map $i \mapsto v_i, 0\leq i \leq 11$. For $(c, b, a) = (8, 9, 10)$ we get $M \cong K_{1, 12}(3^2, 4, 3, 4)$ by the map $0 \mapsto v_0$,  $1 \mapsto v_1$,  $2 \mapsto v_7$,  $3\mapsto v_6$,  $4 \mapsto v_5$,  $5 \mapsto v_4$,  $6 \mapsto v_3$,  $7 \mapsto v_2$, $8\mapsto v_9$,  $9\mapsto v_{11}$,  $10\mapsto v_{10}$,  $11\mapsto v_{8}$ and $M \cong K_{1, 12}(3^2, 4, 3, 4)$ by the map $0 \mapsto v_8$,  $1 \mapsto v_2$,  $2 \mapsto v_9$,  $3\mapsto v_{11}$,  $4 \mapsto v_4$,  $5 \mapsto v_5$,  $6 \mapsto v_{10}$,  $7 \mapsto v_3$, $8\mapsto v_7$,  $9\mapsto v_6$,  $10\mapsto v_{1}$,  $11\mapsto v_{0}$. For other values of $(c, b, a)$ we do not get any object.

When $|V| = 14$, we see that the values of $(c, b, a)$ for which $M$ can be constructed are in
$\{(6, 5, 8)$, $(7, 6, 8)$, $(8, 9, 5)$,  $(8, 9, 6)$, $(8, 9, 10)\}$. Proceeding similarly as in the previous case of $|V| = 12$ we see that using none of these values can $M$ be actually constructed. This proves Lemma \ref{l1}.

%
\bigskip

\noindent\proof{Proof of Lemma \,\,\ref{l2}}:\,\, Let $M$ be a SEM of type $(3, 4, 6, 4)$ on a closed surface of Euler characteristic 0. The notation ${\rm lk}(i) = C_9({\boldsymbol i_{1}}, i_{2}, i_{3}, i_{4}, {\boldsymbol i_{5}}, i_{6}, {\boldsymbol i_{7}}, {\boldsymbol i_{8}}, i_{9})$ for the link of $i$ will mean that $[i, i_{1}, i_{2}, i_{3}, i_{4}, i_{5}]$ forms hexagonal face, $[i, i_{5}, i_{6}, i_{7}]$ and $[i, i_{1}, i_{9}, i_{8}]$ form quadrangular faces and $[i, i_{7}, i_{8}]$ forms triangular face. Let $|V|$ denote the cardinality of vertex set $V(M)$ and $T$, $Q$, $H$ denote the number of triangular, quadrangular and hexagonal faces of $M$, respectively. Then it is easy to see that $T = \frac {|V|}{3}$, $Q = \frac {2|V|}{4}$ and $H = \frac {|V|}{6}$. From these expressions of $T$, $Q$, $H$, we observe that a map in question exists if, $(i)$ $|V(M)|$ is a multiple of 6 and $(ii)$ $|V| \geq 12$. Thus it is sufficient  to assume that $|V| = 12$ and 18. Assume, without loss of generality, that ${\rm lk}(0) = C_9({\bf 1}, 2, 3, 4, {\bf 5}, 6, {\bf 7}, {\bf 8}, 9)$ and ${\rm lk}(1) = C_9({\bf 2}, 3, 4, 5, {\bf 0}, 8, {\bf 9}, {\bf 10}, 11)$. Proceeding similarly as in the previous cases we see that $|V| \neq 12$. So, $|V| = 18$ and in this case we get $M \cong K_{1, 18}(3, 4, 6, 4)$ by the map $i \mapsto v_i$, $0\leq i \leq17$ and $M \cong T_{1, 18}(3, 4, 6, 4)$ by the map $i \mapsto v_i$, $0\leq i \leq 17$.

If $M$ is a SEM of type $(4, 8^2)$ on a surface of Euler characteristic 0, then the notation ${\rm lk}(i) = C_{14}({\boldsymbol i_{1}}, i_{2}, i_{3}, i_{4}, i_{5}, i_{6}, {\boldsymbol i_{7}}, i_{8}, {\boldsymbol i_{9}}, i_{10}, i_{11}, i_{12}, i_{13}, i_{14})$, for the link of $i$ will mean that  $[i, i_{7}, i_{8}, i_{9}]$ forms a quadrangular face and $[i, i_{1}, i_{2}, i_{3}, i_{4}, i_{5}, i_{6}, i_{7}]$, $[i, i_{9}, i_{10}, i_{11}, i_{12}$, $i_{13}, i_{14}, i_1]$ form octagonal faces. Now, let $A$ and $B$ denote the vertex sets of octagonal faces $O_1$ and $O_2$, respectively. Then either $A \cap B = \phi$ or $A \cap B = \{u_1, u_2 \}$, where the set $\{u_1, u_2\}$ forms an edge $u_1u_2$ in the map. Let $Q$ and $O$ denote the number of quadrangular and octagonal faces in $M$, respectively. Then $Q = \frac {|V|}{4}$ and $O = \frac {2|V|}{8}$. Thus the map exists if $|V(M)|$ is a multiple of 4 and $|V| \geq 15$. So, in this case it is sufficient  to assume that $|V| = 16$ or 20. Let ${\rm lk}(0) = C_{14}({\bf 1}, 2, 3, 4, 5, 6, {\bf 7}, 8, {\bf 9}, 10, 11, 12, 13, 14)$. Computing successively, we do not get any example when $|V| = 16$. So, assume $|V| = 20$. In this case we get  $M \cong T_{1, 20}(4, 8^2)$ by the map \,: $0 \mapsto v_0, 1 \mapsto v_1, 2 \mapsto v_{14}, 3 \mapsto v_{13}, 4 \mapsto v_{12}, 5 \mapsto v_{11}, 6 \mapsto v_{10}, 7 \mapsto v_9, 8 \mapsto v_8, 9 \mapsto v_7, 10 \mapsto v_6, 11 \mapsto v_5, 12 \mapsto v_4, 13 \mapsto v_3, 14 \mapsto v_2, 15 \mapsto v_5, 16 \mapsto v_16, 17 \mapsto v_{19}, 18 \mapsto v_{18}, 19 \mapsto v_{17}$.

Let $M$ be a SEM of type $(3^4, 6)$ on a surface of Euler characteristic 0. The notation ${\rm lk}(i) = C_8([i_{1}, i_{2}, i_{3}, i_{4}, i_{5}], i_{6}, i_{7}, i_{8})$ for the link of $i$ will mean that $[i, i_{1}, i_{2}, i_{3}, i_{4}, i_{5}]$ forms hexagonal face and $[i, i_{5}, i_{6}]$, $[i, i_{6}, i_{7}]$, $[i, i_{7}, i_{8}]$ and $[i, i_{1}, i_{8}]$ form triangular faces. Let $|V(M)|$ denote the number of vertices in the set $V(M)$. Let $T(M)$ denote the number of
triangular faces and $H(M)$ denote the number of hexagonal faces in $M$. Then $H(M) = \frac {|V(M)|}{6}$ and $T(M) = \frac {4|V(M)|}{3}$. From this observe that map exists if $|V(M)|$ is a multiple of 6 and $|V| \geq 12$. Thus we have observed that it is sufficient to assume that $|V| = 12$ or 18. Let ${\rm lk}(0) = C_8([1, 2, 3, 4, 5], 6, 7, 8)$. Proceeding as in the previous case we see that $|V| \neq 12$ so $|V| = 18$. In this case we get $M \cong T_{1, 18}(3^4, 6)$ by the map $i \mapsto v_i$, $0 \leq i \leq 17$.

Let $M$ be a SEM of type $(3, 6, 3, 6)$ on a closed surface of Euler characteristic 0. Link of a vertex $i$ of the map is denoted as ${\rm lk}(i) = C_{10}([{\boldsymbol i_1}, i_2, i_3, i_4, {\boldsymbol i_5}], [{\boldsymbol i_6}, i_7, i_8, i_9, {\boldsymbol i_{10}}])$. The representation of the link $(i)$ will mean that $[i, i_1, i_2, i_3, i_4, i_5]$ and $[i, i_6, i_7, i_8, i_9, i_{10}]$ form hexagonal faces and $[i, i_1, i_{10}]$, $[i, i_5, i_6]$ form triangular faces. Let $T$ and $H$ represents the number of triangular faces and hexagonal faces, respectively, then $T = \frac{2|V|}{3}$ and $H = \frac{|V|}{3}$. Thus the map exists if $|V(M)|$ is a multiple of 3 and $|V| \geq 12$. So, we proceed in the following manner. It is sufficient to assume that $|V| = 12$, 15 and 18. Let ${\rm lk}(0) = C_{10}([{\bf 1}, 2, 3, 4, {\bf 5}],[{\bf 6}, 7, 8, 9, {\bf 10}])$. Proceeding as in the previous cases, we see that in this case we do not get any examples. This proves Lemma \ref{l2}.

\noindent{\bf Proof of Theorem} \ref{t1}\,  To show that these maps are non-isomorphic we only need to show that $T_{1, 12}(3^3, 4^2) \not\cong T_{2, 12}(3^3, 4^2) \not\cong T_{3, 12}(3^3, 4^2)$, $K_{1, 12}(3^3, 4^2) \not\cong K_{2, 12}(3^3, 4^2)$ and $T_{1, 14}(3^3, 4^2) \not\cong T_{2, 14}(3^3, 4^2)$. For this we compute the characteristic polynomials of the incidence matrix of the edge graph of respective maps. These are as follows\,:

\smallskip
\noindent $p(EG(T_{1, 12}(3^3, 4^2)) = x^{12} - 30\,x^{10} - 24\,x^9 + 237\,x^8 + 192\,x^7 - 708\,x^6 - 408\,x^5 + 708\,x^4 + 208\,x^3 - 240 \,x^2 $\\

\noindent $p(EG(T_{2, 12}(3^3, 4^2)) = x^{12} - 30\,x^{10} - 32\,x^9 + 237\,x^8 + 360\,x^7 - 484\,x^6 - 696\,x^5 + 516\,x^4 + 368\,x^3 - 240\,x^2$\\

\noindent $p(EG(T_{3, 12}(3^3, 4^2)) = x^{12} - 30\,x^{10} - 32\,x^9 + 231\,x^8 + 384\,x^7 - 388\,x^6 - 960\, x^5 + 63\,x^4 + 896\,x^3 + 258\,x^2 - 288\, - 135$\\

\noindent $p(EG(K_{1, 12}(3^3, 4^2)) = x^{12} - 30\,x^{10} - 24\,x^9 + 243\,x^8 + 192\,x^7 -868\,x^6 - 528\,x^5 + 1527\,x^4 + 576\,x^3 - 1278\,x^2 - 216\,x + 405$\\


\noindent $p(EG(K_{2, 12}(3^3, 4^2)) = x^{12} - 30\,x^{10} - 32\,x^9 + 235\,x^8 + 368\,x^7 - 452\,x^6 - 784\,x^5 + 359\,x^4 + 592\,x^3 - 158\,x^2 - 144\,x + 45$\\


\noindent $p(EG(T_{1, 14}(3^3, 4^2)) = x^{14} - 34\,x^{12} - 30\,x^{11} + 344\,x^{10} +467\,x^9 - 1179\,x^8 - 2119\,x^7 + 597\,x^6 + 2264\,x^5 + 632\,x^4 - 559\,x^3 - 365\,x^2 - 74\,x - 5$\\

\noindent $p(EG(T_{2, 14}(3^3, 4^2)) = x^{14} - 35\,x^{12} -\,28x^{11} + 399\,x^{10} + 420\,x^9 - 2107\,x^8 - 2384\,x^7 + 5544\,x^6 + 6244\,x^5 - 6790\,x^4 - 7112\,x^3 + 3157\,x^2 + 2184\,x - 845$\\
\smallskip

Now the proof of the theorem follows from Lemma \ref{l1} and the fact that if two maps are isomorphic then the characteristic polynomials of the incidence matrices of their edge graphs are identical. \hfill $\Box$

\section{Infinite series of semi-equivelar maps}

Infinite series of all equivelar maps, $i.e.$ maps of types $(3^6)$, $(4^4)$ and $(6^3)$, on the surfaces of torus and
Klein bottle are shown in the figures below (Fig.(4.1)-Fig.(4.5)). For the series of equivelar map of type $(6^3)$ on Klein bottle consider
the dual of the equivelar map of type $(3^6)$ on this surface.

\bigskip

\hrule

\begin{center}


\begin{picture}(0,0)(78,30)
\setlength{\unitlength}{10mm}

\drawpolygon(0,0)(4,0)(4,2)(0,2)

\drawpolygon(5,0)(7,0)(7,2)(5,2)

\drawline[AHnb=0](1,0)(1,2)
\drawline[AHnb=0](2,0)(2,2)
\drawline[AHnb=0](3,0)(3,2)
\drawline[AHnb=0](6,0)(6,2)

\drawline[AHnb=0](0,1)(1,2)
\drawline[AHnb=0](0,1)(4,1)
\drawline[AHnb=0](0,0)(2,2)
\drawline[AHnb=0](1,0)(3,2)
\drawline[AHnb=0](2,0)(4,2)
\drawline[AHnb=0](3,0)(4,1)

\drawline[AHnb=0](5,1)(6,2)
\drawline[AHnb=0](5,1)(7,1)
\drawline[AHnb=0](5,0)(7,2)
\drawline[AHnb=0](6,0)(7,1)

\put(-.1,-.25){\scriptsize {\tiny 2}}
\put(.9,-.25){\scriptsize {\tiny 4}}
\put(1.9,-.25){\scriptsize {\tiny 6}}
\put(2.9,-.25){\scriptsize {\tiny 8}}
\put(3.9,-.25){\scriptsize {\tiny 10}}
\put(4.9,-.25){\scriptsize {\tiny $2n$-2}}
\put(5.9,-.25){\scriptsize {\tiny $2n$}}
\put(7,-.25){\scriptsize {\tiny 2}}

\put(-.2,.91){\scriptsize {\tiny 1}}
\put(.8,1.1){\scriptsize {\tiny 3}}
\put(1.8,1.1){\scriptsize {\tiny 5}}
\put(2.8,1.1){\scriptsize {\tiny 7}}
\put(3.8,1.1){\scriptsize {\tiny 9}}
\put(5.1,.8){\scriptsize {\tiny $2n$-3}}
\put(6.1,.8){\scriptsize {\tiny $2n$-1}}
\put(7.1,.91){\scriptsize {\tiny 1}}

\put(-.2,2.1){\scriptsize {\tiny $2n$-4}}
\put(.8,2.1){\scriptsize {\tiny $2n$-2}}
\put(1.95,2.1){\scriptsize {\tiny $2n$}}
\put(2.95,2.1){\scriptsize {\tiny 2}}
\put(3.95,2.1){\scriptsize {\tiny 4}}

\put(4.8,2.1){\scriptsize {\tiny $2n$-8}}
\put(5.8,2.1){\scriptsize {\tiny $2n$-6}}
\put(6.8,2.1){\scriptsize {\tiny $2n$-4}}
%
%
\put(4.3,1){\scriptsize $\cdots$}

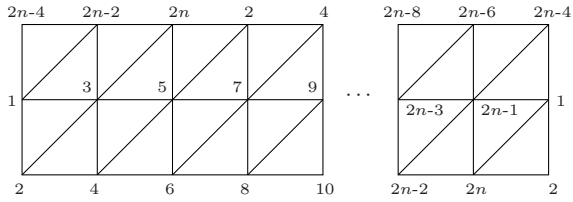
\put(.5,-.8){\scriptsize{\bf Fig.(4.1)\,: equivelar map of type\,-$(3^6)$}}

\put(-.3,2.8){\scriptsize{\bf Equivelar maps on torus}\,: for $n \geq 7$}

\end{picture}
\end{center}

\begin{center}

\begin{picture}(0,0)(-6,20)
\setlength{\unitlength}{9mm}

\drawpolygon(0,0)(4,0)(4,2)(0,2)

\drawpolygon(5,0)(7,0)(7,2)(5,2)

\drawline[AHnb=0](1,0)(1,2)
\drawline[AHnb=0](2,0)(2,2)
\drawline[AHnb=0](3,0)(3,2)
\drawline[AHnb=0](0,1)(4,1)
\drawline[AHnb=0](5,1)(7,1)
\drawline[AHnb=0](6,0)(6,2)
\put(4.3,1){\scriptsize $\cdots$}
\put(-.05,-.25){\scriptsize {\tiny 2}}
\put(.95,-.25){\scriptsize {\tiny 4}}
\put(1.95,-.25){\scriptsize {\tiny 6}}
\put(2.95,-.25){\scriptsize {\tiny 8}}
\put(3.95,-.25){\scriptsize {\tiny 10}}
%
\put(4.8,-.25){\scriptsize {\tiny $2n$-2}}
\put(5.9,-.25){\scriptsize {\tiny $2n$}}
\put(6.95,-.25){\scriptsize {\tiny 2}}
\put(-.15,2.1){\scriptsize {\tiny $2n$-4}}
\put(.85,2.1){\scriptsize {\tiny $2n$-2}}
\put(1.95,2.1){\scriptsize {\tiny $2n$}}
\put(2.95,2.1){\scriptsize {\tiny 2}}
\put(3.95,2.1){\scriptsize {\tiny 4}}
%
\put(4.85,2.1){\scriptsize {\tiny $2n$-8}}
\put(5.85,2.1){\scriptsize {\tiny $2n$-6}}
\put(6.85,2.1){\scriptsize {\tiny $2n$-4}}

\put(-.2,.91){\scriptsize {\tiny 1}}
\put(.8,1.1){\scriptsize {\tiny 3}}
\put(1.8,1.1){\scriptsize {\tiny 5}}
\put(2.8,1.1){\scriptsize {\tiny 7}}
\put(3.8,1.1){\scriptsize {\tiny 9}}
\put(5.1,.8){\scriptsize {\tiny $2n$-3}}
\put(6.1,.8){\scriptsize {\tiny $2n$-1}}
\put(7.1,.91){\scriptsize {\tiny 1}}

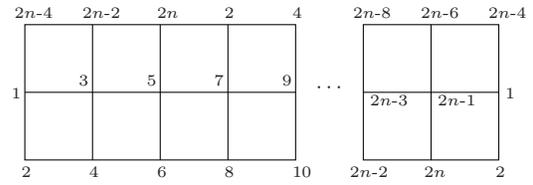
\put(.5,-.8){\scriptsize{\bf Fig.(4.2)\,: equivelar map of type\,-$(4^4)$}}


\end{picture}

\end{center}

\begin{center}


\begin{picture}(0,0)(70,45)
\setlength{\unitlength}{8mm}

\drawpolygon(0,0)(1,-.5)(2,0)(3,-.5)(4,0)(5,-.5)(6,0)(7,-.5)(8,0)(8,1)(7,1.5)(6,1)(5,1.5)(4,1)(3,1.5)(2,1)(1,1.5)(0,1)

\drawpolygon(10,0)(11,-.5)(12,0)(13,-.5)(14,0)(14,1)(13,1.5)(12,1)(11,1.5)(10,1)

\drawline[AHnb=0](2,0)(2,1)
\drawline[AHnb=0](4,0)(4,1)
\drawline[AHnb=0](6,0)(6,1)
\drawline[AHnb=0](12,0)(12,1)

\put(.9,-.8){\scriptsize {\tiny 2}}
\put(2.9,-.8){\scriptsize {\tiny 4}}
\put(4.9,-.8){\scriptsize {\tiny 6}}
\put(6.9,-.8){\scriptsize {\tiny 8}}
\put(10.65,-.8){\scriptsize {\tiny $2n$-2}}
\put(12.8,-.8){\scriptsize {\tiny $2n$}}

\put(.65,1.65){\scriptsize {\tiny $2n$-3}}
\put(2.65,1.65){\scriptsize {\tiny $2n$-1}}
\put(4.95,1.65){\scriptsize {\tiny 1}}
\put(6.95,1.65){\scriptsize {\tiny 3}}
\put(10.65,1.65){\scriptsize {\tiny $2n$-7}}
\put(12.65,1.65){\scriptsize {\tiny $2n$-5}}

\put(-.65,1.15){\scriptsize {\tiny $2n$-4}}
\put(1.65,1.2){\scriptsize {\tiny $2n$-2}}
\put(3.80,1.2){\scriptsize {\tiny $2n$}}
\put(5.90,1.2){\scriptsize {\tiny 2}}
\put(7.90,1.2){\scriptsize {\tiny 4}}

\put(-.05,-.35){\scriptsize {\tiny 1}}
\put(1.95,-.35){\scriptsize {\tiny 3}}
\put(3.95,-.35){\scriptsize {\tiny 5}}
\put(5.95,-.35){\scriptsize {\tiny 7}}
\put(7.95,-.35){\scriptsize {\tiny 9}}
\put(9.6,-.4){\scriptsize {\tiny $2n$-3}}
\put(11.7,-.4){\scriptsize {\tiny $2n$-1}}
\put(13.95,-.4){\scriptsize {\tiny 1}}
\put(9.6,1.2){\scriptsize {\tiny $2n$-8}}
\put(11.7,1.3){\scriptsize {\tiny $2n$-6}}
\put(13.9,1.25){\scriptsize {\tiny $2n$-4}}


\put(9,.5){\scriptsize $\cdots$}

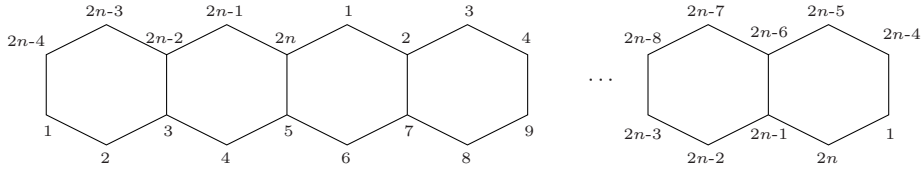
\put(3,-1.5){\scriptsize{\bf Fig.(4.3)\,: equivelar map of type\,-$(6^3)$}}


\end{picture}
\end{center}

\begin{center}


\begin{picture}(0,0)(70,95)
\setlength{\unitlength}{9mm}

\drawpolygon(0,0)(2,0)(2,3)(0,3)
\drawpolygon(3,0)(5,0)(5,3)(3,3)

\drawline[AHnb=0](0,1)(2,1)
\drawline[AHnb=0](0,2)(2,2)
\drawline[AHnb=0](3,1)(5,1)
\drawline[AHnb=0](3,2)(5,2)
\drawline[AHnb=0](1,0)(1,3)
\drawline[AHnb=0](4,0)(4,3)

\drawline[AHnb=0](0,0)(1,1)
\drawline[AHnb=0](1,0)(2,1)
\drawline[AHnb=0](3,0)(4,1)
\drawline[AHnb=0](4,0)(5,1)

\drawline[AHnb=0](0,1)(1,2)
\drawline[AHnb=0](1,1)(2,2)
\drawline[AHnb=0](3,1)(4,2)
\drawline[AHnb=0](4,1)(5,2)

\drawline[AHnb=0](0,2)(1,3)
\drawline[AHnb=0](1,2)(2,3)
\drawline[AHnb=0](3,2)(4,3)
\drawline[AHnb=0](4,2)(5,3)

\put(2.3,1.5){\scriptsize $\cdots$}

\put(-.2,0.05){\scriptsize {\tiny 1}}
\put(-.2,1.05){\scriptsize {\tiny 3}}
\put(-.2,2.05){\scriptsize {\tiny 2}}
\put(-.2,3.05){\scriptsize {\tiny 1}}

\put(.8,0.05){\scriptsize {\tiny 4}}
\put(.8,1.05){\scriptsize {\tiny 6}}
\put(.8,2.05){\scriptsize {\tiny 5}}
\put(.8,3.05){\scriptsize {\tiny 4}}

\put(1.8,0.05){\scriptsize {\tiny 7}}
\put(1.8,1.05){\scriptsize {\tiny 9}}
\put(1.8,2.05){\scriptsize {\tiny 8}}
\put(1.8,3.05){\scriptsize {\tiny 7}}

\put(2.4,0.05){\scriptsize {\tiny $3n$-5}}
\put(2.4,1.05){\scriptsize {\tiny $3n$-3}}
\put(2.4,2.05){\scriptsize {\tiny $3n$-4}}
\put(2.4,3.05){\scriptsize {\tiny $3n$-5}}

\put(3.4,0.05){\scriptsize {\tiny $3n$-2}}
\put(3.4,1.05){\scriptsize {\tiny $3n$}}
\put(3.4,2.05){\scriptsize {\tiny $3n$-1}}
\put(3.4,3.05){\scriptsize {\tiny $3n$-2}}

\put(5.1,0.05){\scriptsize {\tiny 1}}
\put(5.1,1.05){\scriptsize {\tiny 2}}
\put(5.1,2.05){\scriptsize {\tiny 3}}
\put(5.1,3.05){\scriptsize {\tiny 1}}

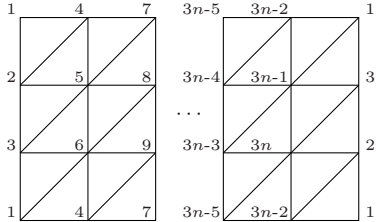
\put(-.5,-.8){\scriptsize{\bf Fig.(4.4)\,: equivelar map of type\,-$(3^6)$}}

\end{picture}

\end{center}

\begin{center}

\begin{picture}(0,0)(-10,85)
\setlength{\unitlength}{9mm}

\drawpolygon(0,0)(2,0)(2,3)(0,3)
\drawpolygon(3,0)(5,0)(5,3)(3,3)

\drawline[AHnb=0](0,1)(2,1)
\drawline[AHnb=0](0,2)(2,2)
\drawline[AHnb=0](3,1)(5,1)
\drawline[AHnb=0](3,2)(5,2)
\drawline[AHnb=0](1,0)(1,3)
\drawline[AHnb=0](4,0)(4,3)

\put(2.3,1.5){\scriptsize $\cdots$}

\put(-.2,0.05){\scriptsize {\tiny 1}}
\put(-.2,1.05){\scriptsize {\tiny 3}}
\put(-.2,2.05){\scriptsize {\tiny 2}}
\put(-.2,3.05){\scriptsize {\tiny 1}}

\put(.8,0.05){\scriptsize {\tiny 4}}
\put(.8,1.05){\scriptsize {\tiny 6}}
\put(.8,2.05){\scriptsize {\tiny 5}}
\put(.8,3.05){\scriptsize {\tiny 4}}

\put(1.8,0.05){\scriptsize {\tiny 7}}
\put(1.8,1.05){\scriptsize {\tiny 9}}
\put(1.8,2.05){\scriptsize {\tiny 8}}
\put(1.8,3.05){\scriptsize {\tiny 7}}

\put(2.4,0.05){\scriptsize {\tiny $3n$-5}}
\put(2.4,1.05){\scriptsize {\tiny $3n$-3}}
\put(2.4,2.05){\scriptsize {\tiny $3n$-4}}
\put(2.4,3.05){\scriptsize {\tiny $3n$-5}}

\put(3.4,0.05){\scriptsize {\tiny $3n$-2}}
\put(3.4,1.05){\scriptsize {\tiny $3n$}}
\put(3.4,2.05){\scriptsize {\tiny $3n$-1}}
\put(3.4,3.05){\scriptsize {\tiny $3n$-2}}

\put(5.1,0.05){\scriptsize {\tiny 1}}
\put(5.1,1.05){\scriptsize {\tiny 2}}
\put(5.1,2.05){\scriptsize {\tiny 3}}
\put(5.1,3.05){\scriptsize {\tiny 1}}

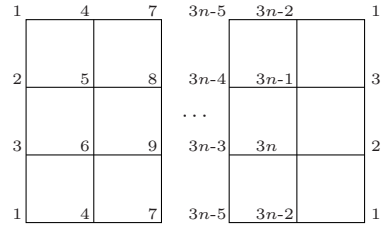
\put(-.5,-.8){\scriptsize{\bf Fig.(4.5)\,: equivelar map of type\,-$(4^4)$}}


\put(-9.8,3.5){\scriptsize {\bf Equivelar maps on Klein bottle\,:} for $n \geq 3$}

\end{picture}

\end{center}

\vspace{9.5cm}

\hrule

\bigskip

An infinite series of all eight types of semi-equivelar maps on torus and Klein bottle can be obtained from the infinite series
of equivelar maps on the respective surfaces by using subdivision and truncation. Here we construct such
series for the torus, and a similar construction will work for the Klein bottle also. To construct the SEMs of types
$(3^3,4^2)$, $(3^3, 4, 3, 4)$, $(3,6,3,6)$, $(4, 8^2)$ and $(4,6,12)$ we consider an equivelar map $M$ of type $(4^4)$ with
the number of vertices 2$n$ for $n\geq 7$ (Fig 4.2). Then the map has exactly $2n$ faces (quadrangles). Now we use the operations
as follows\,:

If we subdivide one layer of quadrangular faces by one diagonal in each we get a SEM of type $(3^3, 4^2)$ (see Fig.(4.6)). Similarly,
subdivision of alternate quadrangular faces of the equivelar map by one diagonal in each such that no two subdivided quadrangular
face share an edge, leads to a SEM of type $(3^2,4,3,4)$ (see Fig.(4.7)). During the construction of these SEMs we see that the number of
vertices remains the same (equal to the number of vertices in the equivelar map).

Truncation of an equivelar map of type $(p, q)$ (each face of the map is a $p$-gon and each vertex belongs to exactly
$q$ faces) with $n$ vertices, along its vertices leads to a SEM of type $(q, (2p)^2)$ with $qn$ vertices. Thus, if we
truncate the equivelar map along its vertices it leads to a SEM of type $(4, 8^2)$ (see Fig.(4.8)) with 8$n$ vertices.

If we subdivide the equivelar map as shown in Fig.(4.9) we get a SEM of type $(3,6,3,6)$. In this process we see easily that
the number of vertices in the SEM is exactly 4$n$ (twice the number of quadrangular faces in the equivelar map). Moreover, if we truncate
this SEM of type $(3,6,3,6)$ along its vertices we get a SEM of type $(4,6,12)$ with (2$n$ + 4.(2$n$)) vertices.

\smallskip

\hrule

\begin{center}
\begin{picture}(0,0)(78,28) 
\setlength{\unitlength}{5mm}

\drawpolygon(0,0)(6,0)(6,4)(0,4)
\drawline[AHnb=0](0,2)(6,2)
\drawline[AHnb=0](2,0)(2,4)

\drawline[AHnb=0](4,0)(4,4)

\drawpolygon(9,0)(13,0)(13,4)(9,4)

\drawline[AHnb=0](9,2)(13,2)
\drawline[AHnb=0](2,2)(0,4)
\drawline[AHnb=0](4,2)(2,4)
\drawline[AHnb=0](6,2)(4,4)

\drawline[AHnb=0](11,0)(11,4)
\drawline[AHnb=0](13,2)(11,4)
\drawline[AHnb=0](11,2)(9,4)

\put(7,1.9){\scriptsize $\cdots$}

\put(-.1,-.4){\scriptsize $a_2$}
\put(1.9,-.4){\scriptsize $a_4$}
\put(3.9,-.4){\scriptsize $a_6$}
\put(5.9,-.4){\scriptsize $a_8$}

\put(8.2,-.4){\scriptsize $a_{2n-2}$}
\put(10.6,-.4){\scriptsize $a_{2n}$}

\put(12.9,-.4){\scriptsize $a_2$}
\put(-.6,1.7){\scriptsize $a_1$}
\put(1.38,1.56){\scriptsize $a_3$}
\put(3.38,1.56){\scriptsize $a_5$}
\put(5.38,1.56){\scriptsize $a_7$}

\put(7.38,1.56){\scriptsize $a_{2n-3}$}
\put(9.4,1.56){\scriptsize $a_{2n-1}$}

\put(13.2,1.7){\scriptsize $a_1$}

\put(-.8,4.2){\scriptsize $a_{2n-4}$}
\put(1.2,4.2){\scriptsize $a_{2n-2}$}
\put(3.7,4.2){\scriptsize $a_{2n}$}
\put(5.7,4.2){\scriptsize $a_2$}

\put(8.2,4.2){\scriptsize $a_{2n-8}$}
\put(10.2,4.2){\scriptsize $a_{2n-6}$}
\put(12.3,4.2){\scriptsize $a_{2n-4}$}

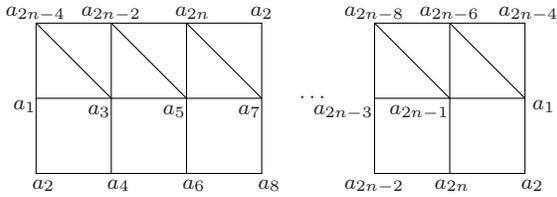
\put(.5,-1.5){\scriptsize{\bf Fig.(4.6)\,: SEM of type\,-$(3^3, 4^2)$, $n\geq 7$}}
%
\put(-.5,5.5){\scriptsize {\bf Infinite series of SEMs on torus}}
\end{picture}

\end{center}

\begin{center}
\begin{picture}(0,0)(-10,19) 
\setlength{\unitlength}{10mm}

\drawpolygon(0,0)(3,0)(3,2)(0,2)

\drawpolygon(4.5,0)(6.5,0)(6.5,2)(4.5,2)

\drawline[AHnb=0](0,1)(3,1)
\drawline[AHnb=0](1,0)(1,2)
\drawline[AHnb=0](2,0)(2,2)

\drawline[AHnb=0](5.5,0)(5.5,2)
\drawline[AHnb=0](4.5,1)(6.5,1)

\drawline[AHnb=0](0,1)(1,2)
\drawline[AHnb=0](2,1)(3,2)
%
\drawline[AHnb=0](1,1)(2,0)

\drawline[AHnb=0](5.5,1)(6.5,0)
\drawline[AHnb=0](4.5,1)(5.5,2)

\put(-.2,-.2){\scriptsize  $a_1$}
\put(.9,-.2){\scriptsize  $a_2$}
\put(1.9,-.2){\scriptsize  $a_3$}
\put(2.9,-.2){\scriptsize  $a_4$}
\put(4.1,-.2){\scriptsize  $a_{n-1}$}
\put(5.35,-.2){\scriptsize  $a_{n}$}
\put(6.65,-.2){\scriptsize  $a_1$}
\put(-.3,.77){\scriptsize  $b_1$}
\put(.65,.77){\scriptsize  $b_2$}
\put(1.65,.77){\scriptsize  $b_3$}
\put(2.65,.77){\scriptsize $b_4$}
\put(3.85,.88){\scriptsize $b_{n-1}$}
%
\put(5.1,.77){\scriptsize $b_{n}$}
\put(6.6,.9){\scriptsize  $b_1$}

\put(-.4,2.15){\scriptsize  $a_{n-3}$}
\put(.6,2.15){\scriptsize  $a_{n-2}$}
\put(1.8,2.15){\scriptsize $a_{n-1}$}
\put(2.8,2.15){\scriptsize $a_n$}

\put(4.1,2.15){\scriptsize  $a_{n-5}$}

\put(5.1,2.15){\scriptsize  $a_{n-4}$}
\put(6.1,2.15){\scriptsize $a_{n-3}$}

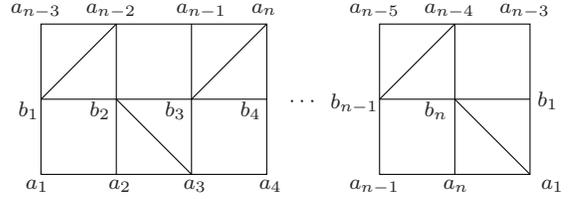
\put(-.2,-.8){\scriptsize{\bf Fig.(4.7)\,: SEM of type\,-$(3^2,4,3,4)$, $n = 2k$, $k\geq4$}}

\put(3.3,.9){\scriptsize $\cdots$}
\end{picture}
\end{center}

\begin{center}

\begin{picture}(0,0)(75,63) 

\setlength{\unitlength}{8.5mm}

\drawpolygon(0,0)(8,0)(8,4)(0,4)

\drawpolygon(10,0)(14,0)(14,4)(10,4)

\drawline[AHnb=0](0,2)(8,2)
\drawline[AHnb=0](2,0)(2,4)
\drawline[AHnb=0](4,0)(4,4)
\drawline[AHnb=0](6,0)(6,4)
\drawline[AHnb=0](12,0)(12,4)
\drawline[AHnb=0](10,2)(14,2)

\drawline[AHnb=0](.5,2)(0,2.5)
\drawline[AHnb=0](.5,2)(0,1.5)

\drawline[AHnb=0](2.5,2)(2,2.5)
\drawline[AHnb=0](2.5,2)(2,1.5)

\drawline[AHnb=0](4.5,2)(4,2.5)
\drawline[AHnb=0](4.5,2)(4,1.5)

\drawline[AHnb=0](6.5,2)(6,2.5)
\drawline[AHnb=0](6.5,2)(6,1.5)

\drawline[AHnb=0](10.5,2)(10,2.5)
\drawline[AHnb=0](10.5,2)(10,1.5)

\drawline[AHnb=0](12.5,2)(12,2.5)
\drawline[AHnb=0](12.5,2)(12,1.5)

\drawline[AHnb=0](1.5,2)(2,2.5)
\drawline[AHnb=0](1.5,2)(2,1.5)

\drawline[AHnb=0](3.5,2)(4,2.5)
\drawline[AHnb=0](3.5,2)(4,1.5)

\drawline[AHnb=0](5.5,2)(6,2.5)
\drawline[AHnb=0](5.5,2)(6,1.5)

\drawline[AHnb=0](7.5,2)(8,2.5)
\drawline[AHnb=0](7.5,2)(8,1.5)

\drawline[AHnb=0](11.5,2)(12,2.5)
\drawline[AHnb=0](11.5,2)(12,1.5)

\drawline[AHnb=0](13.5,2)(14,2.5)
\drawline[AHnb=0](13.5,2)(14,1.5)

\drawline[AHnb=0](0,3.5)(.5,4)
\drawline[AHnb=0](2,3.5)(2.5,4)
\drawline[AHnb=0](4,3.5)(4.5,4)
\drawline[AHnb=0](6,3.5)(6.5,4)
\drawline[AHnb=0](10,3.5)(10.5,4)
\drawline[AHnb=0](12,3.5)(12.5,4)

\drawline[AHnb=0](1.5,4)(2,3.5)
\drawline[AHnb=0](3.5,4)(4,3.5)
\drawline[AHnb=0](5.5,4)(6,3.5)
\drawline[AHnb=0](7.5,4)(8,3.5)

\drawline[AHnb=0](11.5,4)(12,3.5)
\drawline[AHnb=0](13.5,4)(14,3.5)

\drawline[AHnb=0](.5,0)(0,.5)
\drawline[AHnb=0](2.5,0)(2,.5)
\drawline[AHnb=0](4.5,0)(4,.5)
\drawline[AHnb=0](6.5,0)(6,.5)
\drawline[AHnb=0](10.5,0)(10,.5)
\drawline[AHnb=0](12.5,0)(12,.5)

\drawline[AHnb=0](1.5,0)(2,.5)
\drawline[AHnb=0](3.5,0)(4,.5)
\drawline[AHnb=0](5.5,0)(6,.5)
\drawline[AHnb=0](7.5,0)(8,.5)
\drawline[AHnb=0](11.5,0)(12,.5)
\drawline[AHnb=0](13.5,0)(14,.5)

\put(8.7,1.9){\scriptsize $\cdots$}

%

\put(.4,-.2){\scriptsize {\tiny $c_2$}}
\put(2.4,-.2){\scriptsize {\tiny $c_4$}}
\put(4.4,-.2){\scriptsize {\tiny $c_6$}}
\put(6.4,-.2){\scriptsize {\tiny $c_8$}}
\put(10.3,-.2){\scriptsize {\tiny $c_{2n-2}$}}
\put(12.35,-.2){\scriptsize {\tiny $c_{2n}$}}

\put(1.35,-.2){\scriptsize {\tiny $a_4$}}
\put(3.35,-.2){\scriptsize {\tiny $a_6$}}
\put(5.35,-.2){\scriptsize {\tiny $a_8$}}
\put(7.35,-.2){\scriptsize {\tiny $a_{10}$}}
\put(11.35,-.2){\scriptsize {\tiny $a_{2n}$}}
\put(13.35,-.2){\scriptsize {\tiny $a_2$}}

\put(.5,1.75){\scriptsize {\tiny $c_1$}}
\put(2.5,1.75){\scriptsize {\tiny $c_3$}}
\put(4.5,1.75){\scriptsize {\tiny $c_5$}}
\put(6.5,1.75){\scriptsize {\tiny $c_7$}}
\put(10.4,1.75){\scriptsize {\tiny $c_{2n-3}$}}
\put(12.4,1.75){\scriptsize {\tiny $c_{2n-1}$}}

\put(.1,.5){\scriptsize {\tiny $b_2$}}
\put(2.1,.5){\scriptsize {\tiny $b_4$}}
\put(4.1,.5){\scriptsize {\tiny $b_6$}}
\put(6.1,.5){\scriptsize {\tiny $b_8$}}
\put(8.1,.5){\scriptsize {\tiny $b_{10}$}}
\put(10.1,.5){\scriptsize {\tiny $b_{2n-2}$}}
\put(12.1,.5){\scriptsize {\tiny $b_{2n}$}}
\put(14.1,.5){\scriptsize {\tiny $b_2$}}

\put(.1,2.5){\scriptsize {\tiny $b_1$}}
\put(2.1,2.5){\scriptsize {\tiny $b_3$}}
\put(4.1,2.5){\scriptsize {\tiny $b_5$}}
\put(6.1,2.5){\scriptsize {\tiny $b_7$}}
\put(8.1,2.5){\scriptsize {\tiny $b_9$}}
\put(10.1,2.5){\scriptsize {\tiny $b_{2n-3}$}}
\put(12.1,2.5){\scriptsize {\tiny $b_{2n-1}$}}
\put(14.1,2.5){\scriptsize {\tiny $b_1$}}

\put(.1,3.35){\scriptsize {\tiny $d_{2n-4}$}}
\put(2.1,3.35){\scriptsize {\tiny $d_{2n-2}$}}
\put(4.1,3.35){\scriptsize {\tiny $d_{2n}$}}
\put(6.1,3.35){\scriptsize {\tiny $d_2$}}
\put(8.1,3.35){\scriptsize {\tiny $d_4$}}
\put(10.1,3.35){\scriptsize {\tiny $d_{2n-8}$}}
\put(12.1,3.35){\scriptsize {\tiny $d_{2n-6}$}}
\put(14.1,3.35){\scriptsize {\tiny $d_{2n-4}$}}

\put(.1,1.35){\scriptsize {\tiny $d_1$}}
\put(2.1,1.35){\scriptsize {\tiny $d_3$}}
\put(4.1,1.35){\scriptsize {\tiny $d_5$}}
\put(6.1,1.35){\scriptsize {\tiny $d_7$}}
\put(8.1,1.35){\scriptsize {\tiny $d_9$}}
\put(10.1,1.35){\scriptsize {\tiny $d_{2n-3}$}}
\put(12.1,1.35){\scriptsize {\tiny $d_{2n-1}$}}
\put(14.1,1.35){\scriptsize {\tiny $d_1$}}

\put(.2,4.15){\scriptsize {\tiny $c_{2n-4}$}}
\put(2.2,4.15){\scriptsize {\tiny $c_{2n-2}$}}
\put(4.2,4.15){\scriptsize {\tiny $c_{2n}$}}
\put(6.2,4.15){\scriptsize {\tiny $c_2$}}
\put(10.2,4.15){\scriptsize {\tiny $c_{2n-8}$}}
\put(12.2,4.15){\scriptsize {\tiny $c_{2n-6}$}}

\put(1.25,1.8){\scriptsize {\tiny $a_3$}}
\put(3.25,1.8){\scriptsize {\tiny $a_5$}}
\put(5.25,1.8){\scriptsize {\tiny $a_7$}}
\put(7.2,1.8){\scriptsize {\tiny $a_{9}$}}
\put(10.65,2.1){\scriptsize {\tiny $a_{2n-1}$}}
\put(13.15,2.1){\scriptsize {\tiny $a_1$}}

\put(1.1,4.15){\scriptsize {\tiny $a_{2n-2}$}}
\put(3.2,4.15){\scriptsize {\tiny $a_{2n}$}}
\put(5.3,4.15){\scriptsize {\tiny $a_2$}}
\put(7.3,4.15){\scriptsize {\tiny $a_{4}$}}
\put(11.1,4.15){\scriptsize {\tiny $a_{2n-6}$}}
\put(13.1,4.15){\scriptsize {\tiny $a_{2n-4}$}}

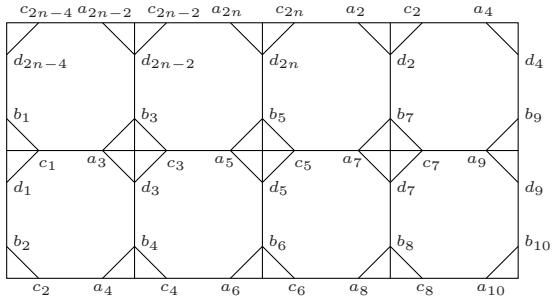
\put(.5,-.8){\scriptsize{\bf Fig.(4.8)\,: SEM of type\,-$(4, 8^2)$}}

\put(14.8,1.9){\scriptsize ($n\geq 6$)}

\end{picture}

\end{center}

\begin{center}

\begin{picture}(0,0)(75,107) 
\setlength{\unitlength}{9mm}

\drawpolygon(0,0)(8,0)(8,4)(0,4)

\drawpolygon(10,0)(14,0)(14,4)(10,4)

\drawline[AHnb=0](0,2)(8,2)
\drawline[AHnb=0](2,0)(2,4)
\drawline[AHnb=0](4,0)(4,4)
\drawline[AHnb=0](6,0)(6,4)
\drawline[AHnb=0](12,0)(12,4)
\drawline[AHnb=0](10,2)(14,2)

\drawline[AHnb=0](.5,0)(2.5,4)
\drawline[AHnb=0](2.5,0)(4.5,4)
\drawline[AHnb=0](4.5,0)(6.5,4)

\drawline[AHnb=0](10.5,0)(12.5,4)

\drawline[AHnb=0](1.5,4)(3.5,0)
\drawline[AHnb=0](3.5,4)(5.5,0)
\drawline[AHnb=0](5.5,4)(7.5,0)
\drawline[AHnb=0](11.5,4)(13.5,0)

\drawline[AHnb=0](0,3)(1.5,0)
\drawline[AHnb=0](10,3)(11.5,0)

\drawline[AHnb=0](6.5,0)(8,3)
\drawline[AHnb=0](12.5,0)(14,3)

\drawline[AHnb=0](0,3)(.5,4)
\drawline[AHnb=0](10,3)(10.5,4)

\drawline[AHnb=0](8,3)(7.5,4)
\drawline[AHnb=0](14,3)(13.5,4)

%
%

\put(.4,-.2){\scriptsize $a_2$}
\put(2.4,-.2){\scriptsize $a_4$}
\put(4.4,-.2){\scriptsize $a_6$}
\put(6.4,-.2){\scriptsize $a_8$}
\put(10.15,-.2){\scriptsize $a_{2n-2}$}
\put(12.4,-.2){\scriptsize $a_{2n}$}

\put(1.35,-.3){\scriptsize $d_4$}
\put(3.35,-.3){\scriptsize $d_6$}
\put(5.35,-.3){\scriptsize $d_8$}
\put(7.35,-.3){\scriptsize $d_{10}$}
\put(11.35,-.3){\scriptsize $d_{2n}$}
\put(13.35,-.3){\scriptsize $d_2$}

\put(.05,4.15){\scriptsize $a_{2n-4}$}
\put(2.05,4.15){\scriptsize $a_{2n-2}$}
\put(4.15,4.15){\scriptsize $a_{2n}$}
\put(6.25,4.15){\scriptsize $a_2$}
\put(10.05,4.15){\scriptsize $a_{2n-8}$}
\put(12.15,4.15){\scriptsize $a_{2n-6}$}

\put(.6,1){\scriptsize $b_3$}
\put(2.6,1){\scriptsize $b_5$}
\put(4.6,1){\scriptsize $b_7$}
\put(6.6,1){\scriptsize $b_9$}

\put(10.05,1){\scriptsize $b_{2n-1}$}
\put(12.6,1){\scriptsize $b_1$}

\put(1,4.15){\scriptsize $d_{2n-2}$}
\put(3.3,4.15){\scriptsize $d_{2n}$}
\put(5.25,4.15){\scriptsize $d_2$}
\put(7.25,4.15){\scriptsize $d_4$}
\put(11.1,4.15){\scriptsize $d_{2n-6}$}
\put(13.1,4.15){\scriptsize $d_{2n-4}$}

\put(.1,2.9){\scriptsize $c_{2n-4}$}
\put(2.1,2.9){\scriptsize $c_{2n-2}$}
\put(4.1,2.9){\scriptsize $c_{2n}$}
\put(6.1,2.9){\scriptsize $c_2$}
\put(8.1,2.9){\scriptsize $c_4$}
\put(10.1,2.9){\scriptsize $c_{2n-8}$}
\put(12.1,2.9){\scriptsize $c_{2n-6}$}
\put(14.1,2.9){\scriptsize $c_{2n-4}$}

\put(.6,2.1){\scriptsize $e_1$}
\put(2.6,2.1){\scriptsize $e_3$}
\put(4.6,2.1){\scriptsize $e_5$}
\put(6.6,2.1){\scriptsize $e_7$}
\put(10.5,2.15){\scriptsize $e_{2n-3}$}
\put(12.5,2.15){\scriptsize $e_{2n-1}$}

\put(1.4,1.65){\scriptsize $f_3$}
\put(3.4,1.65){\scriptsize $f_5$}
\put(5.4,1.65){\scriptsize $f_7$}
\put(7.4,1.65){\scriptsize $f_9$}
\put(11.1,1.7){\scriptsize $f_{2n-1}$}
\put(13.4,1.65){\scriptsize $f_1$}

\put(15,1.9){\scriptsize $(n\geq 6)$}

\put(8.5,1.8){\scriptsize $\cdots$}

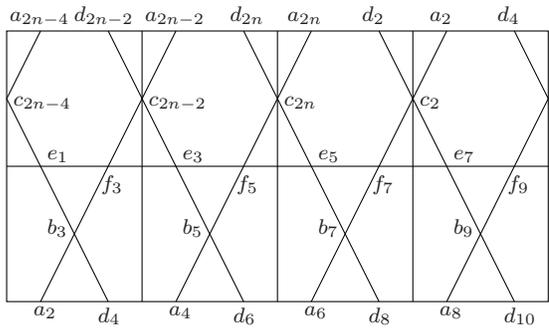
\put(.5,-.8){\scriptsize{\bf Fig.(4.9)\,: SEM of type\,-$(3, 6, 3, 6)$}}
\end{picture}

\end{center}

\vspace{11.3cm}

\hrule

For the remaining SEMs of types\,: $(3, 12^2)$, $(3, 4, 6, 4)$ and $(3^4, 6)$ we apply the same operations to
the equivelar maps of type $(6^3)$, which has 2$n$ vertices and $n$ faces (Fig.(4.3)). This can be seen as follows\,:

If we truncate an equivelar map along its vertices we get a SEM of type $(3, 12^2)$ with $6n$ vertices (see Fig.(4.10)). Now
applying subdivision to the SEM of type $(3, 12^2)$ one can obtain a SEM of type $(3,4,6,4)$ by the following three steps\,: in
the first step replace each edge of triangular faces by a quadrangular face, then in the second step replace each new
edge of these quadrangular faces which is opposite to the edge common in triangular face and quadrangular face by a triangular face.
Finally, in the last step constructing quadrangular faces whose opposite edges are one of the three edges of the triangular faces
(obtained during the second step) we get a SEM of type $(3,4,6,4)$ (see Fig.(4.11)). In this construction we see that
inside each truncated face of the equivelar map there are exactly (6+12) vertices (corresponding to hexagonal face and 12-gonal face).
Thus if a SEM of type (3,4,6,4) is obtained from the equivelar map then it has exactly 3.2$n$ + 18.$n$ vertices.
Moreover, subdividing each quadrangular face of this SEM of type $(3,4,6,4)$ by a diagonal, such that exactly four triangular
faces incident to each vertex, leads to a SEM of type $(3^4, 6)$ (see Fig.(4.12)) with 3.2$n$ + 18.$n$ vertices.

\smallskip

\hrule

\begin{center}

\begin{picture}(0,0)(75,35)
\setlength{\unitlength}{7mm}

\drawpolygon(0,0)(1.5,-1)(3,0)(3,2)(1.5,3)(0,2)
\drawpolygon(3,0)(4.5,-1)(6,0)(6,2)(4.5,3)(3,2)
\drawpolygon(6,0)(7.5,-1)(9,0)(9,2)(7.5,3)(6,2)
\drawpolygon(9,0)(10.5,-1)(12,0)(12,2)(10.5,3)(9,2)
\drawpolygon(14,0)(15.5,-1)(17,0)(17,2)(15.5,3)(14,2)

\put(13,1){\scriptsize $\cdots$}
\put(-.9,5.1){\scriptsize {\bf Infinite series of SEMs on torus}}

\drawline[AHnb=0](.6,-.4)(0,.5)
\drawline[AHnb=0](3.6,-.4)(3,.5)
\drawline[AHnb=0](6.6,-.4)(6,.5)
\drawline[AHnb=0](9.6,-.4)(9,.5)
\drawline[AHnb=0](14.6,-.4)(14,.5)
\drawline[AHnb=0](17.6,-.4)(17,.5)

\drawline[AHnb=0](0,1.4)(.6,2.4)
\drawline[AHnb=0](3,1.4)(3.6,2.4)
\drawline[AHnb=0](6,1.4)(6.6,2.4)
\drawline[AHnb=0](9,1.4)(9.6,2.4)
\drawline[AHnb=0](14,1.4)(14.6,2.4)
\drawline[AHnb=0](17,1.4)(17.6,2.4)

\drawline[AHnb=0](.9,2.6)(2.1,2.6)
\drawline[AHnb=0](3.9,2.6)(5.1,2.6)
\drawline[AHnb=0](6.9,2.6)(8.1,2.6)
\drawline[AHnb=0](9.9,2.6)(11.1,2.6)
\drawline[AHnb=0](14.9,2.6)(16.1,2.6)

\drawline[AHnb=0](.9,-.6)(2.1,-.6)
\drawline[AHnb=0](3.9,-.6)(5.1,-.6)
\drawline[AHnb=0](6.9,-.6)(8.1,-.6)
\drawline[AHnb=0](9.9,-.6)(11.1,-.6)
\drawline[AHnb=0](14.9,-.6)(16.1,-.6)

\drawline[AHnb=0](2.4,-.4)(3,.5)
\drawline[AHnb=0](5.4,-.4)(6,.5)
\drawline[AHnb=0](8.4,-.4)(9,.5)
\drawline[AHnb=0](11.4,-.4)(12,.5)
\drawline[AHnb=0](16.4,-.4)(17,.5)

\drawline[AHnb=0](3,1.4)(2.4,2.4)
\drawline[AHnb=0](6,1.4)(5.4,2.4)
\drawline[AHnb=0](9,1.4)(8.4,2.4)
\drawline[AHnb=0](12,1.4)(11.4,2.4)
\drawline[AHnb=0](17,1.4)(16.4,2.4)

\drawline[AHnb=0](2.4,-.4)(3.6,-.4)
\drawline[AHnb=0](5.4,-.4)(6.6,-.4)
\drawline[AHnb=0](8.4,-.4)(9.6,-.4)
\drawline[AHnb=0](16.4,-.4)(17.6,-.4)

\drawline[AHnb=0](2.4,2.4)(3.6,2.4)
\drawline[AHnb=0](5.4,2.4)(6.6,2.4)
\drawline[AHnb=0](8.4,2.4)(9.6,2.4)
\drawline[AHnb=0](16.4,2.4)(17.6,2.4)

\drawline[AHnb=0](17,2)(18.5,3)
\drawline[AHnb=0](17,0)(18.5,-1)

\drawline[AHnb=0](1.5,-1)(1.5,-1.5)
\drawline[AHnb=0](4.5,-1)(4.5,-1.5)
\drawline[AHnb=0](7.5,-1)(7.5,-1.5)
\drawline[AHnb=0](10.5,-1)(10.5,-1.5)
\drawline[AHnb=0](15.5,-1)(15.5,-1.5)

\drawline[AHnb=0](1.5,3)(1.5,4.2)
\drawline[AHnb=0](4.5,3)(4.5,4.2)
\drawline[AHnb=0](7.5,3)(7.5,4.2)
\drawline[AHnb=0](10.5,3)(10.5,4.2)
\drawline[AHnb=0](15.5,3)(15.5,4.2)

\drawline[AHnb=0](.9,2.6)(1.5,3.7)
\drawline[AHnb=0](3.9,2.6)(4.5,3.7)
\drawline[AHnb=0](6.9,2.6)(7.5,3.7)
\drawline[AHnb=0](9.9,2.6)(10.5,3.7)
\drawline[AHnb=0](14.9,2.6)(15.5,3.7)

\drawline[AHnb=0](2.1,2.6)(1.5,3.7)
\drawline[AHnb=0](5.1,2.6)(4.5,3.7)
\drawline[AHnb=0](8.1,2.6)(7.5,3.7)
\drawline[AHnb=0](11.1,2.6)(10.5,3.7)
\drawline[AHnb=0](16.1,2.6)(15.5,3.7)



%
%
%

%
\put(-.4,.45){\scriptsize {\tiny $a_1$}}
\put(2.5,.55){\scriptsize {\tiny $a_3$}}
\put(5.5,.55){\scriptsize {\tiny $a_5$}}
\put(8.5,.55){\scriptsize {\tiny $a_7$}}
\put(11.5,.55){\scriptsize {\tiny $a_9$}}

\put(14.1,.55){\scriptsize {\tiny $a_{2n-1}$}}
\put(17.1,.55){\scriptsize {\tiny $a_1$}}

\put(.1,1.2){\scriptsize {\tiny $a_{2n-4}$}}
\put(3.1,1.2){\scriptsize {\tiny $a_{2n-2}$}}
\put(6.1,1.2){\scriptsize {\tiny $a_{2n}$}}
\put(9.1,1.2){\scriptsize {\tiny $a_2$}}
\put(11.5,1.2){\scriptsize {\tiny $a_4$}}

\put(14.05,1.25){\scriptsize {\tiny $a_{2n-6}$}}
\put(17.1,1.2){\scriptsize {\tiny $a_{2n-4}$}}

\put(.6,-.3){\scriptsize {\tiny $b_1$}}
\put(3.6,-.3){\scriptsize {\tiny $b_3$}}
\put(6.6,-.3){\scriptsize {\tiny $b_5$}}
\put(9.6,-.3){\scriptsize {\tiny $b_7$}}
\put(14.65,-.35){\scriptsize {\tiny $b_{2n-1}$}}

\put(2.1,-.2){\scriptsize {\tiny $c_3$}}
\put(5.1,-.2){\scriptsize {\tiny $c_5$}}
\put(8.1,-.2){\scriptsize {\tiny $c_7$}}
\put(11.1,-.2){\scriptsize {\tiny $c_9$}}
\put(16.05,-.3){\scriptsize {\tiny $c_1$}}
\put(17.65,-.3){\scriptsize {\tiny $b_1$}}

\put(-.2,2.65){\scriptsize {\tiny $c_{2n-3}$}}
\put(2.95,2.8){\scriptsize {\tiny $c_{2n-1}$}}
\put(6.6,2.8){\scriptsize {\tiny $c_1$}}
\put(9.6,2.8){\scriptsize {\tiny $c_3$}}
\put(13.9,2.8){\scriptsize {\tiny $c_{2n-5}$}}

\put(2.15,2.6){\scriptsize {\tiny $b_{2n-3}$}}
\put(5.1,2.7){\scriptsize {\tiny $b_{2n-1}$}}
\put(8.1,2.7){\scriptsize {\tiny $b_1$}}
\put(11.1,2.7){\scriptsize {\tiny $b_3$}}
\put(16.1,2.7){\scriptsize {\tiny $b_{2n-5}$}}

\put(-.5,2.45){\scriptsize {\tiny $b_{2n-4}$}}
\put(3.6,2.15){\scriptsize {\tiny $b_{2n-2}$}}
\put(6.6,2.15){\scriptsize {\tiny $b_{2n}$}}
\put(9.6,2.15){\scriptsize {\tiny $b_2$}}
\put(13.5,2.45){\scriptsize {\tiny $b_{2n-6}$}}

\put(1.4,2.25){\scriptsize {\tiny $c_{2n-2}$}}
\put(4.85,2.25){\scriptsize {\tiny $c_{2n}$}}
\put(8.05,2.25){\scriptsize {\tiny $c_2$}}
\put(11.05,2.25){\scriptsize {\tiny $c_4$}}
\put(15.4,2.3){\scriptsize {\tiny $c_{2n-4}$}}
\put(17.7,2.2){\scriptsize {\tiny $b_{2n-4}$}}

\put(.6,-.8){\scriptsize {\tiny $c_2$}}
\put(3.6,-.8){\scriptsize {\tiny $c_4$}}
\put(6.6,-.8){\scriptsize {\tiny $c_6$}}
\put(9.6,-.8){\scriptsize {\tiny $c_8$}}
\put(14.35,-.7){\scriptsize {\tiny $c_{2n}$}}

\put(2.15,-.8){\scriptsize {\tiny $b_2$}}
\put(5.15,-.8){\scriptsize {\tiny $b_4$}}
\put(8.15,-.8){\scriptsize {\tiny $b_6$}}
\put(11.15,-.8){\scriptsize {\tiny $b_8$}}
\put(16.2,-.8){\scriptsize {\tiny $b_{2n}$}}

\put(1.6,3.8){\scriptsize {\tiny $a_{2n-3}$}}
\put(4.6,3.8){\scriptsize {\tiny $a_{2n-1}$}}
\put(7.6,3.8){\scriptsize {\tiny $a_1$}}
\put(10.6,3.8){\scriptsize {\tiny $a_3$}}
\put(15.6,3.8){\scriptsize {\tiny $a_{2n-5}$}}

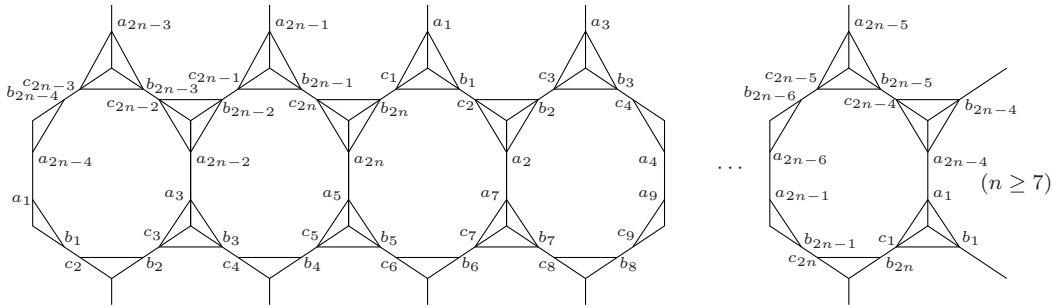
\put(.5,-2.2){\scriptsize{\bf Fig.(4.10)\,: SEM of type\,-$(3, 12^2)$}}

\put(18,.7){\scriptsize $(n \geq 7)$}

\end{picture}

\end{center}

\begin{center}

\begin{picture}(0,0)(75,75)
\setlength{\unitlength}{12mm}

\drawpolygon(0,0)(1,-.7)(2,0)(2,1)(1,1.7)(0,1)
\drawpolygon(2,0)(3,-.7)(4,0)(4,1)(3,1.7)(2,1)
\drawpolygon(4,0)(5,-.7)(6,0)(6,1)(5,1.7)(4,1)
\drawpolygon(6,0)(7,-.7)(8,0)(8,1)(7,1.7)(6,1)

\drawpolygon(9,0)(10,-.7)(11,0)(11,1)(10,1.7)(9,1)

\drawline[AHnb=0](.7,-.5)(1.3,-.5)

\drawline[AHnb=0](2.7,-.5)(3.3,-.5)
\drawline[AHnb=0](4.7,-.5)(5.3,-.5)
\drawline[AHnb=0](6.7,-.5)(7.3,-.5)
\drawline[AHnb=0](9.7,-.5)(10.3,-.5)

\drawline[AHnb=0](.7,1.5)(1.3,1.5)
\drawline[AHnb=0](2.7,1.5)(3.3,1.5)
\drawline[AHnb=0](4.7,1.5)(5.3,1.5)
\drawline[AHnb=0](6.7,1.5)(7.3,1.5)
\drawline[AHnb=0](9.7,1.5)(10.3,1.5)

\drawline[AHnb=0](1.7,1.2)(2.3,1.2)
\drawline[AHnb=0](3.7,1.2)(4.3,1.2)
\drawline[AHnb=0](5.7,1.2)(6.3,1.2)
\drawline[AHnb=0](10.7,1.2)(11.3,1.2)

\drawline[AHnb=0](1.7,-.2)(2.3,-.2)
\drawline[AHnb=0](3.7,-.2)(4.3,-.2)
\drawline[AHnb=0](5.7,-.2)(6.3,-.2)
\drawline[AHnb=0](10.7,-.2)(11.3,-.2)

\drawline[AHnb=0](11,0)(12,-.7)
\drawline[AHnb=0](11,1)(12,1.7)

\drawline[AHnb=0](1.7,-.2)(2,.25)
\drawline[AHnb=0](3.7,-.2)(4,.25)
\drawline[AHnb=0](5.7,-.2)(6,.25)
\drawline[AHnb=0](7.7,-.2)(8,.25)
\drawline[AHnb=0](10.7,-.2)(11,.25)

\drawline[AHnb=0](.3,-.2)(0,.25)
\drawline[AHnb=0](2.3,-.2)(2,.25)
\drawline[AHnb=0](4.3,-.2)(4,.25)
\drawline[AHnb=0](6.3,-.2)(6,.25)
\drawline[AHnb=0](9.3,-.2)(9,.25)
\drawline[AHnb=0](11.3,-.2)(11,.25)

\drawline[AHnb=0](1.7,1.2)(2,.75)
\drawline[AHnb=0](3.7,1.2)(4,.75)
\drawline[AHnb=0](5.7,1.2)(6,.75)
\drawline[AHnb=0](7.7,1.2)(8,.75)
\drawline[AHnb=0](10.7,1.2)(11,.75)

\drawline[AHnb=0](0.3,1.2)(0,.75)
\drawline[AHnb=0](2.3,1.2)(2,.75)
\drawline[AHnb=0](4.3,1.2)(4,.75)
\drawline[AHnb=0](6.3,1.2)(6,.75)
\drawline[AHnb=0](9.3,1.2)(9,.75)
\drawline[AHnb=0](11.3,1.2)(11,.75)

\drawpolygon(.5,1)(.8,1.15)(1.2,1.15)(1.5,1)(1.75,.65)(1.75,.35)(1.5,.05)(1.2,-.1)(.8, -.1)(.5,.05)(.25,.35)(.25,.65)

\drawpolygon(1,.9)(1.35,.7)(1.35,.4)(1,.2)(.65,.4)(.65,.7)

\drawline[AHnb=0](.8,1.15)(1,.9)
\drawline[AHnb=0](1,.9)(1.2,1.15)

\drawline[AHnb=0](1.35,.7)(1.5,1)
\drawline[AHnb=0](1.35,.7)(1.75,.65)

\drawline[AHnb=0](1.35,.4)(1.75,.35)
\drawline[AHnb=0](1.35,.4)(1.5,.05)

\drawline[AHnb=0](1,.2)(1.2,-.1)
\drawline[AHnb=0](1,.2)(.8,-.1)

\drawline[AHnb=0](.65,.4)(.25,.35)
\drawline[AHnb=0](.65,.4)(.5,.05)

\drawline[AHnb=0](.65,.7)(.5,1)
\drawline[AHnb=0](.65,.7)(.25,.65)

\drawline[AHnb=0](.8,1.15)(.72,1.5)
\drawline[AHnb=0](1.2,1.15)(1.3,1.5)

\drawline[AHnb=0](.5,1)(.3,1.2)
\drawline[AHnb=0](.25,.65)(0,.75)

\drawline[AHnb=0](.5,.05)(.3,-.2)
\drawline[AHnb=0](.25,.35)(0,.25)

\drawline[AHnb=0](.8,-.1)(.7,-.5)
\drawline[AHnb=0](1.2,-.1)(1.3,-.5)

\drawline[AHnb=0](1.5,.05)(1.7,-.2)
\drawline[AHnb=0](1.75,.35)(2,.25)

\drawline[AHnb=0](1.75,.65)(2,.75)
\drawline[AHnb=0](1.5,1)(1.7,1.2)

\drawpolygon(2.5,1)(2.8,1.15)(3.2,1.15)(3.5,1)(3.75,.65)(3.75,.35)(3.5,.05)(3.2,-.1)(2.8, -.1)(2.5,.05)(2.25,.35)(2.25,.65)

\drawpolygon(3,.9)(3.35,.7)(3.35,.4)(3,.2)(2.65,.4)(2.65,.7)

\drawline[AHnb=0](2.8,1.15)(3,.9)
\drawline[AHnb=0](3,.9)(3.2,1.15)

\drawline[AHnb=0](3.35,.7)(3.5,1)
\drawline[AHnb=0](3.35,.7)(3.75,.65)

\drawline[AHnb=0](3.35,.4)(3.75,.35)
\drawline[AHnb=0](3.35,.4)(3.5,.05)

\drawline[AHnb=0](3,.2)(3.2,-.1)
\drawline[AHnb=0](3,.2)(2.8,-.1)

\drawline[AHnb=0](2.65,.4)(2.25,.35)
\drawline[AHnb=0](2.65,.4)(2.5,.05)

\drawline[AHnb=0](2.65,.7)(2.5,1)
\drawline[AHnb=0](2.65,.7)(2.25,.65)

\drawline[AHnb=0](2.8,1.15)(2.72,1.5)
\drawline[AHnb=0](3.2,1.15)(3.3,1.5)

\drawline[AHnb=0](2.5,1)(2.3,1.2)
\drawline[AHnb=0](2.25,.65)(2,.75)

\drawline[AHnb=0](2.5,.05)(2.3,-.2)
\drawline[AHnb=0](2.25,.35)(2,.25)

\drawline[AHnb=0](2.8,-.1)(2.7,-.5)
\drawline[AHnb=0](3.2,-.1)(3.3,-.5)

\drawline[AHnb=0](3.5,.05)(3.7,-.2)
\drawline[AHnb=0](3.75,.35)(4,.25)

\drawline[AHnb=0](3.75,.65)(4,.75)
\drawline[AHnb=0](3.5,1)(3.7,1.2)

\drawpolygon(4.5,1)(4.8,1.15)(5.2,1.15)(5.5,1)(5.75,.65)(5.75,.35)(5.5,.05)(5.2,-.1)(4.8, -.1)(4.5,.05)(4.25,.35)(4.25,.65)

\drawpolygon(5,.9)(5.35,.7)(5.35,.4)(5,.2)(4.65,.4)(4.65,.7)

\drawline[AHnb=0](4.8,1.15)(5,.9)
\drawline[AHnb=0](5,.9)(5.2,1.15)

\drawline[AHnb=0](5.35,.7)(5.5,1)
\drawline[AHnb=0](5.35,.7)(5.75,.65)

\drawline[AHnb=0](5.35,.4)(5.75,.35)
\drawline[AHnb=0](5.35,.4)(5.5,.05)

\drawline[AHnb=0](5,.2)(5.2,-.1)
\drawline[AHnb=0](5,.2)(4.8,-.1)

\drawline[AHnb=0](4.65,.4)(4.25,.35)
\drawline[AHnb=0](4.65,.4)(4.5,.05)

\drawline[AHnb=0](4.65,.7)(4.5,1)
\drawline[AHnb=0](4.65,.7)(4.25,.65)

\drawline[AHnb=0](4.8,1.15)(4.72,1.5)
\drawline[AHnb=0](5.2,1.15)(5.3,1.5)

\drawline[AHnb=0](4.5,1)(4.3,1.2)
\drawline[AHnb=0](4.25,.65)(4,.75)

\drawline[AHnb=0](4.5,.05)(4.3,-.2)
\drawline[AHnb=0](4.25,.35)(4,.25)

\drawline[AHnb=0](4.8,-.1)(4.7,-.5)
\drawline[AHnb=0](5.2,-.1)(5.3,-.5)

\drawline[AHnb=0](5.5,.05)(5.7,-.2)
\drawline[AHnb=0](5.75,.35)(6,.25)

\drawline[AHnb=0](5.75,.65)(6,.75)
\drawline[AHnb=0](5.5,1)(5.7,1.2)

\drawpolygon(6.5,1)(6.8,1.15)(7.2,1.15)(7.5,1)(7.75,.65)(7.75,.35)(7.5,.05)(7.2,-.1)(6.8, -.1)(6.5,.05)(6.25,.35)(6.25,.65)

\drawpolygon(7,.9)(7.35,.7)(7.35,.4)(7,.2)(6.65,.4)(6.65,.7)

\drawline[AHnb=0](6.8,1.15)(7,.9)
\drawline[AHnb=0](7,.9)(7.2,1.15)

\drawline[AHnb=0](7.35,.7)(7.5,1)
\drawline[AHnb=0](7.35,.7)(7.75,.65)

\drawline[AHnb=0](7.35,.4)(7.75,.35)
\drawline[AHnb=0](7.35,.4)(7.5,.05)

\drawline[AHnb=0](7,.2)(7.2,-.1)
\drawline[AHnb=0](7,.2)(6.8,-.1)

\drawline[AHnb=0](6.65,.4)(6.25,.35)
\drawline[AHnb=0](6.65,.4)(6.5,.05)

\drawline[AHnb=0](6.65,.7)(6.5,1)
\drawline[AHnb=0](6.65,.7)(6.25,.65)

\drawline[AHnb=0](6.8,1.15)(6.72,1.5)
\drawline[AHnb=0](7.2,1.15)(7.3,1.5)

\drawline[AHnb=0](6.5,1)(6.3,1.2)
\drawline[AHnb=0](6.25,.65)(6,.75)

\drawline[AHnb=0](6.5,.05)(6.3,-.2)
\drawline[AHnb=0](6.25,.35)(6,.25)

\drawline[AHnb=0](6.8,-.1)(6.7,-.5)
\drawline[AHnb=0](7.2,-.1)(7.3,-.5)

\drawline[AHnb=0](7.5,.05)(7.7,-.2)
\drawline[AHnb=0](7.75,.35)(8,.25)

\drawline[AHnb=0](7.75,.65)(8,.75)
\drawline[AHnb=0](7.5,1)(7.7,1.2)

\drawpolygon(9.5,1)(9.8,1.15)(10.2,1.15)(10.5,1)(10.75,.65)(10.75,.35)(10.5,.05)(10.2,-.1)(9.8, -.1)(9.5,.05)(9.25,.35)(9.25,.65)

\drawpolygon(10,.9)(10.35,.7)(10.35,.4)(10,.2)(9.65,.4)(9.65,.7)

\drawline[AHnb=0](9.8,1.15)(10,.9)
\drawline[AHnb=0](10,.9)(10.2,1.15)

\drawline[AHnb=0](10.35,.7)(10.5,1)
\drawline[AHnb=0](10.35,.7)(10.75,.65)

\drawline[AHnb=0](10.35,.4)(10.75,.35)
\drawline[AHnb=0](10.35,.4)(10.5,.05)

\drawline[AHnb=0](10,.2)(10.2,-.1)
\drawline[AHnb=0](10,.2)(9.8,-.1)

\drawline[AHnb=0](9.65,.4)(9.25,.35)
\drawline[AHnb=0](9.65,.4)(9.5,.05)

\drawline[AHnb=0](9.65,.7)(9.5,1)
\drawline[AHnb=0](9.65,.7)(9.25,.65)

\drawline[AHnb=0](9.8,1.15)(9.72,1.5)
\drawline[AHnb=0](10.2,1.15)(10.3,1.5)

\drawline[AHnb=0](9.5,1)(9.3,1.2)
\drawline[AHnb=0](9.25,.65)(9,.75)

\drawline[AHnb=0](9.5,.05)(9.3,-.2)
\drawline[AHnb=0](9.25,.35)(9,.25)

\drawline[AHnb=0](9.8,-.1)(9.7,-.5)
\drawline[AHnb=0](10.2,-.1)(10.3,-.5)

\drawline[AHnb=0](10.5,.05)(10.7,-.2)
\drawline[AHnb=0](10.75,.35)(11,.25)

\drawline[AHnb=0](10.75,.65)(11,.75)
\drawline[AHnb=0](10.5,1)(10.7,1.2)

\drawline[AHnb=0](1,1.7)(1,2.3)
\drawline[AHnb=0](3,1.7)(3,2.3)
\drawline[AHnb=0](5,1.7)(5,2.3)
\drawline[AHnb=0](7,1.7)(7,2.3)
\drawline[AHnb=0](10,1.7)(10,2.3)

\drawline[AHnb=0](.7,1.5)(1,2)
\drawline[AHnb=0](2.7,1.5)(3,2)
\drawline[AHnb=0](4.7,1.5)(5,2)
\drawline[AHnb=0](6.7,1.5)(7,2)
\drawline[AHnb=0](9.7,1.5)(10,2)

\drawline[AHnb=0](1.3,1.5)(1,2)
\drawline[AHnb=0](3.3,1.5)(3,2)
\drawline[AHnb=0](5.3,1.5)(5,2)
\drawline[AHnb=0](7.3,1.5)(7,2)
\drawline[AHnb=0](10.3,1.5)(10,2)

\drawline[AHnb=0](1,-.7)(1,-1)
\drawline[AHnb=0](3,-.7)(3,-1)
\drawline[AHnb=0](5,-.7)(5,-1)
\drawline[AHnb=0](7,-.7)(7,-1)
\drawline[AHnb=0](10,-.7)(10,-1)


%


\put(.18,-.34){\scriptsize {\tiny $b_1$}}
\put(2.18,-.34){\scriptsize {\tiny $b_3$}}
\put(4.18,-.34){\scriptsize {\tiny $b_5$}}
\put(6.18,-.34){\scriptsize {\tiny $b_7$}}
\put(8.8,-.34){\scriptsize {\tiny $b_{n-1}$}}
\put(11.18,-.34){\scriptsize {\tiny $b_1$}}

\put(1.3,-.6){\scriptsize {\tiny $b_2$}}
\put(3.3,-.6){\scriptsize {\tiny $b_4$}}
\put(5.3,-.6){\scriptsize {\tiny $b_6$}}
\put(7.3,-.6){\scriptsize {\tiny $b_8$}}
\put(10.3,-.6){\scriptsize {\tiny $b_{2n}$}}

\put(.55,-.6){\scriptsize {\tiny $c_2$}}
\put(2.55,-.6){\scriptsize {\tiny $c_4$}}
\put(4.55,-.6){\scriptsize {\tiny $c_6$}}
\put(6.55,-.6){\scriptsize {\tiny $c_8$}}
\put(9.45,-.6){\scriptsize {\tiny $c_{2n}$}}

\put(1.7,-.34){\scriptsize {\tiny $c_3$}}
\put(3.7,-.34){\scriptsize {\tiny $c_5$}}
\put(5.7,-.34){\scriptsize {\tiny $c_7$}}
\put(7.7,-.34){\scriptsize {\tiny $c_9$}}
\put(10.7,-.34){\scriptsize {\tiny $c_1$}}

\put(-.22,.2){\scriptsize {\tiny $a_1$}}
\put(11.1,.2){\scriptsize {\tiny $a_1$}}

\put(-.65,.75){\scriptsize {\tiny $a_{2n-4}$}}
\put(11.1,.7){\scriptsize {\tiny $a_{2n-4}$}}

\put(1.3,1.54){\scriptsize {\tiny $b_{2n-3}$}}
\put(3.3,1.54){\scriptsize {\tiny $b_{2n-1}$}}
\put(5.3,1.54){\scriptsize {\tiny $b_{1}$}}
\put(7.3,1.54){\scriptsize {\tiny $b_{3}$}}
\put(10.3,1.54){\scriptsize {\tiny $b_{2n-5}$}}

\put(.07,1.54){\scriptsize {\tiny $c_{2n-3}$}}
\put(2.07,1.54){\scriptsize {\tiny $c_{2n-1}$}}
\put(4.47,1.54){\scriptsize {\tiny $c_{1}$}}
\put(6.47,1.54){\scriptsize {\tiny $c_{3}$}}
\put(9.07,1.54){\scriptsize {\tiny $c_{2n-5}$}}

\put(-.35,1.25){\scriptsize {\tiny $b_{2n-4}$}}
\put(2.15,1.3){\scriptsize {\tiny $b_{2n-2}$}}
\put(4.05,1.3){\scriptsize {\tiny $b_{2n}$}}
\put(6.05,1.3){\scriptsize {\tiny $b_{2}$}}
\put(8.7,1.3){\scriptsize {\tiny $b_{2n-6}$}}
\put(11.35,1.1){\scriptsize {\tiny $b_{2n-4}$}}

\put(1.05,2){\scriptsize {\tiny $a_{2n-3}$}}
\put(3.05,2){\scriptsize {\tiny $a_{2n-1}$}}
\put(5.05,2){\scriptsize {\tiny $a_1$}}
\put(7.05,2){\scriptsize {\tiny $a_3$}}
\put(10.05,2){\scriptsize {\tiny $a_{2n-5}$}}

\put(1.4,1.3){\scriptsize {\tiny $c_{2n-2}$}}
\put(3.65,1.3){\scriptsize {\tiny $c_{2n}$}}
\put(5.65,1.3){\scriptsize {\tiny $c_{2}$}}
\put(7.65,1.3){\scriptsize {\tiny $c_{4}$}}
\put(10.4,1.3){\scriptsize {\tiny $c_{2n-4}$}}


\put(8.3,.5){\scriptsize $\cdots$}

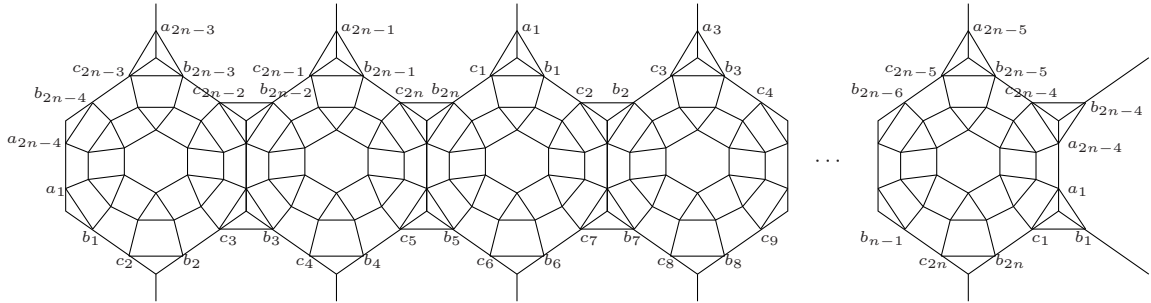
\put(.5,-1.5){\scriptsize{\bf Fig.(4.11)\,: SEM of type\,-$(3, 4, 6, 4)$, $n \geq 7$}}
\end{picture}

\end{center}

\begin{center}

\begin{picture}(0,0)(75,117)
\setlength{\unitlength}{12mm}

\drawpolygon(0,0)(1,-.7)(2,0)(2,1)(1,1.7)(0,1)
\drawpolygon(2,0)(3,-.7)(4,0)(4,1)(3,1.7)(2,1)
\drawpolygon(4,0)(5,-.7)(6,0)(6,1)(5,1.7)(4,1)
\drawpolygon(6,0)(7,-.7)(8,0)(8,1)(7,1.7)(6,1)

\drawpolygon(9,0)(10,-.7)(11,0)(11,1)(10,1.7)(9,1)

\drawline[AHnb=0](.7,-.5)(1.3,-.5)

\drawline[AHnb=0](2.7,-.5)(3.3,-.5)
\drawline[AHnb=0](4.7,-.5)(5.3,-.5)
\drawline[AHnb=0](6.7,-.5)(7.3,-.5)
\drawline[AHnb=0](9.7,-.5)(10.3,-.5)

\drawline[AHnb=0](.7,1.5)(1.3,1.5)
\drawline[AHnb=0](2.7,1.5)(3.3,1.5)
\drawline[AHnb=0](4.7,1.5)(5.3,1.5)
\drawline[AHnb=0](6.7,1.5)(7.3,1.5)
\drawline[AHnb=0](9.7,1.5)(10.3,1.5)

\drawline[AHnb=0](1.7,1.2)(2.3,1.2)
\drawline[AHnb=0](3.7,1.2)(4.3,1.2)
\drawline[AHnb=0](5.7,1.2)(6.3,1.2)
\drawline[AHnb=0](10.7,1.2)(11.3,1.2)

\drawline[AHnb=0](1.7,-.2)(2.3,-.2)
\drawline[AHnb=0](3.7,-.2)(4.3,-.2)
\drawline[AHnb=0](5.7,-.2)(6.3,-.2)
\drawline[AHnb=0](10.7,-.2)(11.3,-.2)

\drawline[AHnb=0](11,0)(11.7,-.5)
\drawline[AHnb=0](11,1)(11.7,1.5)

\drawline[AHnb=0](1.7,-.2)(2,.25)
\drawline[AHnb=0](3.7,-.2)(4,.25)
\drawline[AHnb=0](5.7,-.2)(6,.25)
\drawline[AHnb=0](7.7,-.2)(8,.25)
\drawline[AHnb=0](10.7,-.2)(11,.25)

\drawline[AHnb=0](.3,-.2)(0,.25)
\drawline[AHnb=0](2.3,-.2)(2,.25)
\drawline[AHnb=0](4.3,-.2)(4,.25)
\drawline[AHnb=0](6.3,-.2)(6,.25)
\drawline[AHnb=0](9.3,-.2)(9,.25)
\drawline[AHnb=0](11.3,-.2)(11,.25)

\drawline[AHnb=0](1.7,1.2)(2,.75)
\drawline[AHnb=0](3.7,1.2)(4,.75)
\drawline[AHnb=0](5.7,1.2)(6,.75)
\drawline[AHnb=0](7.7,1.2)(8,.75)
\drawline[AHnb=0](10.7,1.2)(11,.75)

\drawline[AHnb=0](0.3,1.2)(0,.75)
\drawline[AHnb=0](2.3,1.2)(2,.75)
\drawline[AHnb=0](4.3,1.2)(4,.75)
\drawline[AHnb=0](6.3,1.2)(6,.75)
\drawline[AHnb=0](9.3,1.2)(9,.75)
\drawline[AHnb=0](11.3,1.2)(11,.75)

\drawpolygon(.5,1)(.8,1.15)(1.2,1.15)(1.5,1)(1.75,.65)(1.75,.35)(1.5,.05)(1.2,-.1)(.8, -.1)(.5,.05)(.25,.35)(.25,.65)

\drawpolygon(1,.9)(1.35,.7)(1.35,.4)(1,.2)(.65,.4)(.65,.7)

\drawline[AHnb=0](.8,1.15)(1,.9)
\drawline[AHnb=0](1,.9)(1.2,1.15)

\drawline[AHnb=0](1.35,.7)(1.5,1)
\drawline[AHnb=0](1.35,.7)(1.75,.65)

\drawline[AHnb=0](1.35,.4)(1.75,.35)
\drawline[AHnb=0](1.35,.4)(1.5,.05)

\drawline[AHnb=0](1,.2)(1.2,-.1)
\drawline[AHnb=0](1,.2)(.8,-.1)

\drawline[AHnb=0](.65,.4)(.25,.35)
\drawline[AHnb=0](.65,.4)(.5,.05)

\drawline[AHnb=0](.65,.7)(.5,1)
\drawline[AHnb=0](.65,.7)(.25,.65)

\drawline[AHnb=0](.8,1.15)(.72,1.5)
\drawline[AHnb=0](1.2,1.15)(1.3,1.5)

\drawline[AHnb=0](.5,1)(.3,1.2)
\drawline[AHnb=0](.25,.65)(0,.75)

\drawline[AHnb=0](.5,.05)(.3,-.2)
\drawline[AHnb=0](.25,.35)(0,.25)

\drawline[AHnb=0](.8,-.1)(.7,-.5)
\drawline[AHnb=0](1.2,-.1)(1.3,-.5)

\drawline[AHnb=0](1.5,.05)(1.7,-.2)
\drawline[AHnb=0](1.75,.35)(2,.25)

\drawline[AHnb=0](1.75,.65)(2,.75)
\drawline[AHnb=0](1.5,1)(1.7,1.2)

\drawpolygon(2.5,1)(2.8,1.15)(3.2,1.15)(3.5,1)(3.75,.65)(3.75,.35)(3.5,.05)(3.2,-.1)(2.8, -.1)(2.5,.05)(2.25,.35)(2.25,.65)

\drawpolygon(3,.9)(3.35,.7)(3.35,.4)(3,.2)(2.65,.4)(2.65,.7)

\drawline[AHnb=0](2.8,1.15)(3,.9)
\drawline[AHnb=0](3,.9)(3.2,1.15)

\drawline[AHnb=0](3.35,.7)(3.5,1)
\drawline[AHnb=0](3.35,.7)(3.75,.65)

\drawline[AHnb=0](3.35,.4)(3.75,.35)
\drawline[AHnb=0](3.35,.4)(3.5,.05)

\drawline[AHnb=0](3,.2)(3.2,-.1)
\drawline[AHnb=0](3,.2)(2.8,-.1)

\drawline[AHnb=0](2.65,.4)(2.25,.35)
\drawline[AHnb=0](2.65,.4)(2.5,.05)

\drawline[AHnb=0](2.65,.7)(2.5,1)
\drawline[AHnb=0](2.65,.7)(2.25,.65)

\drawline[AHnb=0](2.8,1.15)(2.72,1.5)
\drawline[AHnb=0](3.2,1.15)(3.3,1.5)

\drawline[AHnb=0](2.5,1)(2.3,1.2)
\drawline[AHnb=0](2.25,.65)(2,.75)

\drawline[AHnb=0](2.5,.05)(2.3,-.2)
\drawline[AHnb=0](2.25,.35)(2,.25)

\drawline[AHnb=0](2.8,-.1)(2.7,-.5)
\drawline[AHnb=0](3.2,-.1)(3.3,-.5)

\drawline[AHnb=0](3.5,.05)(3.7,-.2)
\drawline[AHnb=0](3.75,.35)(4,.25)

\drawline[AHnb=0](3.75,.65)(4,.75)
\drawline[AHnb=0](3.5,1)(3.7,1.2)

\drawpolygon(4.5,1)(4.8,1.15)(5.2,1.15)(5.5,1)(5.75,.65)(5.75,.35)(5.5,.05)(5.2,-.1)(4.8, -.1)(4.5,.05)(4.25,.35)(4.25,.65)

\drawpolygon(5,.9)(5.35,.7)(5.35,.4)(5,.2)(4.65,.4)(4.65,.7)

\drawline[AHnb=0](4.8,1.15)(5,.9)
\drawline[AHnb=0](5,.9)(5.2,1.15)

\drawline[AHnb=0](5.35,.7)(5.5,1)
\drawline[AHnb=0](5.35,.7)(5.75,.65)

\drawline[AHnb=0](5.35,.4)(5.75,.35)
\drawline[AHnb=0](5.35,.4)(5.5,.05)

\drawline[AHnb=0](5,.2)(5.2,-.1)
\drawline[AHnb=0](5,.2)(4.8,-.1)

\drawline[AHnb=0](4.65,.4)(4.25,.35)
\drawline[AHnb=0](4.65,.4)(4.5,.05)

\drawline[AHnb=0](4.65,.7)(4.5,1)
\drawline[AHnb=0](4.65,.7)(4.25,.65)

\drawline[AHnb=0](4.8,1.15)(4.72,1.5)
\drawline[AHnb=0](5.2,1.15)(5.3,1.5)

\drawline[AHnb=0](4.5,1)(4.3,1.2)
\drawline[AHnb=0](4.25,.65)(4,.75)

\drawline[AHnb=0](4.5,.05)(4.3,-.2)
\drawline[AHnb=0](4.25,.35)(4,.25)

\drawline[AHnb=0](4.8,-.1)(4.7,-.5)
\drawline[AHnb=0](5.2,-.1)(5.3,-.5)

\drawline[AHnb=0](5.5,.05)(5.7,-.2)
\drawline[AHnb=0](5.75,.35)(6,.25)

\drawline[AHnb=0](5.75,.65)(6,.75)
\drawline[AHnb=0](5.5,1)(5.7,1.2)

\drawpolygon(6.5,1)(6.8,1.15)(7.2,1.15)(7.5,1)(7.75,.65)(7.75,.35)(7.5,.05)(7.2,-.1)(6.8, -.1)(6.5,.05)(6.25,.35)(6.25,.65)

\drawpolygon(7,.9)(7.35,.7)(7.35,.4)(7,.2)(6.65,.4)(6.65,.7)

\drawline[AHnb=0](6.8,1.15)(7,.9)
\drawline[AHnb=0](7,.9)(7.2,1.15)

\drawline[AHnb=0](7.35,.7)(7.5,1)
\drawline[AHnb=0](7.35,.7)(7.75,.65)

\drawline[AHnb=0](7.35,.4)(7.75,.35)
\drawline[AHnb=0](7.35,.4)(7.5,.05)

\drawline[AHnb=0](7,.2)(7.2,-.1)
\drawline[AHnb=0](7,.2)(6.8,-.1)

\drawline[AHnb=0](6.65,.4)(6.25,.35)
\drawline[AHnb=0](6.65,.4)(6.5,.05)

\drawline[AHnb=0](6.65,.7)(6.5,1)
\drawline[AHnb=0](6.65,.7)(6.25,.65)

\drawline[AHnb=0](6.8,1.15)(6.72,1.5)
\drawline[AHnb=0](7.2,1.15)(7.3,1.5)

\drawline[AHnb=0](6.5,1)(6.3,1.2)
\drawline[AHnb=0](6.25,.65)(6,.75)

\drawline[AHnb=0](6.5,.05)(6.3,-.2)
\drawline[AHnb=0](6.25,.35)(6,.25)

\drawline[AHnb=0](6.8,-.1)(6.7,-.5)
\drawline[AHnb=0](7.2,-.1)(7.3,-.5)

\drawline[AHnb=0](7.5,.05)(7.7,-.2)
\drawline[AHnb=0](7.75,.35)(8,.25)

\drawline[AHnb=0](7.75,.65)(8,.75)
\drawline[AHnb=0](7.5,1)(7.7,1.2)

\drawpolygon(9.5,1)(9.8,1.15)(10.2,1.15)(10.5,1)(10.75,.65)(10.75,.35)(10.5,.05)(10.2,-.1)(9.8, -.1)(9.5,.05)(9.25,.35)(9.25,.65)

\drawpolygon(10,.9)(10.35,.7)(10.35,.4)(10,.2)(9.65,.4)(9.65,.7)

\drawline[AHnb=0](9.8,1.15)(10,.9)
\drawline[AHnb=0](10,.9)(10.2,1.15)

\drawline[AHnb=0](10.35,.7)(10.5,1)
\drawline[AHnb=0](10.35,.7)(10.75,.65)

\drawline[AHnb=0](10.35,.4)(10.75,.35)
\drawline[AHnb=0](10.35,.4)(10.5,.05)

\drawline[AHnb=0](10,.2)(10.2,-.1)
\drawline[AHnb=0](10,.2)(9.8,-.1)

\drawline[AHnb=0](9.65,.4)(9.25,.35)
\drawline[AHnb=0](9.65,.4)(9.5,.05)

\drawline[AHnb=0](9.65,.7)(9.5,1)
\drawline[AHnb=0](9.65,.7)(9.25,.65)

\drawline[AHnb=0](9.8,1.15)(9.72,1.5)
\drawline[AHnb=0](10.2,1.15)(10.3,1.5)

\drawline[AHnb=0](9.5,1)(9.3,1.2)
\drawline[AHnb=0](9.25,.65)(9,.75)

\drawline[AHnb=0](9.5,.05)(9.3,-.2)
\drawline[AHnb=0](9.25,.35)(9,.25)

\drawline[AHnb=0](9.8,-.1)(9.7,-.5)
\drawline[AHnb=0](10.2,-.1)(10.3,-.5)

\drawline[AHnb=0](10.5,.05)(10.7,-.2)
\drawline[AHnb=0](10.75,.35)(11,.25)

\drawline[AHnb=0](10.75,.65)(11,.75)
\drawline[AHnb=0](10.5,1)(10.7,1.2)

\drawline[AHnb=0](1,1.7)(1,2.3)
\drawline[AHnb=0](3,1.7)(3,2.3)
\drawline[AHnb=0](5,1.7)(5,2.3)
\drawline[AHnb=0](7,1.7)(7,2.3)
\drawline[AHnb=0](10,1.7)(10,2.3)

\drawline[AHnb=0](.7,1.5)(1,2)
\drawline[AHnb=0](2.7,1.5)(3,2)
\drawline[AHnb=0](4.7,1.5)(5,2)
\drawline[AHnb=0](6.7,1.5)(7,2)
\drawline[AHnb=0](9.7,1.5)(10,2)

\drawline[AHnb=0](1.3,1.5)(1,2)
\drawline[AHnb=0](3.3,1.5)(3,2)
\drawline[AHnb=0](5.3,1.5)(5,2)
\drawline[AHnb=0](7.3,1.5)(7,2)
\drawline[AHnb=0](10.3,1.5)(10,2)

\drawline[AHnb=0](1,-.7)(1,-1)
\drawline[AHnb=0](3,-.7)(3,-1)
\drawline[AHnb=0](5,-.7)(5,-1)
\drawline[AHnb=0](7,-.7)(7,-1)
\drawline[AHnb=0](10,-.7)(10,-1)

%
%
%

\put(.18,-.34){\scriptsize {\tiny $b_1$}}
\put(2.18,-.34){\scriptsize {\tiny $b_3$}}
\put(4.18,-.34){\scriptsize {\tiny $b_5$}}
\put(6.18,-.34){\scriptsize {\tiny $b_7$}}
\put(8.8,-.34){\scriptsize {\tiny $b_{n-1}$}}
\put(11.18,-.34){\scriptsize {\tiny $b_1$}}

\put(1.3,-.6){\scriptsize {\tiny $b_2$}}
\put(3.3,-.6){\scriptsize {\tiny $b_4$}}
\put(5.3,-.6){\scriptsize {\tiny $b_6$}}
\put(7.3,-.6){\scriptsize {\tiny $b_8$}}
\put(10.3,-.6){\scriptsize {\tiny $b_{2n}$}}

\put(.55,-.6){\scriptsize {\tiny $c_2$}}
\put(2.55,-.6){\scriptsize {\tiny $c_4$}}
\put(4.55,-.6){\scriptsize {\tiny $c_6$}}
\put(6.55,-.6){\scriptsize {\tiny $c_8$}}
\put(9.45,-.6){\scriptsize {\tiny $c_{2n}$}}

\put(1.7,-.34){\scriptsize {\tiny $c_3$}}
\put(3.7,-.34){\scriptsize {\tiny $c_5$}}
\put(5.7,-.34){\scriptsize {\tiny $c_7$}}
\put(7.7,-.34){\scriptsize {\tiny $c_9$}}
\put(10.7,-.34){\scriptsize {\tiny $c_1$}}

\put(-.22,.2){\scriptsize {\tiny $a_1$}}
\put(11.1,.2){\scriptsize {\tiny $a_1$}}

\put(-.65,.75){\scriptsize {\tiny $a_{2n-4}$}}
\put(11.1,.7){\scriptsize {\tiny $a_{2n-4}$}}

\put(1.3,1.54){\scriptsize {\tiny $b_{2n-3}$}}
\put(3.3,1.54){\scriptsize {\tiny $b_{2n-1}$}}
\put(5.3,1.54){\scriptsize {\tiny $b_{1}$}}
\put(7.3,1.54){\scriptsize {\tiny $b_{3}$}}
\put(10.3,1.54){\scriptsize {\tiny $b_{2n-5}$}}

\put(.07,1.54){\scriptsize {\tiny $c_{2n-3}$}}
\put(2.07,1.54){\scriptsize {\tiny $c_{2n-1}$}}
\put(4.47,1.54){\scriptsize {\tiny $c_{1}$}}
\put(6.47,1.54){\scriptsize {\tiny $c_{3}$}}
\put(9.07,1.54){\scriptsize {\tiny $c_{2n-5}$}}

\put(-.35,1.25){\scriptsize {\tiny $b_{2n-4}$}}
\put(2.05,1.3){\scriptsize {\tiny $b_{2n-2}$}}
\put(4.05,1.3){\scriptsize {\tiny $b_{2n}$}}
\put(6.1,1.3){\scriptsize {\tiny $b_{2}$}}
\put(8.7,1.3){\scriptsize {\tiny $b_{2n-6}$}}
\put(11.35,1.1){\scriptsize {\tiny $b_{2n-4}$}}

\put(1.05,2){\scriptsize {\tiny $a_{2n-3}$}}
\put(3.05,2){\scriptsize {\tiny $a_{2n-1}$}}
\put(5.05,2){\scriptsize {\tiny $a_1$}}
\put(7.05,2){\scriptsize {\tiny $a_3$}}
\put(10.05,2){\scriptsize {\tiny $a_{2n-5}$}}

\put(1.4,1.3){\scriptsize {\tiny $c_{2n-2}$}}
\put(3.65,1.3){\scriptsize {\tiny $c_{2n}$}}
\put(5.65,1.3){\scriptsize {\tiny $c_{2}$}}
\put(7.65,1.3){\scriptsize {\tiny $c_{4}$}}
\put(10.4,1.3){\scriptsize {\tiny $c_{2n-4}$}}


\put(8.3,.5){\scriptsize $\cdots$}

\drawline[AHnb=0](.65,.4)(.25,.65)
\drawline[AHnb=0](.65,.7)(.8,1.15)
\drawline[AHnb=0](1,.9)(1.5,1)
\drawline[AHnb=0](1.35,.7)(1.75,.35)
\drawline[AHnb=0](1.35,.4)(1.2,-.1)
\drawline[AHnb=0](1,.2)(.5,.05)

\drawline[AHnb=0](2.65,.4)(2.25,.65)
\drawline[AHnb=0](2.65,.7)(2.8,1.15)
\drawline[AHnb=0](3,.9)(3.5,1)
\drawline[AHnb=0](3.35,.7)(3.75,.35)
\drawline[AHnb=0](3.35,.4)(3.2,-.1)
\drawline[AHnb=0](3,.2)(2.5,.05)

\drawline[AHnb=0](4.65,.4)(4.25,.65)
\drawline[AHnb=0](4.65,.7)(4.8,1.15)
\drawline[AHnb=0](5,.9)(5.5,1)
\drawline[AHnb=0](5.35,.7)(5.75,.35)
\drawline[AHnb=0](5.35,.4)(5.2,-.1)
\drawline[AHnb=0](5,.2)(4.5,.05)

\drawline[AHnb=0](6.65,.4)(6.25,.65)
\drawline[AHnb=0](6.65,.7)(6.8,1.15)
\drawline[AHnb=0](7,.9)(7.5,1)
\drawline[AHnb=0](7.35,.7)(7.75,.35)
\drawline[AHnb=0](7.35,.4)(7.2,-.1)
\drawline[AHnb=0](7,.2)(6.5,.05)

\drawline[AHnb=0](1.2,1.15)(.7,1.5)
\drawline[AHnb=0](.5,1)(0,.75)
\drawline[AHnb=0](.25,.35)(.3,-.17)
\drawline[AHnb=0](.8,-.1)(1.3,-.5)
\drawline[AHnb=0](1.5,.05)(2,.25)
\drawline[AHnb=0](1.75,.65)(1.7,1.2)

\drawline[AHnb=0](3.2,1.15)(2.7,1.5)
\drawline[AHnb=0](2.5,1)(2,.75)
\drawline[AHnb=0](2.25,.35)(2.3,-.17)
\drawline[AHnb=0](2.8,-.1)(3.3,-.5)
\drawline[AHnb=0](3.5,.05)(4,.25)
\drawline[AHnb=0](3.75,.65)(3.7,1.2)

\drawline[AHnb=0](5.2,1.15)(4.7,1.5)
\drawline[AHnb=0](4.5,1)(4,.75)
\drawline[AHnb=0](4.25,.35)(4.3,-.17)
\drawline[AHnb=0](4.8,-.1)(5.3,-.5)
\drawline[AHnb=0](5.5,.05)(6,.25)
\drawline[AHnb=0](5.75,.65)(5.7,1.2)

\drawline[AHnb=0](7.2,1.15)(6.7,1.5)
\drawline[AHnb=0](6.5,1)(6,.75)
\drawline[AHnb=0](6.25,.35)(6.3,-.17)
\drawline[AHnb=0](6.8,-.1)(7.3,-.5)
\drawline[AHnb=0](7.5,.05)(8,.25)
\drawline[AHnb=0](7.75,.65)(7.7,1.2)

\drawline[AHnb=0](10.2,1.15)(9.7,1.5)
\drawline[AHnb=0](9.5,1)(9,.75)
\drawline[AHnb=0](9.25,.35)(9.3,-.17)
\drawline[AHnb=0](9.8,-.1)(10.3,-.5)
\drawline[AHnb=0](10.5,.05)(11,.25)
\drawline[AHnb=0](10.75,.65)(10.7,1.2)
\drawline[AHnb=0](9.65,.4)(9.25,.65)
\drawline[AHnb=0](9.65,.7)(9.8,1.15)
\drawline[AHnb=0](10,.9)(10.5,1)
\drawline[AHnb=0](10.35,.7)(10.75,.35)
\drawline[AHnb=0](10.35,.4)(10.2,-.1)
\drawline[AHnb=0](10,.2)(9.5,.05)

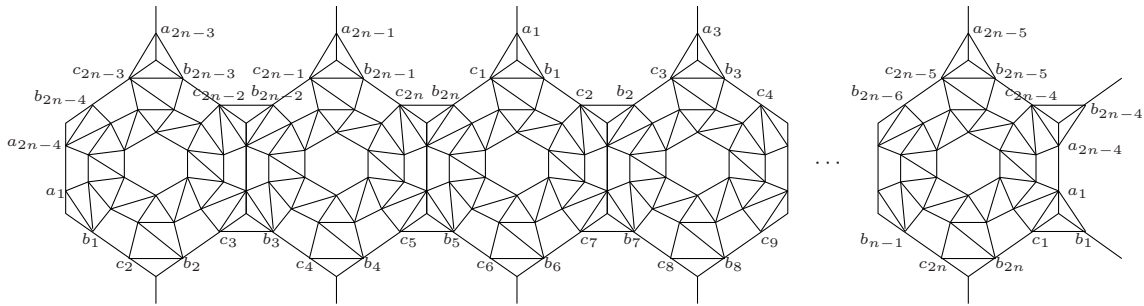
\put(.5,-1.5){\scriptsize{\bf Fig.(4.12)\,: SEM of type\,-$(3^4, 6)$, $n\geq 7$}}

\end{picture}

\end{center}

\vspace{13.5cm}

\hrule

The torus is a double cover of the Klein bottle therefore each map $M$ on the surface of Klein bottle can be lifted to the
surface of the torus by a suitable double cover map. The following SEMs $T_{1, 28}(3^3, 4^2)$, $T_{2,24}(3^3, 4^2)$,
$T_{24}(3^2, 4, 3, 4)$, $T_{1,20}(3^3, 4^2)$, $T_{1,24}(3^3, 4^2)$ and $T_{1, 36}(3, 4, 6, 4)$ are the double covers of the
SEMs $K_{1, 14}(3^3, 4^2)$, $K_{2, 12}(3^3, 4^2)$, $K_{12}(3^3, 4, 3, 4)$, $K_{1,10}(3^3, 4^2)$, $K_{1, 12}(3^3, 4^2)$ and
$K_{1, 18}(3, 6, 4, 6)$ respectively. This can be seen easily by considering double covering maps\,: $\theta\,: T_{1, 28}(3^3, 4^2) \,
\longrightarrow \, K_{1, 14}(3^3, 4^2)$, $\gamma\,: T_{2,24}(3^3, 4^2) \, \longrightarrow \, K_{2,12}(3^3, 4^2)$,
$\phi\,: T_{24}(3^2, 4, 3, 4) \, \longrightarrow \, K_{12}(3^2, 4, 3, 4)$, $\alpha\,: T_{1,20}(3^3, 4^2) \, \longrightarrow \, K_{1,10}(3^3, 4^2)$,
$\beta\,: T_{1,24}(3^3, 4^2) \, \longrightarrow \, K_{1,12}(3^3, 4^2)$ and $\varphi\,: T_{36}(3, 4, 6, 4) \, \longrightarrow \, K_{18}(3, 4, 6, 4)$
such that $\theta\{i, \, i+14\} = v_i$ (for $0\leq i \leq 13$), $\gamma\{i, \, i+12\} = v_i$ (for $0\leq i \leq 11$),
$\phi\{i, \, i+12\} = v_i$ (for $0\leq i \leq 11$), $\alpha\{i, \, i+10\} = v_i$ (for $0\leq i \leq 9$),
$\beta\{i, \, i+12\} = v_i$ (for $0\leq i \leq 11$)  and $\varphi\{i, \, i+18\} = v_i$ (for $0\leq i \leq 17$).

\smallskip

\hrule

\begin{center}

\begin{picture}(0,0)(75,38) 
\setlength{\unitlength}{8.5mm}

\drawpolygon(0,0)(7,0)(7,4)(0,4)

\drawline[AHnb=0](0,1)(7,1)
\drawline[AHnb=0](0,2)(7,2)
\drawline[AHnb=0](0,3)(7,3)

\drawline[AHnb=0](1,0)(1,4)
\drawline[AHnb=0](2,0)(2,4)
\drawline[AHnb=0](3,0)(3,4)
\drawline[AHnb=0](4,0)(4,4)
\drawline[AHnb=0](5,0)(5,4)
\drawline[AHnb=0](6,0)(6,4)

\drawline[AHnb=0](0,1)(1,2)
\drawline[AHnb=0](1,1)(2,2)
\drawline[AHnb=0](2,1)(3,2)
\drawline[AHnb=0](3,1)(4,2)
\drawline[AHnb=0](4,1)(5,2)
\drawline[AHnb=0](5,1)(6,2)
\drawline[AHnb=0](6,1)(7,2)

\drawline[AHnb=0](1,3)(0,4)
\drawline[AHnb=0](2,3)(1,4)
\drawline[AHnb=0](3,3)(2,4)
\drawline[AHnb=0](4,3)(3,4)
\drawline[AHnb=0](5,3)(4,4)
\drawline[AHnb=0](6,3)(5,4)
\drawline[AHnb=0](7,3)(6,4)

\put(-.2,-.35){\scriptsize 19}
\put(.8,-.35){\scriptsize 18}
\put(1.8,-.35){\scriptsize 13}
\put(2.95,-.35){\scriptsize 6}
\put(3.95,-.35){\scriptsize 0}
\put(4.8,-.35){\scriptsize 17}
\put(5.8,-.35){\scriptsize 26}
\put(6.9,-.35){\scriptsize 19}

\put(-.2,4.15){\scriptsize 19}
\put(.8,4.15){\scriptsize 18}
\put(1.8,4.15){\scriptsize 13}
\put(2.95,4.15){\scriptsize 6}
\put(3.95,4.15){\scriptsize 0}
\put(4.8,4.15){\scriptsize 17}
\put(5.8,4.15){\scriptsize 26}
\put(6.9,4.15){\scriptsize 19}

\put(-.43,.9){\scriptsize 22}
\put(-.43,1.9){\scriptsize 21}
\put(-.43,2.9){\scriptsize 20}

\put(7.1,.9){\scriptsize 22}
\put(7.1,1.9){\scriptsize 21}
\put(7.1,2.9){\scriptsize 20}

\put(.55,.7){\scriptsize 23}
\put(1.55,.7){\scriptsize 10}
\put(2.65,.7){\scriptsize 7}
\put(3.65,.7){\scriptsize 1}
\put(4.55,.7){\scriptsize 16}
\put(5.55,.7){\scriptsize 25}

\put(.55,2.7){\scriptsize 14}
\put(1.55,2.7){\scriptsize 3}
\put(2.55,2.7){\scriptsize 12}
\put(3.7,2.7){\scriptsize 5}
\put(4.7,2.7){\scriptsize 4}
\put(5.55,2.7){\scriptsize 27}

\put(1.1,2.1){\scriptsize 15}
\put(2.1,2.1){\scriptsize 2}
\put(3.1,2.1){\scriptsize 11}
\put(4.1,2.1){\scriptsize 8}
\put(5.1,2.1){\scriptsize 9}
\put(6.1,2.1){\scriptsize 24}

\put(2.5,-1){\scriptsize $T_{1,28}(3^3, 4^2)$}

\end{picture}

\end{center}


\begin{center}

\begin{picture}(0,0)(-2,28) 
\setlength{\unitlength}{9mm}

\drawpolygon(0,0)(8,0)(8,3)(0,3)

\drawline[AHnb=0](0,1)(8,1)
\drawline[AHnb=0](0,2)(8,2)

\drawline[AHnb=0](1,0)(1,3)
\drawline[AHnb=0](2,0)(2,3)
\drawline[AHnb=0](3,0)(3,3)
\drawline[AHnb=0](4,0)(4,3)
\drawline[AHnb=0](5,0)(5,3)
\drawline[AHnb=0](6,0)(6,3)
\drawline[AHnb=0](7,0)(7,3)

\drawline[AHnb=0](0,1)(1,0)
\drawline[AHnb=0](0,2)(1,1)
\drawline[AHnb=0](0,3)(1,2)

\drawline[AHnb=0](2,1)(3,0)
\drawline[AHnb=0](2,2)(3,1)
\drawline[AHnb=0](2,3)(3,2)

\drawline[AHnb=0](4,0)(5,1)
\drawline[AHnb=0](4,1)(5,2)
\drawline[AHnb=0](4,2)(5,3)

\drawline[AHnb=0](6,0)(7,1)
\drawline[AHnb=0](6,1)(7,2)
\drawline[AHnb=0](6,2)(7,3)

\put(-.2,-.3){\scriptsize 11}
\put(.9,-.3){\scriptsize 6}
\put(1.9,-.3){\scriptsize 7}
\put(2.9,-.3){\scriptsize 8}
\put(3.8,-.3){\scriptsize 16}
\put(4.8,-.3){\scriptsize 15}
\put(5.8,-.3){\scriptsize 14}
\put(6.8,-.3){\scriptsize 22}
\put(7.9,-.3){\scriptsize 11}

\put(-.2,3.1){\scriptsize 11}
\put(.9,3.1){\scriptsize 6}
\put(1.9,3.1){\scriptsize 7}
\put(2.9,3.1){\scriptsize 8}
\put(3.8,3.1){\scriptsize 16}
\put(4.8,3.1){\scriptsize 15}
\put(5.8,3.1){\scriptsize 14}
\put(6.8,3.1){\scriptsize 22}
\put(7.9,3.1){\scriptsize 11}

\put(-.3,1.9){\scriptsize 4}
\put(-.3,.9){\scriptsize 5}

\put(8.1,1.9){\scriptsize 4}
\put(8.1,.9){\scriptsize 5}

\put(1.1,2.1){\scriptsize 3}
\put(1.1,1.1){\scriptsize 0}

\put(3.1,2.1){\scriptsize 10}
\put(3.1,1.1){\scriptsize 9}

\put(5.1,2.1){\scriptsize 18}
\put(5.1,1.1){\scriptsize 12}

\put(7.1,2.1){\scriptsize 20}
\put(7.1,1.1){\scriptsize 21}

\put(1.75,1.7){\scriptsize 2}
\put(1.75,.7){\scriptsize 1}

\put(3.6,1.7){\scriptsize 23}
\put(3.6,.7){\scriptsize 17}

\put(5.6,1.7){\scriptsize 19}
\put(5.6,.7){\scriptsize 13}

\put(3.5,-1){\scriptsize $T_{2, 24}(3^3, 4^2)$}

\end{picture}
\end{center}


\begin{center}

\begin{picture}(0,0)(77,75) 
\setlength{\unitlength}{7.8mm}

\drawpolygon(0,0)(6,0)(6,4)(0,4)

\drawline[AHnb=0](0,1)(6,1)
\drawline[AHnb=0](0,2)(6,2)
\drawline[AHnb=0](0,3)(6,3)

\drawline[AHnb=0](1,0)(1,4)
\drawline[AHnb=0](2,0)(2,4)
\drawline[AHnb=0](3,0)(3,4)
\drawline[AHnb=0](4,0)(4,4)
\drawline[AHnb=0](5,0)(5,4)

\drawline[AHnb=0](2,0)(1,1)
\drawline[AHnb=0](3,1)(4,0)
\drawline[AHnb=0](5,1)(6,0)
\drawline[AHnb=0](0,1)(1,2)
\drawline[AHnb=0](2,1)(3,2)
\drawline[AHnb=0](4,1)(5,2)
\drawline[AHnb=0](1,3)(2,2)
\drawline[AHnb=0](3,3)(4,2)
\drawline[AHnb=0](5,3)(6,2)
\drawline[AHnb=0](0,3)(1,4)
\drawline[AHnb=0](2,3)(3,4)
\drawline[AHnb=0](4,3)(5,4)

\put(-.2,-.33){\scriptsize 6}
\put(.9,-.33){\scriptsize 7}
\put(1.9,-.33){\scriptsize 1}
\put(2.8,-.33){\scriptsize 18}
\put(3.8,-.33){\scriptsize 19}
\put(4.8,-.33){\scriptsize 13}
\put(6.1,-.33){\scriptsize 6}

\put(-.2,4.1){\scriptsize 6}
\put(.9,4.1){\scriptsize 7}
\put(1.9,4.1){\scriptsize 1}
\put(2.8,4.1){\scriptsize 18}
\put(3.8,4.1){\scriptsize 19}
\put(4.8,4.1){\scriptsize 13}
\put(6,4.1){\scriptsize 6}

\put(-.3,.9){\scriptsize 5}
\put(-.3,1.9){\scriptsize 8}
\put(-.3,2.9){\scriptsize 9}

\put(6.1,.9){\scriptsize 5}
\put(6.1,1.9){\scriptsize 8}
\put(6.1,2.9){\scriptsize 9}

\put(.6,2.65){\scriptsize 11}
\put(.7,.65){\scriptsize 0}

\put(2.55,2.65){\scriptsize 17}
\put(2.55,.65){\scriptsize 21}

\put(4.55,2.65){\scriptsize 14}
\put(4.55,.65){\scriptsize 22}

\put(1.1,1.7){\scriptsize 4}
\put(3.1,1.7){\scriptsize 20}
\put(5.1,1.7){\scriptsize 15}

\put(1.7,1.1){\scriptsize 2}
\put(3.55,1.1){\scriptsize 23}

\put(1.55,3.1){\scriptsize 10}
\put(3.55,3.1){\scriptsize 12}

\put(2.1,2.1){\scriptsize 3}
\put(4.1,2.1){\scriptsize 16}

\put(1.7,-1){\scriptsize $T_{24}(3^2, 4, 3, 4)$}

\end{picture}

\end{center}


\begin{center}

\begin{picture}(0,0)(21,67) 
\setlength{\unitlength}{7.8mm}

\drawpolygon(0,0)(5,0)(5,4)(0,4)
\drawline[AHnb=0](0,1)(5,1)
\drawline[AHnb=0](0,2)(5,2)
\drawline[AHnb=0](0,3)(5,3)

\drawline[AHnb=0](1,0)(1,4)
\drawline[AHnb=0](2,0)(2,4)
\drawline[AHnb=0](3,0)(3,4)
\drawline[AHnb=0](4,0)(4,4)

\drawline[AHnb=0](0,0)(1,1)
\drawline[AHnb=0](1,0)(2,1)
\drawline[AHnb=0](2,0)(3,1)
\drawline[AHnb=0](3,0)(4,1)
\drawline[AHnb=0](4,0)(5,1)

\drawline[AHnb=0](0,3)(1,2)
\drawline[AHnb=0](1,3)(2,2)
\drawline[AHnb=0](2,3)(3,2)
\drawline[AHnb=0](3,3)(4,2)
\drawline[AHnb=0](4,3)(5,2)

\put(-.3,-.3){\scriptsize 18}
\put(-.45,.9){\scriptsize 17}
\put(-.45,1.9){\scriptsize 16}
\put(-.45,2.9){\scriptsize 15}
\put(-.3,4.1){\scriptsize 18}

\put(5.05,-.3){\scriptsize 18}
\put(5.05,.9){\scriptsize 17}
\put(5.05,1.9){\scriptsize 16}
\put(5.05,2.9){\scriptsize 15}
\put(5.05,4.1){\scriptsize 18}

\put(.9,-.3){\scriptsize 9}
\put(1.9,-.3){\scriptsize 7}
\put(2.9,-.3){\scriptsize 11}
\put(3.9,-.3){\scriptsize 12}

\put(.9,4.1){\scriptsize 9}
\put(1.9,4.1){\scriptsize 7}
\put(2.9,4.1){\scriptsize 11}
\put(3.9,4.1){\scriptsize 12}

\put(.75,3.1){\scriptsize 4}
\put(1.75,3.1){\scriptsize 6}
\put(2.56,3.1){\scriptsize 10}
\put(3.56,3.1){\scriptsize 13}

\put(.75,1.7){\scriptsize 0}
\put(1.75,1.7){\scriptsize 3}
\put(2.75,1.7){\scriptsize 5}
\put(3.54,1.7){\scriptsize 14}

\put(1.05,1.1){\scriptsize 1}
\put(2.05,1.1){\scriptsize 2}
\put(3.05,1.1){\scriptsize 8}
\put(4.05,1.1){\scriptsize 19}

\put(1.5,-1){\scriptsize $T_{1, 20}(3^3, 4^2)$}

\end{picture}

\end{center}


\begin{center}

\begin{picture}(0,0)(-28,58) 
\setlength{\unitlength}{7.8mm}

\drawpolygon(0,0)(6,0)(6,4)(0,4)
\drawline[AHnb=0](0,1)(6,1)
\drawline[AHnb=0](0,2)(6,2)
\drawline[AHnb=0](0,3)(6,3)

\drawline[AHnb=0](1,0)(1,4)
\drawline[AHnb=0](2,0)(2,4)
\drawline[AHnb=0](3,0)(3,4)
\drawline[AHnb=0](4,0)(4,4)
\drawline[AHnb=0](5,0)(5,4)

\drawline[AHnb=0](0,0)(1,1)
\drawline[AHnb=0](1,0)(2,1)
\drawline[AHnb=0](2,0)(3,1)
\drawline[AHnb=0](3,0)(4,1)
\drawline[AHnb=0](4,0)(5,1)
\drawline[AHnb=0](5,0)(6,1)

\drawline[AHnb=0](0,3)(1,2)
\drawline[AHnb=0](1,3)(2,2)
\drawline[AHnb=0](2,3)(3,2)
\drawline[AHnb=0](3,3)(4,2)
\drawline[AHnb=0](4,3)(5,2)
\drawline[AHnb=0](5,3)(5,4)

\put(-.3,-.3){\scriptsize 23}
\put(-.45,.9){\scriptsize 19}
\put(-.45,1.9){\scriptsize 18}
\put(-.45,2.9){\scriptsize 22}
\put(-.3,4.1){\scriptsize 23}

\put(6.05,-.3){\scriptsize 23}
\put(6.05,.9){\scriptsize 19}
\put(6.05,1.9){\scriptsize 18}
\put(6.05,2.9){\scriptsize 22}
\put(6.05,4.1){\scriptsize 23}

\put(.9,-.3){\scriptsize 4}
\put(1.9,-.3){\scriptsize 5}
\put(2.9,-.3){\scriptsize 10}
\put(3.9,-.3){\scriptsize 21}
\put(4.9,-.3){\scriptsize 20}

\put(.9,4.1){\scriptsize 4}
\put(1.9,4.1){\scriptsize 5}
\put(2.9,4.1){\scriptsize 10}
\put(3.9,4.1){\scriptsize 21}
\put(4.9,4.1){\scriptsize 20}

\put(.75,3.1){\scriptsize 9}
\put(1.75,3.1){\scriptsize 8}
\put(2.56,3.1){\scriptsize 11}
\put(3.56,3.1){\scriptsize 16}
\put(4.56,3.1){\scriptsize 17}

\put(.75,1.7){\scriptsize 3}
\put(1.75,1.7){\scriptsize 0}
\put(2.75,1.7){\scriptsize 6}
\put(3.56,1.7){\scriptsize 14}
\put(4.56,1.7){\scriptsize 13}

\put(1.05,1.1){\scriptsize 2}
\put(2.05,1.1){\scriptsize 1}
\put(3.05,1.1){\scriptsize 7}
\put(4.05,1.1){\scriptsize 15}
\put(5.05,1.1){\scriptsize 12}

\put(1.8,-1){\scriptsize $T_{1, 24}(3^3, 4^2)$}

\end{picture}

\end{center}


\begin{center}

\begin{picture}(0,0)(70,90) 
\setlength{\unitlength}{8mm}

\drawpolygon(0,0)(1,-.5)(2,0)(3,0)(4,-.5)(5,0)(6,0)(7,-.5)(8,0)(9,0)(10,-.5)(11,0)(12,0)(13,-.5)(14,0)(15,0)(16,-.5)(17,0)
(17,1)(16,1.5)(15,1)(14,1)(13,1.5)(12,1)(11,1)(10,1.5)(9,1)(8,1)(7,1.5)(6,1)(5,1)(4,1.5)(3,1)(2,1)(1,1.5)(0,1)

\drawline[AHnb=0](2,0)(2,1)
\drawline[AHnb=0](3,0)(3,1)

\drawline[AHnb=0](5,0)(5,1)
\drawline[AHnb=0](6,0)(6,1)

\drawline[AHnb=0](8,0)(8,1)
\drawline[AHnb=0](9,0)(9,1)

\drawline[AHnb=0](11,0)(11,1)
\drawline[AHnb=0](12,0)(12,1)

\drawline[AHnb=0](14,0)(14,1)
\drawline[AHnb=0](15,0)(15,1)

\drawpolygon(1,1.5)(.6,2.2)(1.4,2.2)
\drawpolygon(4,1.5)(3.6,2.2)(4.4,2.2)
\drawpolygon(7,1.5)(6.6,2.2)(7.4,2.2)
\drawpolygon(10,1.5)(9.6,2.2)(10.4,2.2)
\drawpolygon(13,1.5)(12.6,2.2)(13.4,2.2)
\drawpolygon(16,1.5)(15.6,2.2)(16.4,2.2)

\drawpolygon(2,1)(3,1)(2.5,1.7)
\drawpolygon(5,1)(6,1)(5.5,1.7)
\drawpolygon(8,1)(9,1)(8.5,1.7)
\drawpolygon(11,1)(12,1)(11.5,1.7)
\drawpolygon(14,1)(15,1)(14.5,1.7)
\drawpolygon(17,1)(18,1)(17.5,1.7)

\drawline[AHnb=0](1.4,2.2)(2.5,1.7)
\drawline[AHnb=0](4.4,2.2)(5.5,1.7)
\drawline[AHnb=0](7.4,2.2)(8.5,1.7)
\drawline[AHnb=0](10.4,2.2)(11.5,1.7)
\drawline[AHnb=0](13.4,2.2)(14.5,1.7)
\drawline[AHnb=0](16.4,2.2)(17.5,1.7)

\drawline[AHnb=0](-.5,1.7)(.6,2.2)
\drawline[AHnb=0](2.5,1.7)(3.6,2.2)
\drawline[AHnb=0](5.5,1.7)(6.6,2.2)
\drawline[AHnb=0](8.5,1.7)(9.6,2.2)
\drawline[AHnb=0](11.5,1.7)(12.6,2.2)
\drawline[AHnb=0](14.5,1.7)(15.6,2.2)

\drawline[AHnb=0](-.5,1.7)(0,1)

\drawline[AHnb=0](17,0)(18,0)
\drawline[AHnb=0](18,0)(18,1)

\put(.83,-.9){\scriptsize 11}
\put(3.9,-.9){\scriptsize 8}
\put(6.9,-.9){\scriptsize 4}
\put(9.83,-.9){\scriptsize 28}
\put(12.83,-.9){\scriptsize 31}
\put(15.9,-.9){\scriptsize 23}

\put(-.2,-.4){\scriptsize 15}
\put(1.8,-.4){\scriptsize 10}
\put(4.8,-.4){\scriptsize 13}
\put(7.9,-.4){\scriptsize 5}
\put(10.8,-.4){\scriptsize 29}
\put(13.8,-.4){\scriptsize 26}
\put(16.8,-.4){\scriptsize 22}

\put(2.9,-.4){\scriptsize 9}
\put(5.9,-.4){\scriptsize 3}
\put(8.8,-.4){\scriptsize 24}
\put(11.8,-.4){\scriptsize 30}
\put(14.8,-.4){\scriptsize 18}
\put(17.9,-.4){\scriptsize 15}

\put(1.7,.7){\scriptsize 6}
\put(4.55,.7){\scriptsize 12}
\put(7.7,.7){\scriptsize 0}
\put(10.55,.7){\scriptsize 33}
\put(13.55,.7){\scriptsize 27}
\put(16.55,.7){\scriptsize 21}

\put(.1,.7){\scriptsize 14}
\put(3.1,.7){\scriptsize 17}
\put(6.1,.7){\scriptsize 2}
\put(9.1,.7){\scriptsize 25}
\put(12.1,.7){\scriptsize 34}
\put(15.1,.7){\scriptsize 19}
\put(18.1,.7){\scriptsize 14}

\put(.9,1.1){\scriptsize 7}
\put(3.8,1.1){\scriptsize 16}
\put(6.9,1.1){\scriptsize 1}
\put(9.8,1.1){\scriptsize 32}
\put(12.8,1.1){\scriptsize 35}
\put(15.8,1.1){\scriptsize 20}

\put(-.7,1.9){\scriptsize 31}
\put(2.35,1.9){\scriptsize 23}
\put(5.35,1.9){\scriptsize 11}
\put(8.4,1.9){\scriptsize 8}
\put(11.4,1.9){\scriptsize 4}
\put(14.35,1.9){\scriptsize 28}
\put(17.4,1.9){\scriptsize 31}

\put(.4,2.3){\scriptsize 26}
\put(3.4,2.3){\scriptsize 22}
\put(6.4,2.3){\scriptsize 10}
\put(9.4,2.3){\scriptsize 13}
\put(12.5,2.3){\scriptsize 5}
\put(15.4,2.3){\scriptsize 29}

\put(1.2,2.3){\scriptsize 18}
\put(4.2,2.3){\scriptsize 15}
\put(7.3,2.3){\scriptsize 9}
\put(10.3,2.3){\scriptsize 3}
\put(13.2,2.3){\scriptsize 24}
\put(16.2,2.3){\scriptsize 30}

\put(7.5,-1.5){\scriptsize $T_{1, 36}(3,4,6,4)$}

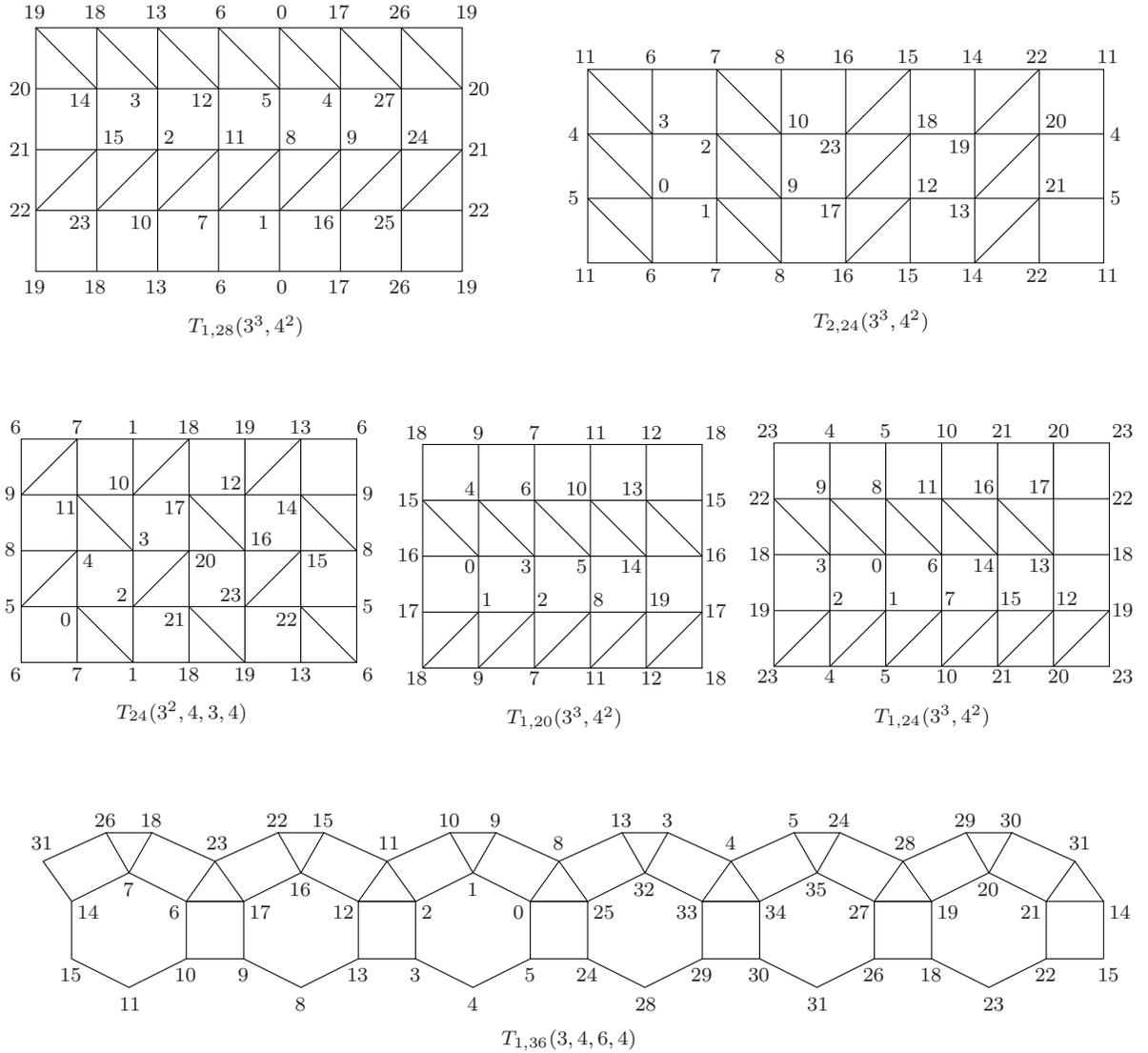
\put(-1.3,-2.5){\scriptsize {\bf Fig.(4.13)\,: Double covers of\,: $K_{1, 14}(3^3, 4^2)$, $K_{2, 12}(3^3, 4^2)$,
$K_{12}(3^2,4,3,4)$, $K_{1, 10}(3^3, 4^2)$, $K_{1, 12}(3^3, 4^2)$, $K_{18}(3, 4, 6, 4)$}}

\end{picture}
\end{center}

\vspace{11cm}

\hrule

\newpage

\begin{rem} \label{remark}
Although we have employed the characteritic polynomial to check whether two maps are isomorphic or not, one can employ a more geometric
technique in the following way: we can assign to each loop starting at a vertex and the number of crossings (crossing edges only
transversally). Further more we can associate to each homotopy class the minimum number of crossings that are needed to represent the homotopy class by such a loop. Similarly we obtain a minimal basis, that is to say, a basis by the set of two elements $\{a, b\}$ with the minimum number of crossings.

Let $A$ and $B$ denote the set which contain the loops the type $a$ and $b$ respectively, then for $T_{1, 12}(3^3, 4^2)$ if we choose any
base point, say $v_7$, then $A = \{C_4(v_7, v_{10}, v_5, v_6)\}$ and $B = \{C_{5}(v_7, v_{11}, v_4, v_5, v_6)$,
$C_5(v_7, v_{11}, v_4, v_0, v_6)$, $C_5(v_7, v_{11}, v_4, v_0, v_1)$, $C_5(v_7, v_{11}, v_3, v_0, v_6)$, $C_5(v_7, v_{11}, v_3, v_0, v_1)$,
$C_5(v_7, v_{11}, v_3, v_2, v_1)$, $C_5(v_7, v_{11}, v_4, v_5, v_{10})$, $C_5(v_7, v_{11}, v_4, v_8, v_{10})$, $C_5(v_7, v_{11}, v_4, v_8, v_1)$, $C_5(v_7, v_{11}, v_9, v_8, v_{10})$, $C_5(v_7, v_{11}, v_9, v_8, v_1)$, $C_5(v_7, v_{11}, v_9, v_2, v_1)$, $C_5(v_7, v_{10}, v_3, v_0, v_6)$,
$C_5(v_7, v_{10}, v_3, v_0, v_1)$, $C_5(v_7, v_{10}, v_3, v_2, v_1)$, $C_5(v_7, v_{10}, v_5, v_2, v_1)$, $C_5(v_7, v_6, v_9, v_8, v_{10})$,
$C_5(v_7, v_6, v_9, v_8, v_1)$, $C_5(v_7, v_6, v_9, v_2, v_1)$, $C_5(v_7, v_6, v_5, v_2, v_1)\}$. Now we see that any combination of the
set $\{a, b\}$ has 3 ($a$ at $v_{10}$, $v_5$ and $v_6$) and 4 crossings ($b$ at 4 vertices of the loop other than $v_7$).
Similarly, for $T_{2, 12}$ the minimal basis $\{a, b\}$ has 2 and 4 crossings (considering, $v_8$ as a base point and
then $a = C_3(v_8, v_9, v_{10})$, $b = C_6(v_8, v_4, v_3, v_2, v_{10})$). For $T_{3, 12}(3^3, 4^2)$ the minimal basis $\{a, b\}$ has 2
and 4 crossings (considering, $v_8$ as a base point and then $a = C_3(v_8, v_9, v_{10})$, $b = C_5(v_8, v_9, v_1, v_0, v_5)$).
Thus, $T_{1, 12}(3^3, 4^2) \ncong T_{2, 12}(3^3, 4^2)$ and $T_{1, 12}(3^3, 4^2) \ncong T_{3, 12}(3^3, 4^2)$.

In case of $T_{1, 14}(3^3, 4^2)$ we see a minimal loop (an element of the minimal basis $\{a, b\}$) having 3 crossings (considering $v_5$ as a base point
we have a minimal loop $C_4(v_5, v_{12}, v_{11}, v_{10})$), while for $T_{2, 14}(3^3, 4^2)$ the minimal loop (an element of the minimal basis
$\{a, b\}$), has 4 crossings (considering $v_5$ as a base point we have a minimal loop $C_5(v_5, v_{11}, v_{10}, v_9, v_8)$). Thus,
$T_{1, 14}(3^3, 4^2) \ncong T_{2, 14}(3^3, 4^2)$.

Similarly, for $K_{1, 12}(3^3, 4^2)$ the minimal loop has 3 crossings (considering $v_{10}$ as a base point we have a minimal loop
$C_4(v_{10}, v_{11}, v_{8}, v_{9})$), while for $T_{2, 14}(3^3, 4^2)$ the minimal loop has only 2 crossings (considering $v_{11}$ as a
base point we have $a = C_5(v_{11}, v_{4}, v_{5})$). Thus, $K_{1, 12}(3^3, 4^2) \ncong K_{2, 14}(3^3, 4^2)$.

%
%

\end{rem}

\section{Acknowledgements} We are thankful to the anonymous referee(s) for her/his(their) suggestions and comments which enlightened and educated us on the subject and presentation of the work. This lead to significant improvements in the content and presentation of the article. Geometric arguments in Remark \ref{remark} were suggested by the referee. Part of this work was completed when AKU was visiting Indian Institute of Science. We would like to thank Prof. B. Datta for his support. AKU also acknowledges financial support from DST-SERB grant number $SR/S4/MS:717/10$ and AKT acknowledges support from CSIR(SRF) award no. 09/1023(0003)/2010-EMR-I.

\end{document}